\documentclass[12pt]{article}
\usepackage{amsmath,amssymb,amsthm, amscd, pb-diagram}

\newtheorem{Prop}[subsection]{Proposition}
\newtheorem{Lem}[subsection]{Lemma}
\newtheorem{Exam}[subsection]{Example} 
\newtheorem{Exams}[subsection]{Examples}

\newtheorem{Cor}[subsection]{Corollary}

\newtheorem{Claim}[subsubsection]{Claim}

\begin{document}

\begin{center}
{\bf Birational geometry for the covering of a nilpotent orbit closure}
\end{center}
\vspace{0.4cm}

\begin{center}
{\bf Yoshinori Namikawa}
\end{center}
\vspace{0.2cm}

\begin{abstract}
Let $\mathfrak{g}$ be a complex classical simple Lie algebra and let $O$ be a nilpotent orbit of $\mathfrak{g}$. 
The fundamental group $\pi_1(O)$ is finite. Take the universal covering $\pi^0: X^0 \to O$. Then $\pi^0$ 
extends to a finite cover $\pi: X \to \bar{O}$. By using the Kirillov-Kostant form $\omega_{KK}$ on $O$, 
the normal affine variety $X$ becomes a conical symplectic variety. In this article we give an explicit construction of 
a {\bf Q}-factorial terminalization of $X$. \vspace{0.3cm}

\end{abstract}

\begin{center}
{\bf Introduction}
\end{center} 

Let $\mathfrak{g}$ be a complex semisimple Lie algebra and let 
$O$ be a nilpotent orbit of $\mathfrak{g}$. Then $O$ admits a 
symplectic 2-form $\omega_{KK}$ called the Kirillov-Kostant 2-form. 
In general $O$ is not simply connected, but $\pi_1(O)$ is finite.  
Let $\pi_0: X^0 \to O$ be a finite etale covering. The function field $\mathbf{C}(X^0)$ of 
$X^0$ is a finite algebraic extension of $\mathbf{C}(O)$. Let $\bar{O}$ be the closure of $O$ in $\mathfrak{g}$. The 
normalization $X$ of $\bar{O}$ in $\mathbf{C}(X^0)$ determines a finite 
covering $\pi: X \to \bar{O}$. Then $X^0$ is contained in $X$ as a Zariski open subset so that $\pi\vert_{X^0} = \pi_0$. If $G$ is a simply connected complex semisimple Lie 
group with $Lie(G) = \mathfrak{g}$, then the adjoint $G$-action on $O$ (resp. $\bar{O}$)  
lifts to a $G$ action on $X^0$ (resp. $X$). The 2-form  
$\omega := (\pi^0)^*\omega_{KK}$ is a $G$-invariant symplectic 2-form on $X^0$. 
The pair $(X, \omega)$ is then a symplectic variety in the sense of Beauville \cite{Be}. 
Brylinski and Kostant \cite{B-K} extensively studied such varieties in the context of 
shared orbits. A nilpotent orbit closure $\bar{O}$ has a natural scaling $\mathbf{C}^*$-action for which $\omega_{KK}$ has weight $1$. Let $s: \mathbf{C}^* 
\to \mathrm{Aut}(\bar{O})$ be a homomorphism determined by the scaling 
action. In general this $\mathbf{C}^*$-action does not lift to a $\mathbf{C}^*$-action on $X$. But, if we instead define a new $\mathbf{C}^*$-action on $\bar{O}$ 
by the composite $\sigma$ of $$\mathbf{C}^* \to \mathbf{C}^* \:\: (t \to t^2), \:\:\: 
\mathrm{and} \:\:\: s: \mathbf{C}^* \to \mathrm{Aut}(\bar{O}),$$
then the new $\mathbf{C}^*$-action $\sigma$ always lifts to a $\mathbf{C}^*$-action on $X$ by 
\cite{B-K}, \S 1. By definition, $wt(\omega) = 2$ with respect to this $\mathbf{C}^*$-action. 
Therefore $(X, \omega)$ is a conical symplectic $G$-variety with $wt(\omega) = 2$.   
   
A main purpose of this article is to study the birational geometry for the resolutions 
of $(X, \omega)$. A crepant projective resolution $f: Y \to X$ of $X$ is, by definition, 
a projective birational morphism $f$ from a nonsingular variety $Y$ to $X$ such 
that $K_Y = f^*K_X$. In general $X$ does not have a crepant projective resolution. 
But, instead, $X$ always has a nice crepant projective partial resolution $f: Y \to X$  called a {\em {\bf Q}-factorial terminalization}  by \cite{BCHM}. A {\bf Q}-factorial 
terminalization $f$ is, by definition, a projective birational morphism from a normal 
variety $Y$ to $X$ such that $Y$ has only {\bf Q}-factorial terminal singularities and $K_Y = f^*K_X$. 
It is natural to expect that such a {\bf Q}-factorial terminalization can be constructed very explicitly in a group theoretic manner when $(X, \omega)$ is the above. 

When $X^0$ is the nilpotent orbit $O$ itself, $X$ is nothing but the normalization $\tilde{O}$ of $\bar{O}$. Namikawa \cite{Na 1} and Fu \cite{Fu 1} respectively constructed 
a {\bf Q}-factorial terminalization of $\tilde{O}$ quite explicitly when $\mathfrak{g}$ is a classical simple Lie algebra and when $\mathfrak{g}$ is an exceptional simple Lie algebra
(See also \cite{Lo} for a unified treatment). 
However, when $\mathrm{deg}(\pi) > 1$, this problem has not yet been enough pursued  
though there is an interesting observation by \cite{Fu 2}, Theorem 1.4 for the case $\mathrm{deg}(\pi) $ is odd.     

In this article we construct a {\bf Q}-factorial terminalization of $X$ when $O$ is a nilpotent orbit of a classical simple Lie algebra $\mathfrak{g}$ and $X^0$ is the 
universal covering of $O$. We shall explain here a basic idea for the construction. 
Let $Q \subset G$ be a parabolic subgroup of $G$ and let $Q = U\cdot L$ be a Levi 
decomposition of $Q$ by the unipotent radical $U$ and a Levi subgroup $L$. 
Correspondingly the Lie algebra $\mathfrak{q}$ decomposes 
$\mathfrak{q} = \mathfrak{n} \oplus \mathfrak{l}$ as a direct sum of  $\mathfrak{n} := Lie(U)$ and $\mathfrak{l} := Lie(L)$. Let $O'$ be a nilpotent orbit of $\mathfrak{l}$. 
Then there is a unique nilpotent orbit $O$ of $\mathfrak{g}$ such that $O$ meets 
$\mathfrak{n} + O'$ in a Zariski open subset of $\mathfrak{n} + O'$. In such a case 
we say that $O$ is induced from $O'$ and write $O = \mathrm{Ind}^{\mathfrak g}_{\mathfrak l}(O')$.  There is a generically finite map 
$$\mu: G \times^Q(\mathfrak{n} + \bar{O}') \to \bar{O} \:\:\: ([g, z] \to Ad_g(z)),$$ which we call a generalized Springer map. Let $(X')^0 \to O'$ be an etale covering (which is not 
necessarily the universal covering) and let $X' \to \bar{O}'$ be the associated finite 
cover.  Then we can consider the space $\mathfrak{n} + X'$ which is, by definition, a product of an affine space $\mathfrak{n}$ and the affine variety $X'$. There is a 
finite cover $\mathfrak{n} + X' \to \mathfrak{n} + \bar{O}'$. If the $Q$-action on $\mathfrak{n} + \bar{O}'$ 
lifts to a Q-action on $\mathfrak{n} + X'$, then we can make $G \times^Q(\mathfrak{n} 
+ X')$ and get a commutative diagram 
\begin{equation} 
\begin{CD} 
G \times^Q(\mathfrak{n} + X') @>{\mu'}>> Z \\ 
@V{\pi'}VV @VVV \\ 
G \times^Q (\mathfrak{n} + \bar{O}') @>{\mu}>> \bar{O},       
\end{CD} 
\end{equation}  
where $Z$ is the Stein factorization of $\mu \circ \pi'$. 
For an arbitrary nilpotent orbit $O$ of a classical Lie algebra $\mathfrak{g}$, we will give an explicit algorithm for finding $Q$, $O'$ and $X'$   
such that  

(1) $O = \mathrm{Ind}^{\mathfrak g}_{\mathfrak l}(O')$, 

(2) $X'$ has only {\bf Q}-factorial terminal singularities, and 

(3) the Q-action on $\mathfrak{n} + \bar{O}'$ lifts to a Q-action on 
$\mathfrak{n} + X'$ and the finite covering $Z \to \bar{O}$ in the diagram coincides with the finite covering $\pi: X \to \bar{O}$ associated with the universal covering $X^0$ of $O$.  

Then $\mu'$ gives a {\bf Q}-factorial terminalization of 
$X$. As for the existence of such ($Q$, $O'$, $X'$), there is another interesting view point \cite{Lo}. Along the argument of \cite{Lo}, Matvieievskyi \cite{Ma}, Corollary 4.3 assures the existence in a more general situation 
by using the universal Poisson deformation of $X$ \cite{Na 5}.   

Once we get an explicit {\bf Q}-factorial terminalization $G \times^Q(\mathfrak{n} + X')$ of $X$, it is easy to see if it is nonsingular or not (cf. Lemma \ref{Lemma (1.6)}).  If it is nonsingular, $\mu'$ is a crepant 
projective resolution of $X$. 
If it is singular, then any other {\bf Q}-factorial terminalizations of $X$ are also singular by \cite{Na 3}, Corollary 25; hence $X$ has no crepant projective 
resolution in such a case.  

In the subsequent article \cite{Part 2}, we will determine the Weyl group $W(X)$ of $X$ and count how many different {\bf Q}-factorial terminalizations $X$ has.  \vspace{0.2cm}

\setcounter{section}{-1}
\section{Preliminaries}


(P.1) {\bf Symplectic varieties}

Let $X$ be a normal variety over $\mathbf{C}$. 
Let $\omega$ be a regular 2-form on $X_{reg}$. 
Then $(X, \omega)$ is a {\em symplectic variety} if 

(a) $\omega$ is non-degenerate and $d$-closed, and 

(b) for a resolution $f: \tilde{X} \to X$ of $X$, the 2-form $f^*\omega$ 
on $f^{-1}(X_{reg})$ extends to a regular 2-form on $\tilde{X}$. 

A symplectic variety $X$ has only canonical singularities; hence it has only 
rational singularities.  The following properties of a symplectic variety 
will be frequently used in this article.

\begin{Prop} \label{Proposition (0.1)} Let $(X, \omega)$ be a symplectic variety of $\dim 2n$ and 
let $f: \tilde{X} \to X$ be a resolution of $X$. Write $K_{\tilde X} = f^*K_X + 
\sum a_iE_i$ with $f$-exceptional prime divisors $E_i$ (each coefficient $a_i$ being called 
the discrepancy of $E_i$). Let $E_{i_0}$ be an 
$f$-exceptional prime divisor with $a_{i_0} = 0$.  Then $\dim f (E_{i_0}) = 2n-2$. \end{Prop}

{\em Proof}. Put $S := f(E_{i_0})$. We blow up $\tilde{X}$ further to get 
a projective birational morphism $\nu: Z \to \tilde{X}$ such that 
$F := (f \circ \nu)^{-1}(S)$ is a simple normal crossing divisor of $Z$. 
$F$ contains the proper transform $F_0$ of $E_{i_0}$ by $\nu$ as an 
irreducible component.  
By definition $\omega$ lifts to a regular 2-form $\tilde{\omega}$ on $Z$. 
Since the discrepancy of $F_0$ is zero, $\tilde{\omega}$ is a non-degenerate 2-form on $Z$ at a general point $p \in F_0$. 
We may assume that $p$ is a smooth point of $F_0$.  Let $\phi$ be  the defining equation of $F_0$ at $p$. Then one can take a system of local parameters $\phi$, $\phi_2$, ... 
, $\phi_{2n}$ of $Z$ at $p$ in such a way that 
$$\tilde{\omega}(p) = d\phi \wedge d\phi_2 + ... + d\phi_{2n-1}\wedge d\phi_{2n} \in \wedge^2 (T^*Z)_p.$$ 
 This implies that  
$\wedge^{n-1} \tilde{\omega}\vert _{F_0} \ne 0$. According to \cite{Fr}, $\S 1$, we put $\hat{\Omega}^i_F := \Omega^i_F/\tau_i$ for $i \geq 0$, where 
$\tau_i$ is the 
subsheaf of $\Omega^i_F$ consisting of the sections supported on $\mathrm{Sing}(F)$. In particular, $\hat{\Omega}^0_F = 
\mathcal{O}_F$. 
Then $\tilde{\omega}\vert_F$ determines a nonzero element of 
$H^0(F, \hat{\Omega}^2_F)$.  Take a (non-empty) smooth open set $U$ of $S$ so that the fiber $F_x$ of $(f \circ \nu)\vert_F : F \to S$ over $x \in U$ 
is a simple normal crossing variety for any $x \in U$ and that $F$ is locally a product of 
$F_x$ and $U$. 
Replace $S$ by $U$ and $F$ by $(f \circ \nu\vert_F)^{-1}(U)$.  
Define $$\mathcal{G} := \mathrm{Ker}[\hat{\Omega}^2_F \to \hat{\Omega}^2_{F/S}].$$ 
Then there are exact sequences 
$$ 0 \to \mathcal{G} \to \hat{\Omega}^2_F \to \hat{\Omega}^2_{F/S} \to 0$$ 
$$ 0 \to  (f \circ \nu)\vert_F^*\Omega^2_S \to \mathcal{G} \to 
(f \circ \nu)\vert_F^*\Omega^1_S \otimes \hat{\Omega}^1_{F/S} \to 0.$$
This can be checked as follows. 
Let $F_0$, ..., $F_k$ be the irreducible components of $F$. For each $p$ with $p \le k$, we put 
$$F^{[p]} := \coprod_{0 \le i_0 < i_1 <... < i_p \leq k}F_{i_0} \cap \cdot \cdot \cdot \cap F_{i_p}.$$
Denote by $\mu_p : F^{[p]} \to F$ the natural map. For simplicity, we write $\alpha$ for 
$(f\circ \nu)\vert_F$, and $\alpha_p$ for  
$(f\circ \nu)\vert_F \circ \mu_p$.    
Note that $F^{[p]}$ is a disjoint union of smooth varieties. If we put 
$$\mathcal{G}_p := \mathrm{Ker}[\Omega^2_{F^{[p]}} \to \Omega^2_{F^{[p]}/S}],$$
then there are exact sequences 
$$0 \to \mathcal{G}_p \to \Omega^2_{F^{[p]}} \to \Omega^2_{F^{[p]}/S} \to 0$$
$$0 \to \alpha_p^*\Omega^2_S \to \mathcal{G}_p \to \alpha_p^*\Omega^1_S \otimes \Omega^1_{F^{[p]}/S} \to 0.$$
By taking $(\mu_p)_*$ of these exact sequences, we get the exact sequences 
$$0 \to (\mu_p)_*\mathcal{G}_p \to (\mu_p)_*\Omega^2_{F^{[p]}} \to (\mu_p)_*\Omega^2_{F^{[p]}/S} \to 0$$ 
$$0 \to \alpha^*\Omega^2_S \otimes (\mu_p)_*\mathcal{O}_{F^{[p]}} \to (\mu_p)_*\mathcal{G}_p \to 
\alpha^*\Omega^1_S \otimes (\mu_p)_*\Omega^1_{F^{[p]}/S} \to 0$$
On the other hand, the sheaves $\hat{\Omega}^i_F$ and $\hat{\Omega}^i_{F/S}$ respectively have 
resolutions (see \cite{Fr}, Proposition (1.5) for details):  
$$0 \to \hat{\Omega}^i_F \to (\mu_0)_*\Omega^i_{F^{[0]}} \to (\mu_1)_*\Omega^i_{F^{[1]}} \to ... $$ 
$$0 \to \hat{\Omega}^i_{F/S} \to (\mu_0)_*\Omega^i_{F^{[0]}/S} \to (\mu_1)_*\Omega^i_{F^{[1]}/S} \to ... $$

By combining these exact sequences, we finally get the following commutative diagrams with exact rows and 
exact columns

\begin{equation} 
\begin{CD}
... @.  ... @. ... @. ... @.  ...\\
@. @AAA @AAA  @AAA  @. \\  
0 @>>> (\mu_1)_*\mathcal{G}_1 @>>> (\mu_1)_*\Omega^2_{F^{[1]}} @>>> (\mu_1)_*\Omega^2_{F^{[1]}/S} @>>> 0 \\
@.  @AAA @AAA @AAA @.  \\ 
0  @>>> (\mu_0)_*\mathcal{G}_0 @>>> (\mu_0)_*\Omega^2_{F^{[0]}} @>>> (\mu_0)_*\Omega^2_{F^{[0]}/S} @>>> 0 \\ 
@. @AAA @AAA @AAA  @. \\ 
0 @>>> \mathcal{G} @>>> \hat{\Omega}^2_F @>>> \hat{\Omega}^2_{F/S} @>>> 0\\
@. @AAA @AAA  @AAA @. \\  
@.  0 @. 0 @. 0 @.  @.       
\end{CD} 
\end{equation}

\begin{equation} 
\begin{CD}
... @.  ... @. ... @. ... @.  ...\\
@. @AAA @AAA  @AAA  @. \\  
0 @>>> \alpha^*\Omega^2_S \otimes (\mu_1)_*\mathcal{O}_{F^{[1]}}  @>>> (\mu_1)_*\mathcal{G}_1 @>>> 
 \alpha^*\Omega^1_S \otimes (\mu_1)_*\Omega^1_{F^{[1]}/S} @>>> 0 \\
@.  @AAA @AAA @AAA @.  \\ 
0  @>>> \alpha^*\Omega^2_S \otimes (\mu_0)_*\mathcal{O}_{F^{[0]}}  @>>> (\mu_0)_*\mathcal{G}_0 @>>> 
\alpha^*\Omega^1_S \otimes (\mu_0)_*\Omega^1_{F^{[0]}/S}  @>>> 0 \\ 
@. @AAA @AAA @AAA  @. \\ 
0 @>>> \alpha^*\Omega^2_S \otimes \mathcal{O}_F  @>>> \mathcal{G} @>>> \alpha^*\Omega^1_S \otimes \hat{\Omega}^1_{F/S}
 @>>> 0\\
@. @AAA @AAA  @AAA @. \\  
@.  0 @. 0 @. 0 @.  @.       
\end{CD} 
\end{equation}

The exact sequences at the bottom rows are nothing but the desired exact sequences.  
 
We derive a contradiction by assuming that $i :=  \mathrm{Codim}_X S \geq 3$. 
$\tilde{\omega}\vert_F \in H^0(F, \hat{\Omega}^2_F)$ cannot be written as 
the pull-back of a 2-form $\omega_S$ on $S$. In fact, since $\dim S \le 2n -3$, 
$\wedge^{n-1}\omega_S = 0$. If $\tilde{\omega}\vert_F$ is the pull-back of $\omega_S$,   
then $\wedge^{n-1} \tilde{\omega}\vert_F = 0$, 
which means that $\wedge^{n-1} \tilde{\omega}\vert_{F_0} = 0$. This contradicts 
that $\wedge^{n-1} \tilde{\omega}\vert_{F_0} \ne 0$.   
Since $\tilde{\omega}\vert_F \in H^0(F, \hat{\Omega}^2_F)$ is not the pull-back 
of a 2-form on $S$, we see that $H^0(F_x, \hat{\Omega}^2_{F_x}) \ne 0$ or 
$H^0(F_x, \hat{\Omega}^1_{F_x}) \ne 0$ for a general point 
$x \in S$ by the exact sequences above. 
On the other hand, since $(X, x)$ is a rational singularity, we have 
$H^0(F_x, \hat{\Omega}^p_{F_x}) = 0$ for all $p > 0$ by \cite{Na 2}, Lemma (1.2). 
This is a contradiction. Therefore $\mathrm{Codim}_X(S) = 2$.  
$\square$ 

\begin{Cor} \label{Corollary (0.2)} Let $\pi: Y \to X$ be a crepant partial resolution of a 
symplectic variety $X$ of  $\dim 2n$ and let $E$ be a $\pi$-exceptional prime divisor. Then $\dim \pi (E) = 2n-2$. \end{Cor} 

{\em Proof}. Let $\nu: \tilde{X} \to Y$ be a resolution and let $\tilde{E}$ 
be the proper transform of $E$ by $\nu$. Put $f = \pi \circ \nu$. 
Then the discrepancy of $\tilde{E}$ is $0$. By Proposition \ref{Proposition (0.1)}  
$\dim f(\tilde{E}) = 2n-2$. $\square$ \vspace{0.2cm}

\begin{Cor} \label{Corollary (0.3)} Let $(X, \omega)$ be a symplectic variety of $\dim 2n$ with 
$\mathrm{Codim}_X\mathrm{Sing}(X) \geq 4$. Then $X$ has only 
terminal singularities. \end{Cor} 

{\em Proof}. Let $f : \tilde{X} \to X$ be a resolution of $X$ such that 
$f\vert_{f^{-1}(X_{reg})}: f^{-1}(X_{reg}) \cong X_{reg}$. Assume that 
$X$ has worse singularities than terminal singularities. Then there is an $f$-exceptional prime 
divisor $F$ such that its discrepancy is $0$. By Proposition \ref{Proposition (0.1)}, $\dim f(F) = 2n-2$. 
This contradicts that $\mathrm{Codim}_X\mathrm{Sing}(X) \geq 4$.  
$\square$ 

Let $(X, \omega)$ be an affine symplectic variety with $R = \Gamma (X, \mathcal{O}_X)$. Assume that $R$ is a positively graded ring 
$R = \oplus_{i \geq 0} R_i$ with $R_0 = \mathbf{C}$. This means that $X$ has 
a $\mathbf{C}^*$-action such that the closed point $0 \in X$ corresponding to 
the maximal ideal $m_R := \oplus_{i > 0}R_i$ is a unique fixed point of the 
$\mathbf{C}^*$-action. If $\omega$ is homogeneous with respect to this 
$\mathbf{C}^*$-action (i.e. $t^*\omega = t^l\omega$ for some integer $l$ and 
for $t \in \mathbf{C}^*$), then we call $(X, \omega)$ a {\em conical symplectic 
variety}.  By the property (b) of a symplectic variety, the {\em weight} $l$ is 
a positive integer. 
\vspace{0.2cm}

(P.2) {\bf Induced orbits and generalized Springer maps} 

Let $Q \subset G$ be a parabolic subgroup of $G$ and let $Q = U\cdot L$ be a Levi 
decomposition of $Q$ by the unipotent radical $U$ and a Levi subgroup $L$. 
Correspondingly the Lie algebra $\mathfrak{q}$ decomposes 
$\mathfrak{q} = \mathfrak{n} \oplus \mathfrak{l}$ as a direct sum of  $\mathfrak{n} := Lie(U)$ and $\mathfrak{l} := Lie(L)$. Let $O'$ be a nilpotent orbit of $\mathfrak{l}$. 
Then there is a unique nilpotent orbit $O$ of $\mathfrak{g}$ such that $O$ meets 
$\mathfrak{n} + O'$ in a Zariski open subset of $\mathfrak{n} + O'$. In such a case 
we say that $O$ is induced from $O'$ and write $O = \mathrm{Ind}^{\mathfrak g}_{\mathfrak l}(O')$.  There is a generically finite map 
$$\mu: G \times^Q(\mathfrak{n} + \bar{O}') \to \bar{O} \:\:\: ([g, z] \to Ad_g(z)),$$ which we call a generalized Springer map. 
Let $\mu^0: \mu^{-1}(O) \to O$ be the induced map and 
let $\omega_{KK}$ be the Kirillov-Kostant form on $O$. Put $Y^0 := G \times^Q(\mathfrak{n} + O')$.  
Then $(\mu^0)^*\omega_{KK}$ extends to a symplectic 2-form 
on $Y^0$ by \cite{Na 6}, Proposition 4.2 (see also \cite{Na 1}, Lemma (1.2.4)). By taking the wedge product of the symplectic forms, we get 
$K_{Y^0} = (\mu\vert_{Y^0})^*K_{\bar O}$. Take the normalization $\tilde{O}'$ of $\bar{O}'$ and put 
$Y := G \times^Q(\mathfrak{n} + \tilde{O}')$. Let $\tilde{O}$ be the normalization of $\bar{O}$. Then the induced 
map $\mu^n:  Y \to \tilde{O}$ satisfies $K_Y = (\mu^n)^*K_{\tilde O}$. In particular, when $\mu$ is birational,   
$Y$ is a crepant partial resolution of $\tilde{O}$. 
\vspace{0.2cm}    

(P.3) {\bf Nilpotent orbits and their finite coverings}
 
Let $O \subset \mathfrak{g}$ be a nilpotent orbit of a complex semisimple Lie algebra $\mathfrak{g}$. Take a simply connected complex algebraic 
group $G$ with $Lie(G) = \mathfrak{g}$.  
Let $\pi_0: X^0 \to O$ be a finite etale covering. Then the $G$-action on $O$ extends to a $G$-action on $X^0$ and $\pi_0$ is a $G$-equivariant covering.  
Let $\bar{O}$ be the closure 
of $O \subset \mathfrak{g}$.  We can extend $\pi_0$ to a finite covering map $\pi: X \to \bar{O}$. Since $X = \mathrm{Spec}\Gamma (X^0, \mathcal{O}_{X^0})$ and $G$ acts 
on $X^0$, $G$ acts on $X$ in such a way that $\pi$ is a $G$-equivariant.

For a point $x \in O$, let $G^x \subset G$ be the stabilizer 
group of $x$. We take the universal covering of $O$ as $X^0$ and pick a point 
$\tilde{x} \in X^0$ such that $\pi_0(\tilde{x}) = x$. Let  
$G^{\tilde x}$ be the stabilizer group of $\tilde{x}$ for the $G$-action on $X^0$. Then $G^{\tilde x}$ coincides with the  
identity component of $G^x$. Therefore $\pi_1(O) \cong 
G^x/(G^x)^0$. By the Jacobson-Morozov theorem we take an $sl(2)$-triple $x$, $y$ and 
$h$ in $\mathfrak{g}$. We denote by $\phi$ the map $sl(2) \to \mathfrak{g}$ determined by the $sl(2)$-triple. Put $$\mathfrak{g}^{\phi} := \{z \in \mathfrak{g} \: \vert \: [z, x] = [z, y] 
= [z, h] = 0\}.$$  Obviously $\mathfrak{g}^{\phi} \subset \mathfrak{g}^x$. Let $\mathfrak{u}^x$ 
be the nilradical of $\mathfrak{g}^x$. Note that $\mathfrak{u}^x$ is the Lie algebra of the unipotent radical of $G^x$. 
By Barbasch-Vogan and Kostant (cf. \cite{C-M}, Lemma 
3.7.3), there is a direct sum decomposition $\mathfrak{g}^x = \mathfrak{u}^x \oplus \mathfrak{g}^{\phi}$. Correspondingly we have a semi-direct product $G^x = U^x \cdot G^{\phi}$. 
Moreover, the inclusion $G^{\phi} \to G^x$ induces an isomorphism 
$$G^{\phi}/(G^{\phi})^0 \cong G^x/(G^x)^0,$$ which, in particular, means that 
$\pi_1(O) \cong G^{\phi}/(G^{\phi})^0$.

A nilpotent orbit closure $\bar{O}$ has a natural scaling $\mathbf{C}^*$-action for which $\omega_{KK}$ has weight $1$. Let $s: \mathbf{C}^* 
\to \mathrm{Aut}(\bar{O})$ be a homomorphism determined by the scaling 
action. In general this $\mathbf{C}^*$-action does not lift to a $\mathbf{C}^*$-action on $X$. But, if we instead define a new $\mathbf{C}^*$-action on $\bar{O}$ by the composite $\sigma$ of $$\mathbf{C}^* \to \mathbf{C}^* \:\: (t \to t^2), \:\:\: 
\mathrm{and} \:\:\: s: \mathbf{C}^* \to \mathrm{Aut}(\bar{O}),$$
then the new $\mathbf{C}^*$-action $\sigma$ always lifts to a $\mathbf{C}^*$-action on $X$ by 
\cite{B-K}, \S 1. By definition, $wt(\omega) = 2$ with respect to this $\mathbf{C}^*$-action. 
$\pi$ is etale in codimension one because $\mathrm{Codim}_{\bar{O}}\bar{O} - O \geq 2$. The map $\pi$ factors through the normalization $\tilde{O}$ of $\bar{O}$. 
Put $\omega := \pi_0^*\omega_{KK}$. Then this means that $(X, \omega)$ is a 
symplectic variety. In fact, take a resolution $f: Y \to \bar{O}$ so that $f^{-1}(O) \cong O$, and make a 
commutative diagram 
\begin{equation} 
\begin{CD} 
X \times_{\tilde O}Y  @>>> Y \\ 
@V{f'}VV @V{f}VV \\ 
X @>{\pi}>> \bar{O}       
\end{CD} 
\end{equation}  
Then $Z^0 := X^0 \times_O f^{-1}(O)$ is isomorphically mapped to $X^0$ by $f'$. 
Let $Z$ be the closure of $Z^0$ in $X \times_{\bar O}Y$, and let $\tilde{Z}$ be a resolution of 
$Z$. There is a commutative diagram 
\begin{equation} 
\begin{CD} 
\tilde{Z}  @>{\tilde \pi}>> Y \\ 
@V{\tilde f}VV @V{f}VV \\ 
X @>{\pi}>> \bar{O}       
\end{CD} 
\end{equation} 
where $\tilde{Z}$ and $Y$ are both smooth varieties.  
Since $\tilde{O}$ has symplectic singularities, $\omega_{KK}$ extends to a 2-form $\omega_Y$ on $Y$. 
Then $\tilde{\pi}^*\omega_Y$ is a 2-form on $\tilde{Z}$ and we see that the 2-form $\omega$ on $X^0$ 
extends to the 2-form $\tilde{\pi}^*\omega_Y$ on $\tilde{Z}$. This means that $(X, \omega)$ is a symplectic variety. 
Since $wt(\omega) = 2$, $(X, \omega)$ is a conical symplectic 
variety. \vspace{0.2cm}

\section{$\mathfrak{g} = sl(d)$}

A nilpotent orbit $O$ of $sl(d)$ is uniquely determined by its Jordan type. 
If the Jordan normal form of  $x \in O$ has $j_1$ Jordan blocks of size $d_1$, 
$j_2$ Jordan blocks of size $d_2$, ..., and $j_k$ Jordan blocks of size $d_k$,  
then the Jordan type of $O$ is a partition $[d_1^{j_1}, ..., d_k^{j_k}]$ of $d$. 
In the remainder we assume that $d_1 > d_2 > ... > d_k$. We indicate by 
$O_{[d_1^{j_1}, ..., d_k^{j_k}]}$ the nilpotent orbit with Jordan type $[d_1^{j_1}, ..., d_k^{j_k}]$.     

We first consider the case when $O$ is the regular nilpotent orbit $O_{[d]}$ of $sl(d)$. 
In this case $\bar{O}_{[d]}$ is the nilpotent cone of $sl(d)$. \vspace{0.2cm}

\begin{Prop} \label{Proposition (1.1)}  

(1) $\pi_1(O_{[d]}) \cong \mathbf{Z}/d\mathbf{Z}$. 

(2) \em $X$ is $\mathbf{Q}$-factorial for any etale covering $\pi_0: X^0 \to O_{[d]}$. \end{Prop}

{\em Proof}.  Take $x \in O_{[d]}$ and consider an $sl(2)$-triple $\phi$. 

(1) As already remarked above, $\pi_1(O_{[d]}) = G^{\phi}/(G^{\phi})^0$. 
By Springer and Steinberg (cf. \cite{C-M}, Theorem 6.1.3), one has 
$$G^{\phi} \cong \{(\zeta, ..., \zeta) \in GL(1)^d \: \vert \: \zeta^d = 1\} (\cong \mathbf{Z}/d\mathbf{Z}). $$  

(2)  It is enough to prove that $\mathrm{Pic}(X^0)$ is a finite group for $X^0$.  $X^0$ can be written as $G/H$ with a subgroup $H$ with $(G^x)^0 \subset H \subset G^x$. Then we have an exact sequence (cf. \cite{KVV}, Proposition 3.2) 
$$ \chi (H) \to \mathrm{Pic}(G/H) \to \mathrm{Pic}(G),$$ 
where $\chi (H) = \mathrm{Hom}(H, \mathbf{C}^*)$.  
Since $\mathrm{Pic}(G)$ is finite, we need to show that $\chi (H)$ is finite. 
By the exact sequence $$1 \to U^x \to (G^x)^0 \to (G^{\phi})^0 \to 1$$
we have an exact sequence $$\chi ((G^{\phi})^0) \to \chi ((G^x)^0) \to 
\chi (U^x).$$ Since $(G^{\phi})^0 = 1$, the 1-st term is zero. The 3-rd term is 
zero because $U^x$ is unipotent. Hence $\chi ((G^x)^0) = 0$. 
Now, by the exact sequence $$ \chi (H/(G^x)^0) \to \chi (H) \to \chi ((G^x)^0)$$ 
we see that $\chi (H)$ is finite because $\chi (H/(G^x)^0)$ is finite.

\begin{Prop} \label{Proposition (1.2)}  
Assume that $\pi_0: X^0 \to O_{[d]}$ is the universal covering. Then $\mathrm{Codim}_X\mathrm{Sing}(X) \geq 4$. 
In particular, $X$ has only terminal singularities. \end{Prop} 

{\em Proof}. Take a point $z$ from the subregular nilpotent orbit $O_{[d-1, 1]}$ and 
let $S$ be a (complex analytic) transverse slice for $O_{[d-1,1]} \subset \bar{O}_{[d]}$ at $z$ (cf. \cite{K-P 1}, Theorem 3.2). 
Then $S$ is a surface with an $A_{d-1}$-singularity at $z$. This means that the complex analytic germ $(S, z)$ is a 2 dimensional quotient singularity $V_d$, where  
$$V_d := (\mathbf{C}^2/(\mathbf{Z}/d\mathbf{Z}), 0).$$ 
Here $\bar{1} \in \mathbf{Z}/d\mathbf{Z}$ acts on $\mathbf{C}^2$ by 
$$(x, y) \to (\zeta x, \zeta^{-1} y)$$ with $\zeta$ a primitive $d$-th root of unity. 
Note that, for any $e$ with $e \vert d$, the quotient map $(\mathbf{C}^2, 0) \to V_d$ 
factorizes as $(\mathbf{C}^2, 0) \to V_e \to V_d$.    

The inclusion map $S - \{z\} \to O_{[d]}$ induces a homomorphism $\pi_1(S - \{z\}) 
\to \pi_1(O_{[d]})$. We prove that it is an isomorphism. Suppose to the contrary. Then 
$\pi^{-1}(S)$ splits into more than one connected components, each of which 
is a copy of $V_e$ with some divisor $e (\ne d)$ of $d$; namely,    
$$\pi^{-1}(S) = V_e^{(1)} \sqcup ... \sqcup V_e^{(e)}, $$ 
$$\pi_0^{-1}(S - \{z\}) = (V_e^{(1)} - \{0\}) \sqcup ... \sqcup (V_e^{(e)} - \{0\}). $$
If we put $f := d/e$, then $V_e^{(i)} \to S$ is a cyclic cover of degree $f$.   

First we shall construct a cyclic covering $\nu: Y \to O_{[d]}$ of degree $e$ in such a way that 
$\nu$ is etale not only over $O_{[d]}$ but also over $O_{[d-1, 1]}$. 
Take a point $p \in O_{[d]}$. Then $\pi^{-1}(p)$ consists of exactly $d$ points $\{p_1, ..., p_d\}$
and $\pi_1(O_{[d]})$ acts on them. 
We may renumber these points so that  

(a) each $V_e^{(i)}$ contains $p_j$ with $j = i \:\: (\mathrm{mod}\: e)$, and 

(b) $\pi_1(O_{[d]}) ( =\mathbf{Z}/d\mathbf{Z})$ acts on these $d$ points so that the 
generator $\bar{1}$ acts as the  
permutation $(1,2, ..., d)$; namely, $p_1 \to p_2$, $p_2 \to p_3$, ..., $p_{d-1} \to p_d$, 
$p_d \to p_1$. \vspace{0.2cm}
  
The natural surjection $\mathbf{Z}/d\mathbf{Z} \to \mathbf{Z}/e\mathbf{Z}$ determines 
a $\mathbf{Z}/e\mathbf{Z}$-covering $\nu: Y \to \bar{O}_{[d]}$. 
By definition $\nu^{-1}(S)$ splits up into $e$ copies of $V_d$. This is the desired cyclic 
covering. 

On the other hand, $\bar{O}_{[d]}$ can be resolved by the cotangent bundle $W := T^*(SL(d)/B)$ of 
the full flag variety $SL(d)/B$. We call this resolution the Springer resolution and denote it by 
$s: W \to \bar{O}_{[d]}$. The Springer resolution $s$ is a symplectic resolution;  hence, $s$ is a semi-small map. Put $Z := \bar{O}_{[d]} - O_{[d]} - O_{[d-1, 1]}$.
Then, by the semi-smallness of $s$, we have $\mathrm{Codim}_{W}s^{-1}(Z) 
\geq 2$. Since $\pi_1(W) = \{1\}$, we have $\pi_1(W - s^{-1}(Z)) = \{1\}$. 
By the construction $\nu : Y \to \bar{O}_{[d]}$ is etale over $\bar{O}_{[d]} - Z$. 
Then one can construct a (connected) non-trivial etale cover of $W - s^{-1}(Z)$ by pulling 
back $\nu$ by the map $W - s^{-1}(Z) \to \bar{O}_{[d]} - Z$. This contradicts that 
$\pi_1(W - s^{-1}(Z)) = \{1\}$. 
\vspace{0.2cm}

In the following we construct a {\bf Q}-factorial terminalization of $X$ for an arbitrary etale cover $\pi_0 :X^0 \to O_{[d]}$. Let us begin with the simplest cases. 

\begin{Exams} \label{Examples (1.3)} \end{Exams}
 (1) When $\mathfrak{g} = sl(2)$, $\pi_1(O_{[2]}) = \mathbf{Z}/2\mathbf{Z}$. 
Then $X = \mathbf{C}^2$ for the universal covering $X^0 \to O_{[2]}$, and $\pi: \mathbf{C}^2 
\to \bar{O}_{[2]}$ is the quotient map of $\mathbf{Z}/2\mathbf{Z}$ by the action $(x_1, x_2) 
\to (-x_1, -x_2)$.  

(2) When $\mathfrak{g} = sl(4)$, $\pi_1(O_{[4]}) = \mathbf{Z}/4\mathbf{Z}.$ 
When $X^0$ is the universal cover of $O_{[4]}$, we already know that $X$ has only  
$\mathbf{Q}$-factorial terminal singularities by Proposition \ref{Proposition (1.1)}, (2) and Proposition \ref{Proposition (1.2)}.
 
Next assume that $X^0$ is a double cover of $O_{[4]}$. 
The nilpotent cone $\bar{O}_{[4]}$ has a Springer resolution $T^*(SL(4)/Q_{1,1,1,1})$. 
Here $Q_{1, 1,1,1}$ is a parabolic subgroup of $SL(4)$ stabilizing a flag $0 \subset F^3 
\subset F^2 \subset F^1 \subset \mathbf{C}^4$ with $\dim Gr^i_F = 1$ for all $i$. 
Let $\mathfrak{n}_{1,1,1,1}$ be the nilradical of $Q_{1,1,1,1}$. By using the Killing form of $sl(4)$, we have an identification  
$$T^*(SL(4)/Q_{1,1,1,1}) \cong SL(4) \times^{Q_{1,1,1,1}} \mathfrak{n}_{1,1,1,1}.$$
Let $Q_{2,2}$ be the parabolic subgroup of $SL(4)$ stabilizing the flag $0 \subset F^2 \subset \mathbf{C}^4$. Then $Q_{2,2}$ has a Levi decomposition $$Q_{2,2} = U \cdot L,$$ where $U$ is 
the unipotent radical of $Q_{2,2}$ and $L$ is a Levi part of $Q_{2,2}$. In our case 
$$U = \{\left(\begin{array}{cc} 
I_2 & * \\
0 & I_2 
\end{array}\right)\}.$$
$$L = \{\left(\begin{array}{cc} 
A & 0 \\
0 & B 
\end{array}\right) \: \vert \: \mathrm{det}(A)\mathrm{det}(B) = 1 \}.$$

Corresponding to the decomposition we have a direct sum decomposition of Lie algebras
$$ \mathfrak{q}_{2,2} = \mathfrak{n} \oplus \mathfrak{l}.$$ In our case 
$$\mathfrak{l} = sl(2)^{\oplus 2} \oplus \mathfrak{z}, $$ where $\mathfrak{z}$ is a 
1-dimensional center of $\mathfrak{l}$. Take a nilpotent orbit closure $\bar{O}_{[2]} 
\times \bar{O}_{[2]}$ in $sl(2)^{\oplus 2} \oplus \mathfrak{z}$. 
The parabolic subgroup $Q_{2,2}$ acts on $\mathfrak{q}_{2,2}$. Since $\mathfrak{n}$ is an ideal of $\mathfrak{q}_{2,2}$, $\mathfrak{n}$ is stable under the $Q_{2,2}$-action. On the other hand, 
$\mathfrak{l}$ is not stable. Let $z \in \mathfrak{l}$ and $q \in Q_{2,2}$. Then $Ad_q(z)$ 
decomposes into the sum of  the nilradical part $(Ad_q(z))_n$ and the Levi part
$(Ad_q(z))_l$.  Let $Q_{1,1} \to Q_{1,1}/U 
= L$ be the quotient map and let $\bar{q} \in L$ be the image of $q \in Q_{2,2}$ by this map. 
The Levi subgroup $L$ acts on $\mathfrak{l}$ by the adjoint action.   
Then $$(Ad_q(z))_l = Ad_{\bar{q}}(z).$$ In particular, if $z \in \bar{O}_{[2]} \times \bar{O}_{[2]}$, 
then $(Ad_q(z))_l \in \bar{O}_{[2]} \times \bar{O}_{[2]}$. 
Therefore $Q_{2,2}$ acts on $\mathfrak{n} + (\bar{O}_{[2]} \times \bar{O}_{[2]})$. 
Then $SL(4) \times^{Q_{2,2}} (\mathfrak{n} + \bar{O}_{[2]} \times \bar{O}_{[2]})$ gives a crepant partial resolution of $\bar{O}_{[4]}$. Moreover, the Springer resolution of $O_{[4]}$ factors through this partial resolution: 
$$SL(4) \times^{Q_{1,1,1,1}} \mathfrak{n}_{1,1,1,1} \to SL(4) \times^{Q_{2,2}} (\mathfrak{n} + \bar{O}_{[2]} \times \bar{O}_{[2]}) \to \bar{O}_{[4]}.$$
We want to make a double cover $X'$ of $SL(4) \times^{Q_{2,2}} (\mathfrak{n} + \bar{O}_{[2]} \times \bar{O}_{[2]})$ so that the diagram  

\[
  \begin{diagram}
    \node{X'} \arrow{e,t}{\mu} \arrow{s,l}{\pi'} \node{X} \arrow{s,l}{\pi} \\
    \node{SL(4) \times^{Q_{2,2}} (\mathfrak{n} + \bar{O}_{[2]} \times \bar{O}_{[2]})} \arrow{e} \node{\bar{O}_{[4]}}
  \end{diagram}
\]
commutes and $\mu$ gives a $\mathbf{Q}$-factorial terminalization of $X$.  
By (1) we have a finite covering $$\mathbf{C}^2 \times \mathbf{C}^2 \to \bar{O}_{[2]} \times 
\bar{O}_{[2]}$$ of degree 4.  Then $\mathbf{Z}/2\mathbf{Z}$ acts on $\mathbf{C}^4 ( = 
\mathbf{C}^2 \times \mathbf{C}^2)$ by $(x_1, x_2, y_1, y_2) \to (-x_1, -x_2, -y_1, -y_2)$.
We denote by $\mathbf{C}^4/\langle +1,-1 \rangle$ the quotient space.   
The covering map above then factors through $\mathbf{C}^4 /\langle +1,-1 \rangle$:  $$\mathbf{C}^4  \to \mathbf{C}^4 /\langle +1,-1 \rangle \to \bar{O}_{[2]} 
\times \bar{O}_{[2]}.$$ 
We want to make $\mathfrak{n} + \mathbf{C}^4 /\langle +1,-1 \rangle $ into 
a $Q_{2,2}$-space and to define $X'$ to be $SL(4) \times^{Q_{2,2}} (\mathfrak{n} + 
\mathbf{C}^4 /\langle +1,-1 \rangle )$. 
\vspace{0.2cm}

\begin{Claim} \label{Claim (1.3.1)} The adjoint action of $L$ on $\bar{O}_{[2]} \times \bar{O}_{[2]}$ 
lifts to an action on $\mathbf{C}^4 /\langle +1,-1 \rangle $. \end{Claim} 

{\em Proof}. As already remarked above, 
$$L = \{\left(\begin{array}{cc} 
A & 0 \\
0 & B 
\end{array}\right) \: \vert \: \mathrm{det}(A)\mathrm{det}(B) = 1 \}.$$  
Define a subgroup $T$ of $L$ by 
$$T = \{\left(\begin{array}{cc} 
\lambda I_2 & 0 \\
0 & \lambda^{-1}I_2 
\end{array}\right) \: \vert \: \lambda \in \mathbf{C}^* \}.$$  
We identify $SL(2) \times SL(2)$ with a subgroup of $L$ by  
$$\{\left(\begin{array}{cc} 
A & 0 \\
0 & B 
\end{array}\right) \: \vert \: A, B \in SL(2) \}.$$
Then the inclusion map $SL(2) \times SL(2) \subset L$ induces an isomorphism 
$$SL(2) \times SL(2)/\langle +1, -1 \rangle \cong L/T.$$
$SL(2) \times SL(2)$ naturally acts on $\mathbf{C}^4$; hence 
$SL(2) \times SL(2)/\langle +1, -1 \rangle$ acts on $\mathbf{C}^4/\langle +1,-1 \rangle$. 
As a consequence, $L/T$ acts on $\mathbf{C}^4/\langle +1,-1 \rangle$. In particular, 
$L$ acts on $\mathbf{C}^4/\langle +1,-1 \rangle$, which is a lift of the adjoint action 
of $L$ on $\bar{O}_{[2]} \times \bar{O}_{[2]}$. $\square$
\vspace{0.2cm}

Now $Q_{2,2}$ acts on the space $\mathfrak{n} + \mathbf{C}^4/\langle +1,-1 \rangle$ as 
follows. Take a point $z + v$ from the space. Here $z \in \mathfrak{n}$ and 
$v \in \mathbf{C}^4/\langle +1,-1 \rangle$. We denote by $\bar{v}$ the image of 
$v$ by the map $\mathbf{C}^4/\langle +1,-1 \rangle \to \bar{O}_{[2]} \times \bar{O}_{[2]}$.  
For $q \in Q_{2,2}$ we denote by $\bar{q} 
\in L$ the image of $q$ by the map $Q_{2,2} \to Q_{2,2}/U = L$. 
We define 
$$ q\cdot (z + v) := (Ad_q(z + \bar{v}))_n + \bar{q} \cdot v \in \mathfrak{n} + \mathbf{C}^4/\langle +1,-1 \rangle.$$
Here $\bar{q} \in L$ acts on $v \in \mathbf{C}^4/\langle +1,-1 \rangle$ as described in 
Claim \ref{Claim (1.3.1)}.   
Let us consider the composed map   
$$SL(4) \times^{Q_{2,2}} (\mathfrak{n} + 
\mathbf{C}^4 /\langle +1,-1 \rangle ) \to  
SL(4) \times^{Q_{2,2}} (\mathfrak{n} + \bar{O}_{[2]} \times \bar{O}_{[2]}) 
\to \bar{O}_{[4]}.$$
The Stein factorization of this map is nothing but $X$. 
As a consequence we have a commutative diagram 

\begin{equation} 
\begin{CD} 
SL(4) \times^{Q_{2,2}} (\mathfrak{n} + 
\mathbf{C}^4 /\langle +1,-1 \rangle )  @>{\mu}>> X \\ 
@V{\pi'}VV @V{\pi}VV \\ 
SL(4) \times^{Q_{2,2}} (\mathfrak{n} + \bar{O}_{[2]} \times \bar{O}_{[2]}) @>>>  \bar{O}_{[4]}     
\end{CD} 
\end{equation} 

We can generalize this construction to more general situations. 
Let us consider the regular nilpotent orbit $O_{[d]}$ of $sl(d)$. 
Assume that $e$ is a divisor of $d$. We put $f := d/e$. By Proposition \ref{Proposition (1.1)}, (1) $\pi_1(O_{[d]}) = 
\mathbf{Z}/d\mathbf{Z}$. The surjective homomorphism $\mathbf{Z}/d\mathbf{Z} \to 
\mathbf{Z}/e\mathbf{Z}$ determines an etale cover $X^0 \to O_{[d]}$ of degree $e$. 
We will construct a {\bf Q}-factorial terminalization of $X$ by the same idea of 
Examples \ref{Examples (1.3)}. Let $Q_{e,e, ..., e}$ be a parabolic subgroup of $SL(d)$ with flag type $(e, ..., e)$. 
Let $Q_{e,..., e} = U \cdot L$ be a Levi decomposition where 
   
$$U = \{\left(\begin{array}{ccccc} 
I_e & * & * & * & *\\
0 & I_e & * & * & *\\ 
... & ... & ... & ... & ... \\ 
0 & 0 & 0 & I_e & * \\ 
0 & 0 & 0 & 0 & I_e 
\end{array}\right)\},$$

$$L = \{\left(\begin{array}{ccccc} 
A_1 & 0 & 0 & 0 & 0\\
0 & A_2 & 0 & 0 & 0\\ 
... & ... & ... & ... & ... \\ 
0 & 0 & 0 & A_{f-1} & 0 \\ 
0 & 0 & 0 & 0 & A_f 
\end{array}\right) \: \vert \:  A_i \in GL(e), \: \mathrm{det}(A_1) \cdot \cdot \cdot  \mathrm{det}(A_f) = 1 \}.$$
We have $$ \mathfrak{q}_{e, ..., e} = \mathfrak{n} \oplus \mathfrak{l}, $$
$$\mathfrak{l} = sl(e)^{\oplus f} \oplus \mathfrak{z},$$ where 
$\mathfrak{z}$ is the $f-1$ dimensional center of $\mathfrak{l}$. 
We take a nilpotent orbit closure $\bar{O}_{[e]}^{\times f}$ 
inside $\mathfrak{l}$. Then 
$$SL(d) \times^{Q_{e, ..., e}} (\mathfrak{n} + (\bar{O}_{[e]}^{\times f})$$ is a crepant partial resolution of $\bar{O}_{[d]}$. 
Let $X_{[e]} \to \bar{O}_{[e]}$ be the finite covering of degree $e$ corresponding to 
the universal covering of $O_{[e]}$. The adjoint action of $SL(e)$ on $\bar{O}_{[e]}$ extends 
to an action on $X_{[e]}$. The center $Z$ of $SL(e)$ is written as 
$$\{\zeta I_e \in SL(e) \: \vert \: \zeta^e = 1\}.$$ 
Then $Z$ acts effectively on $X_{[e]}$ as covering transformations of 
$X_{[e]} \to \bar{O}_{[e]}$. In fact,  let $x \in O_{[e]}$ and put $G := SL(e)$. 
Then the universal cover $X_{[e]}^0$ of $O_{[e]}$ is written as 
$G/(G^x)^0$. Now $G^x/(G^x)^0 \cong G^{\phi}/(G^{\phi})^0$. By the description of 
$G^{\phi}$ in the proof of Proposition \ref{Proposition (1.1)}, (1) we see that $G^{\phi} = Z$ and 
$(G^{\phi})^0 = \{1\}$. Since $G^x/(G^x)^0$ acts effectively on $G/(G^x)^0$, we see 
that $Z$ acts effectively on $G/(G^x)^0$; hence, acts effectively on $X_{[e]}$. 
Then $Z^{\times f}$ acts on $X_{[e]}^{\times f}$. Define a subgroup 
$S(Z^{\times f})$ of $Z^{\times f}$ by  
$$S(Z^{\times f}) := \{(\zeta_1I_e, ..., \zeta_fI_e) \in Z^{\times f} \: \vert \: \zeta_1\cdot \cdot \zeta_f = 1\}.$$
Notice that $S(Z^{\times f}) \cong (\mathbf{Z}/e\mathbf{Z})^{\oplus f-1}$.  
Then the map $X_{[e]}^{\times f} \to \bar{O}_{[e]}^{\times f}$ factors through $X_{[e]}^{\times  f}/S(Z^{\times f})$: 
$$X_{[e]}^{\times f} \to X_{[e]}^{\times f}/S(Z^{\times f}) \to \bar{O}_{[e]}^{\times f}.$$ By definition the 2-nd map is a $\mathbf{Z}/e\mathbf{Z}$-Galois 
covering. \vspace{0.2cm}

\begin{Claim} \label{Claim (1.3.2)} The adjoint action of $L$ on $\bar{O}_{[e]}^{\times f}$ lifts to an 
action on $X_{[e]}^{\times f}/S(Z^{\times f})$. \end{Claim}

{\em Proof}. 
Define a subgroup $T$ of $L$ by  
$$T := \{\left(\begin{array}{ccccc} 
\zeta_1I_e & 0 & 0 & 0 & 0\\
0 & \zeta_2I_e & 0& 0 & 0\\ 
... & ... & ... & ... & ... \\ 
0 & 0 & 0 & \zeta_{f-1}I_e & 0 \\ 
0 & 0 & 0 & 0 & \zeta_fI_e 
\end{array}\right)\: \vert \: \zeta_i \in \mathbf{C}^*, \: \zeta_1\cdot \cdot \cdot \zeta_f = 1\}.$$ 
We identify $SL(e)^{\times f}$ with a subgroup of $L$ defined by      
$$\{\left(\begin{array}{ccccc} 
A_1 & 0 & 0 & 0 & 0\\
0 & A_2 & 0 & 0 & 0\\ 
... & ... & ... & ... & ... \\ 
0 & 0 & 0 & A_{f-1} & 0 \\ 
0 & 0 & 0 & 0 & A_f 
\end{array}\right) \: \vert \:  A_i \in SL(e) \: \forall  i \}.$$ 
Then the inclusion $SL(e)^{\times f} \to L$ induces an isomorphism 
$SL(e)^{\times f}/S(Z^{\times f}) \cong L/T$. As $SL(e)^{\times f}$ acts on $X_{[e]}^{\times f}$,  $L/T$ acts on $X_{[e]}^{\times f}/S(Z^{\times f})$. Hence $L$ acts on 
$X_{[e]}^{\times f}/S(Z^{\times f})$.  $\square$
\vspace{0.2cm}

Now $Q_{e, ..., e}$ acts on $\mathfrak{n} + X_{[e]}^{\times f}/S(Z^{\times f})$ as follows. 
Take a point $z + v \in \mathfrak{n} + X_{[e]}^{\times f}/S(Z^{\times f})$. Here $z \in \mathfrak{n}$ and 
$v \in X_{[e]}^{\times f}/S(Z^{\times f})$. We denote by $\bar{v}$ the image of 
$v$ by the map $X_{[e]}^{\times f}/S(Z^{\times f}) \to \bar{O}_{[2]}^{\times f}$.  
For $q \in Q_{e, ..., e}$ we denote by $\bar{q} 
\in L$ the image of $q$ by the map $Q_{e,..., e} \to Q_{e,...,e}/U = L$. 
We define 
$$ q\cdot (z + v) := (Ad_q(z + \bar{v}))_n + \bar{q} \cdot v \in \mathfrak{n} + X_{[e]}^{\times f}/S(Z^{\times f}).$$
Here $\bar{q} \in L$ acts on $v \in X_{[e]}^{\times f}/S(Z^{\times f})$ as described in 
Claim \ref{Claim (1.3.2)}.

Recall that $\pi: X \to \bar{O}_{[d]}$ is a finite $\mathbf{Z}/e\mathbf{Z}$-cover.  
We have a commutative diagram 
\begin{equation} 
\begin{CD} 
SL(d) \times^{Q_{e, ..., e}} (\mathfrak{n} + X_{[e]}^{\times f}/S(Z^{\times f})) @>{\mu}>> X \\ 
@V{\pi'}VV @V{\pi}VV \\ 
SL(d) \times^{Q_{e, ..., e}} (\mathfrak{n} + \bar{O}_{[2]}^{\times f}) @>{s}>>  \bar{O}_{[d]}     
\end{CD} 
\end{equation}     
Here $\mu$ is the Stein factorization of $s \circ \pi'$. 
$X_{[e]}^{\times f}$ is $\mathbf{Q}$-factorial by Propositon \ref{Proposition (1.1)}, (2). 
Then $X_{[e]}^{\times f}/S(Z^{\times f})$ is also $\mathbf{Q}$-factorial by the 
following lemma. 

\begin{Lem} \label{Lemma (1.4)} Let $f: V \to W$ be a finite Galois covering of normal varieties.  
If $V$ is {\bf Q}-factorial, then $W$ is also {\bf Q}-factorial. \end{Lem} 

{\em Proof}. Let $D$ be a prime Weil divisor of $W$.  We need to show that $mD$ is a Cartier divisor for a suitable $m$. Put $E := f^{-1}(D)$ and regard it as a reduced Weil divisor. By the assumption $rE$ is a Cartier divisor on $V$ for some $r > 0$.  Take a point $x \in W$.
By \cite{Mum}, Lecture 10, Lemma B, one can choose an open neighborhood 
$U$ of $x \in W$ such that $rE$ is a principal divisor of $f^{-1}(U)$. We take a defining 
equation $\varphi$of $rE\vert_{f^{-1}(U)}$. Let $G$ be the Galois group of $f$. Then 
$$\Phi := \prod_{g \in G} \varphi^g$$ is a local equation of the Cartier divisor $\vert G \vert r E$. Since $\Phi$ is $G$-invariant, it can be regarded as an element of 
$\Gamma (U, \mathcal{O}_W)$. Then $\Phi$ is a local equation of some multiple of 
$D$ on $U$.  $\square$
\vspace{0.2cm}

The variety $SL(d) \times^{Q_{e, ..., e}} (\mathfrak{n} + X_{[e]}^{\times f}/S(Z^{\times f}))$ is a 
fiber bundle over $SL(d)/Q_{e, ..., e}$ with a typical fiber $\mathfrak{n} + X_{[e]}^{\times f}/S(Z^{\times f})$. By \cite{Ha}, II, Proposition 6.6, we have 
$$\mathrm{Cl}( X_{[e]}^{\times f}/S(Z^{\times f}) \cong \mathrm{Cl}(\mathfrak{n} + X_{[e]}^{\times f}/S(Z^{\times f})),$$ where $\mathrm{Cl}$ denotes the divisor class 
group.   By using this we see that  $\mathfrak{n} + X_{[e]}^{\times f}/S(Z^{\times f})$ is {\bf Q}-factorial since $X_{[e]}^{\times f}/S(Z^{\times f})$ is {\bf Q}-factorial,. Then $SL(d) \times^{Q_{e, ..., e}} (\mathfrak{n} + X_{[e]}^{\times f}/S(Z^{\times f}))$ is also {\bf Q}-factorial by the following lemma.  \vspace{0.2cm}

\begin{Lem}\label{Lemma (1.5)} Let $f: V \to T$ be an etale fiber bundle over a nonsingular variety $T$ with a typical fiber $Y$. Assume that 

 (1) $Y$ is a {\bf Q}-factorial normal variety. 

 (2) $Y$ has only rational singularities with $\mathrm{Codim}_Y\mathrm{Sing}(Y) \geq 
3$.  

Then $V$ is also $\mathbf{Q}$-factorial. \end{Lem}

{\em Proof}. Take a closed point $v \in V$ and put $t = f(v)$. 
Replace $T$ by a suitable open neighborhood of $t$. Put $n = \dim T$. Then 
one has a sequence of nonsingular subvarieties $\{t\} \subset T_1 \subset T_2 
\subset ... \subset T_{n-1} \subset T_n = T$ with $\dim T_i = i$. 
Put $V_i := V \times_T T_i$. Then we get a sequence 
$V_0 (= Y) \subset V_1 \subset V_2 \subset ... \subset V_n = V$. Here each $V_i$ is a 
Cartier divisor of $V_{i+1}$.  Since rational singularities are Cohen-Macaulay, we can apply \cite{Ko-Mo}, Corollary (12.1.9) for $Y \subset V_1$ 
to see that $V_1$ is {\bf Q}-factorial around $V_0$. Now $V_1$ satisfies the 
conditions (1) and (2). Therefore we can apply [ibid, Corollary (12.1.9)] repeatedly for $V_i \subset V_{i+1}$ and finally see that $V$ is {\bf Q}-factorial around $V_{n-1}$. 
In particular, $V$ is {\bf Q}-factorial around $f^{-1}(t) \subset V$. Since $v$ is an arbitrary closed point, $V$ is {\bf Q}-factorial. $\square$      
\vspace{0.2cm}   

\begin{Lem}\label{Lemma (1.6)} For $e > 2$, $X_{[e]}$ is singular. For $e > 1$ and $f > 1$, $X_{[e]}^{\times f}/S(Z^{\times f})$ 
is singular. \end{Lem} 

{\em Proof}. Assume that $X_{[e]}$ is smooth for $e > 2$. Since $X_{[e]} \to \bar{O}_{[e]}$ 
is a finite quotient map, $\bar{O}_{[e]}$ would be a symplectic quotient singularity. 
Then the closure of any symplectic leaf of $\bar{O}_{[e]}$ would be again a quotient singularity.   
On the other hand, the closure $\bar{O}_{[2, 1^{e-2}]}$ of the minimal nilpotent orbit $O_{[2, 1^{e-2}]}$ is not 
a quotient singularity if $e > 2$. In fact, the Springer resolution $s$ of $\bar{O}_{[2, 1^{e-2}]}$ 
is given by the cotangent bundle $T^*\mathbf{P}^{e-1}$ of $\mathbf{P}^{e-1}$. The exceptional 
locus of the Springer resolution is the zero locus of the cotangent bundle, which has codimension $e-1 (\geq 2)$.  
This means that $\bar{O}_{[2, 1^{e-2}]}$ is not $\mathbf{Q}$-factorial. 
In fact,  take an $s$-ample effective divisor $H$ on $T^*\mathbf{P}^{e-1}$ and consider $s_*H$. If $\bar{O}_{[2, 1^{e-2}]}$ is $\mathbf{Q}$-factorial, then $ms_*H$ is Cartier for some $m > 0$.  Then $s^*(ms_*H)$ and $mH$ coincides 
outside $\mathrm{Exc}(s)$, which has codimension $> 1$. Hence, they coincide 
on $T^*\mathbf{P}^{e-1}$. Take a curve $C \subset \mathrm{Exc}(s)$ so that 
$s(C)$ is a point. Then $(s^*(ms_*H). C) = 0$. On the other hand, $(mH, C) > 0$ 
because $H$ is $s$-ample. This is a contradiction.   
In particular, $\bar{O}_{[2, 1^{e-2}]}$ is not a quotient singularity. Therefore $X_{[e]}$ is singular. 
When $e > 2$, we can prove similarly that $X_{[e]}^{\times f}/S(Z^{\times f})$ is singular by using 
the quotient map $X_{[e]}^{\times f}/S(Z^{\times f}) \to \bar{O}_{[e]}^{\times f}$. When $e = 2$, $X_{[e]} = \mathbf{C}^2$. 
But, if $f > 1$, then any non-zero element of $S(Z^{\times f})$ has the fixed locus of codimension $\geq 4$. Hence
$X_{[2]}^{\times f}/S(Z^{\times f})$ is singular.  $\square$
\vspace{0.2cm}

Since $X_{[e]}$ has terminal singularities, the product $X_{[e]}^{\times f}$ has terminal singularities.  The fixed locus of any nonzero element of $S(Z^{\times f})$ has codimension $\geq 4$; hence $X_{[e]}^{\times f}/S(Z^{\times f})$ has only terminal singularities. 

\begin{Lem}\label{Corollary (1.7)} Let $O_{[d]} \subset sl(d)$ be the regular nilpotent orbit with $d > 2$. 
Assume that $X^0 \to O_{[d]}$ is an etale covering of degree $e > 1$ and let $X \to 
\bar{O}_{[d]}$ be the associated finite cover of $\bar{O}_{[d]}$. Then the map 
 $$SL(d) \times^{Q_{e, ..., e}} (\mathfrak{n} + X_{[e]}^{\times f}/S(Z^{\times f})) \to 
X$$ is a {\bf Q}-factorial terminalization of $X$. In particular, $X$ has 
no crepant resolutions. \end{Lem}  

We next consider the nilpotent orbit $O_{[d^i]}$ of $sl(di)$. Put $G = SL(di)$. As before , take an element $x$ from $O_{[d^i]}$ and fix an $sl(2)$-triple $\phi$ containing $x$. 
Then $$G^{\phi} = \{(A, ..., A) \in GL(i)^{\times d} \: \vert \: \mathrm{det}(A)^d 
= 1\}. $$ $G^{\phi}$ has exactly $d$ connected components.  
Let $\zeta$ be a primitive $d$-th root of unity. Then each connected component is given by $$\{(A, ..., A) \in GL(i)^{\times d} \: \vert \: \mathrm{det}(A) = \zeta^i\},\:  i = 0,1, ...,d-1.$$  
Note that $G^{\phi}/(G^{\phi})^0 \cong \mathbf{Z}/d\mathbf{Z}$.  

\begin{Prop}\label{Proposition (1.8)} 

(1) $\pi_1(O_{[d^i]}) \cong \mathbf{Z}/d\mathbf{Z}$. 

(2) $X$ is $\mathbf{Q}$-factorial for any etale covering $\pi_0: X^0 \to O_{[d^i]}$. \end{Prop}

{\em Proof}. (1) We have already seen that (1) holds. 

(2) We can prove (2) in the same manner as in Proposition (1.1), (2).  Note that 
$\chi ((G^{\phi})^0) = \{0\}$ because $(G^{\phi})^0 \cong SL(i)^{\times d}$. $\square$
\vspace{0.2cm}

The closure $\bar{O}_{[d^i]}$ contains the largest orbit $O_{[d^{i-1}, d-1, 1]}$ in 
$\bar{O}_{[d^i]} - O_{[d^i]}$. Note that $\bar{O}_{[d^i]}$ has $A_{d-1}$-surface singularities along 
$O_{[d^{i-1}, d-1, 1]}$ (cf. \cite{K-P 1}, Theorem 3.2). Moreover, $\bar{O}_{[d^i]}$ has a Springer resolution. Since the universal 
covering of $O_{[d^i]}$ is a cyclic covering of degree $d$, we can prove the following in the same way as Proposition \ref{Proposition (1.2)}. 

\begin{Prop}\label{Proposition (1.9)}  
Assume that $\pi_0: X^0 \to O_{[d^i]}$ is the universal covering. Then $\mathrm{Codim}_X\mathrm{Sing}(X) \geq 4$. 
In particular, $X$ has only terminal singularities.
\end{Prop}
 
For a partition $[d_1^{j_1}, ..., d_k^{j_k}]$ of $j_1d_1 + ... + j_kd_k$, consider the nilpotent orbit $O_{[d_1^{j_1}, ..., d_k^{j_k}]} \subset sl(j_1d_1 + ... + j_kd_k)$. Let $d := gcd(d_1, ..., d_k)$. By \cite{C-M}, Corollary 6.1.6, we have $\pi_1(O_{[d_1^{j_1}, ..., d_k^{j_k}]}) \cong \mathbf{Z}/d\mathbf{Z}$. Let $\pi^0: X^0 \to O_{[d_1^{j_1}, ..., d_k^{j_k}]}$ be the universal covering and let $\pi: X \to \bar{O}_{[d_1^{j_1}, ..., d_k^{j_k}]}$ be the associated finite cover. We will construct a {\bf Q}-factorial terminalization 
of $X$ by using Propositon \ref{Proposition (1.9)}.

\begin{Exam}\label{Example (1.10)} \end{Exam} Let $Q_{i_1d, ..., i_rd} \subset SL((i_1 + ... +i_r)d)$ be parabolic subgroup 
of flag type $(i_1d, ..., i_rd)$. We assume that $gcd(i_1, ..., i_r) = 1$. 
Let $Q_{i_1d, ..., i_rd} = U \cdot L$ be a Levi decomposition 
with   
$$ U = \{\left(\begin{array}{ccccc} 
I_{i_1d} & * & .* & ... & * \\
0 & I_{i_2d} & * & ... & * \\ 
... & ... & ... & ... & ... \\ 
0 & 0 & ... & I_{i_{r-1}d} & * \\ 
0 & 0 & ... & 0 & I_{i_rd} 
\end{array}\right)\} $$ 
and 
$$ L = \{\left(\begin{array}{ccccc} 
A_1 & 0 & .0 & ... & 0 \\
0 & A_2 & 0 & ... & 0 \\ 
... & ... & ... & ... & ... \\ 
0 & 0 & ... & A_{r-1} & 0 \\ 
0 & 0 & ... & 0 & A_r 
\end{array}\right) \: \vert A_1 \in GL(i_1d), \: ..., \: A_r \in GL(i_rd), \: \mathrm{det}(A_1)\cdot \cdot
\cdot \mathrm{det}(A_r) = 1\}. $$
The subgroup of $L$ consisting of the matrices with $A_1 \in SL(i_1d)$, ..., $A_r \in 
SL(i_rd)$ is isomorphic to $SL(i_1d) \times ... \times SL(i_rd)$.  In the remainder 
we identify $SL(i_1d) \times ... \times SL(i_rd)$ with this subgroup of $L$. 
Let us consider the product of the nilpotent orbits $\bar{O}_{[d^{i_1}]} \times ... \times 
\bar{O}_{[d^{i_r}]}$ in $sl(i_1d) \times ... \times sl(i_rd)$. Let $X_{[d^{i_{\alpha}}]} \to \bar{O}_{[d^{i_{\alpha}}]}$ 
be the finite cover associated with the universal covering of $O_{[d^{i_{\alpha}}]}$ for each 
$1 \le \alpha \le r$. Then we have a 
finite covering $$f: X_{[d^{i_1}]} \times ... \times X_{[d^{i_r}]} \to \bar{O}_{[d^{i_1}]} \times ... \times \bar{O}_{[d^{i_r}]}.$$
We consider the finite subgroup $\mu_{i_1d} \times ... \times \mu_{i_rd}$ of $SL(i_1d) \times ... \times SL(i_rd)$ consisting of the 
elements $(t_1I_{i_1d}, ..., t_rI_{i_rd})$ with $t_1^{i_1d} = ... = t_r^{i_rd} = 1$. 
Moreover let $H'$ be the subgroup of $\mu_{i_1d} \times ... \times \mu_{i_rd}$ determined  
by $t_1^{i_1}... t_r^{i_r} = 1$. Define a surjection 
$$\alpha: \mu_{i_1d} \times ... \times \mu_{i_rd} \to \mu_d \times ... \times \mu_d$$ by 
$(t_1, ..., t_r) \to (t_1^{i_1}, ..., t_r^{i_r})$. Similarly, define a surjection 
$$\beta: \mu_d \times ... \times \mu_d \to \mu_d$$ by $(t_1, ..., t_r) \to t_1... t_r$. 
Then $H'$ is nothing but the kernel of the composed map $$\mu_{i_1d} \times ... \times \mu_{i_rd} \to \mu_d.$$
     
\begin{equation} 
\begin{CD}
@.  0 @. 0 @. @.  @.  \\
@. @VVV @VVV   @. \\  
0 @>>> \mu_{i_1} \times ... \times \mu_{i_r} @>{\iota}>> H' @>>> \mathrm{Coker}(\iota) \\
@.  @VVV @VVV @.   \\ 
0  @>>> \mu_{i_1d} \times ... \times \mu_{i_rd} @>{\cong}>> \mu_{i_1d} \times ... \times \mu_{i_rd} @>>> 0  @.\\ 
@. @V{\alpha}VV @V{\beta \circ \alpha}VV  @. \\ 
\mathrm{Ker}(\beta) @>>> \mu_d \times ... \times \mu_d @>{\beta}>> \mu_d @>>>  0 \\
@. @VVV @VVV  @. \\  
@.  0 @. 0 @. @.  @.       
\end{CD} 
\end{equation}
Here $\mathrm{Ker}(\alpha) = \mu_{i_1} \times ... \times \mu_{i_r}$ and 
$\iota$ is the natural inclusion. 
By the snake lemma, $\mathrm{Coker}(\iota) \cong \mathrm{Ker}(\beta)$. We put $H := \mathrm{Ker}(\beta)$.  
The subgroup $\mu_{i_1d} \times ... \times \mu_{i_rd}$ of $SL(i_1d) \times ... \times SL(i_rd)$ acts on $X_{[d^{i_1}]} \times ... \times X_{[d^{i_r}]}$ as covering transformations of $f$.  But $\mathrm{Ker}(\alpha)$ acts 
trivially on $X_{[d^{i_1}]} \times ... \times X_{[d^{i_r}]}$.  This means that 
$\mu_d \times ... \mu_d$ acts on  $X_{[d^{i_1}]} \times ... \times X_{[d^{i_r}]}$, which is nothing 
but the Galois group of the map $f$. Since $\mathrm{Coker}(\iota) \cong 
\mathrm{Ker}(\beta)$, the finite covering $$(X_{[d^{i_1}]} \times ... \times X_{[d^{i_r}]}) /H \to \bar{O}_{[d^{i_1}]} \times ... \times \bar{O}_{[d^{i_r}]}$$ is a cyclic covering of degree 
$d$.  

\begin{Claim}\label{Claim (1.10.1)} The adjoint action of $L$ on $\bar{O}_{[d^{i_1}]} \times ... \times 
\bar{O}_{[d^{i_r}]}$ lifts to an action on  $(X_{[d^{i_1}]} \times ... \times X_{[d^{i_r}]}) /H$. 
\end{Claim}

{\em Proof}. Define a subgroup $T$ of $L$ by 
$$ T = \{\left(\begin{array}{ccccc} 
t_1I_{i_1d} & 0 & .0 & ... & 0 \\
0 & t_2I_{[i_2d]} & 0 & ... & 0 \\ 
... & ... & ... & ... & ... \\ 
0 & 0 & ... & t_{r-1}I_{i_{r-1}d} & 0 \\ 
0 & 0 & ... & 0 & t_rI_{i_rd} 
\end{array}\right) \: \vert \: t_1^{i_1}...t_r^{i_r} = 1 \}. $$
Then the inclusion $SL(i_1d) \times ... \times SL(i_rd) \to L$ induces an 
isomorphism $(SL(i_1d) \times ... \times SL(i_rd))/H \cong L/T$. Since 
$(SL(i_1d) \times ... \times SL(i_rd))/H$ acts on $(X_{[d^{i_1}]} \times ... \times X_{[d^{i_r}]}) /H$, 
$L$ acts on $(X_{[d^{i_1}]} \times ... \times X_{[d^{i_r}]})/H$, which gives a lift of the 
adjoint action. $\square$ \vspace{0.2cm} 

By Claim \ref{Claim (1.10.1)}, $Q_{i_1d, ..., i_rd}$ acts on $\mathfrak{n} +  (X_{[d^{i_1}]} \times ... \times X_{[d^{i_r}]})/H$. Then we have a cyclic covering 
$$SL((i_1 + ... + i_r)d) \times^{Q_{i_1d, ..., i_rd}}(\mathfrak{n} + (X_{[d^{i_1}]} \times ... \times X_{[d^{i_r}]})/H)$$ $$\to  SL((i_1 + ... + i_r)d) \times^{Q_{i_1d, ..., i_rd}}(\mathfrak{n} + \bar{O}_{[d^{i_1}]} 
\times ... \times \bar{O}_{[d^{i_r}]})$$ of degree $d$. \hspace{12.0cm} $\square$ \vspace{0.2cm}

Now look at a nilpotent orbit $O_{[d_1^{j_1}, ..., d_k^{j_k}]} \subset sl(j_1d_1 + ... + j_kd_k)$.
Put $d := \mathrm{gcd}(d_1, ..., d_k)$. Then the dual partition of $[d_1^{j_1}, ..., d_k^{j_k}]$ can be written as 
the form $[i_1^{d}, ...., i_r^{d}]$ by using suitable set of positive integers $\{i_1, ... i_r\}$. Here the same number may 
possibly appears more than once in $\{i_1, ..., i_r\}$. For example, the partition $[9,6]$ has 
the dual partition $[2^3, 2^3, 1^3]$. Note that $3 = \mathrm{gcd}(9,6)$. The dual partition 
of $i_j^{d}$ ($j = 1, ..., r$) equals $[d^{i_j}]$ ($j = 1, ..., r$). Now let us consider the 
product of the nilpotent orbit closures $\bar{O}_{[d^{i_1}]} \times ... \times \bar{O}_{[d^{i_r}]}$ 
in $sl(i_1d) \times ... \times sl(i_rd)$. Then $SL((i_1 + ... + i_r)d) \times^{Q_{i_1d, ..., i_rd}}(\mathfrak{n} + \bar{O}_{[d^{i_1}]} 
\times ... \times \bar{O}_{[d^{i_r}]})$ gives a crepant partial resolution 
of $\bar{O}_{[d_1^{j_1}, ..., d_k^{j_k}]}$. 

{\bf Construction of a Q-factorial terminalization}:  Let $X \to  \bar{O}_{[d_1^{j_1}, ..., d_k^{j_k}]}$ be 
the finite covering associated with the universal covering of $O_{[d_1^{j_1}, ..., d_k^{j_k}]}$. 
By using Example \ref{Example (1.10)}, we have a commutative diagram 
\begin{equation} 
\begin{CD} 
SL((i_1 + ... + i_r)d) \times^{Q_{i_1d, ..., i_rd}}(\mathfrak{n} + (X_{[d^{i_1}]} \times ... \times X_{[d^{i_r}]})/H) @>{\mu}>> X \\ 
@V{\pi'}VV @V{\pi}VV \\ 
SL((i_1 + ... + i_r)d) \times^{Q_{i_1d, ..., i_rd}}(\mathfrak{n} + \bar{O}_{[d^{i_1}]} 
\times ... \times \bar{O}_{[d^{i_r}]}) @>>> \bar{O}_{[d_1^{j_1}, ..., d_k^{j_k}]}      
\end{CD} 
\end{equation}
Since $X_{[d^{i_1}]} \times ... \times X_{[d^{i_r}]}$ has terminal singularities and 
the fixed locus of each nonzero element of $H$ has codimension $\geq 4$, 
$(X_{[d^{i_1}]} \times ... \times X_{[d^{i_r}]})/H$ has only terminal singularities. 
Since $X_{[d^{i_1}]} \times ... \times X_{[d^{i_r}]}$ is {\bf Q}-factorial, $(X_{[d^{i_1}]} \times ... \times X_{[d^{i_r}]})/H$ is {\bf Q}-factorial by Lemma \ref{Lemma (1.4)}. 
Therefore $\mu$ gives a {\bf Q}-factorial terminalization of $X$.  
 
Next assume that $X \to \bar{O}_{[d_1^{j_1}, ..., d_k^{j_k}]}$ is the finite covering associated 
with an etale covering of $O_{[d_1^{j_1}, ..., d_k^{j_k}]}$ of degree $e$. We put 
$f := d/e$. 
In this case we take instead the product of nilpotent orbit closures $\bar{O}_{[e^{i_1}]}^{\times f} \times ... \times \bar{O}_{[e^{i_r}]}^{\times f}$ and consider the finite covering 
$$ X_{[e^{i_1}]}^{\times f} \times ... \times X_{[e^{i_r}]}^{\times f} \to 
\bar{O}_{[e^{i_1}]}^{\times f} \times ... \times \bar{O}_{[e^{i_r}]}^{\times f}.$$
The situation being the same as the previous case, we define similarly a subgroup $H$ of 
$\mu_{i_1e}^{\times f} \times ... \times \mu_{i_re}^{\times}$. 
Then we have a commutative diagram 
\begin{equation} 
\begin{CD} 
SL((i_1 + ... + i_r)d) \times^{Q_{(i_1e)^f, ..., (i_re)^f}}(\mathfrak{n} + (X_{[e^{i_1}]}^{\times f} \times ... \times X_{[e^{i_r}]}^{\times f})/H) @>{\mu}>> X \\ 
@V{\pi'}VV @V{\pi}VV \\ 
SL((i_1 + ... + i_r)d) \times^{Q_{(i_1e)^f, ..., (i_re)^f}}(\mathfrak{n} + \bar{O}_{[e^{i_1}]}^{\times f} 
\times ... \times \bar{O}_{[e^{i_r}]}^{\times f}) @>>> \bar{O}_{[d_1^{j_1}, ..., d_k^{j_k}]}      
\end{CD} 
\end{equation}
and $\mu$ gives a {\bf Q}-factorial terminalization.  \vspace{0.2cm}

\begin{Cor}\label{Corollary (1.11)} Let $\pi: X \to \bar{O}$ be the finite covering associated with a nontrivial etale covering of a nilpotent orbit $O$ of $sl(d)$. Then $X$ has 
no crepant resolutions except when $\pi$ is the double covering $\mathbf{C}^2 \to 
\bar{O}_{[2]} \subset sl(2)$. \end{Cor}

\section{$\mathfrak{g} = sp(2n)$}


We write a partition $\mathbf{p}$ of $2n$ as $[d^{r_d}, (d-1)^{r_{d-1}}, ..., 2^{r_2}, 1^{r_1}]$ 
with $r_d \ne 0$. 
Other $r_i$ may possibly be zero; in such a case $i$ does not appear in the 
partition. If $r_i > 0$, then we call $i$ a member of the partition. Each nilpotent 
orbit of $sp(2n)$ is uniquely determined by its Jordan type. Such a Jordan type is a 
partition $\mathbf{p}$ of $2n$ with all odd members having even multiplicities (cf. \cite{C-M}, \S 5). 
Conversely, for each  such partition $\mathbf{p}$ of $2n$, there is  a nilpotent orbit with 
Jordan type $\mathbf{p}$. Let $b$ be the number of distinct even members of 
$\mathbf{p}$. Then we have \vspace{0.2cm}

\begin{Prop}\label{Proposition (2.1)}  
(1) $\pi_1(O_{\mathbf p}) \cong (\mathbf{Z}/2\mathbf{Z})^{\oplus b}$. 

(2) Assume that $r_i \ne 2$ for all even members $i$ of $\mathbf{p}$. 
Then $X$ is {\bf Q}-factorial for any etale covering $X^0 \to O_{\mathbf p}$. 
\end{Prop} 

{\em Proof}. Put $G = Sp(2n)$. Let $x \in O_{\mathbf p}$ and take an $sl(2)$-triple 
$\phi$ in $sp(2n)$ containing $x$. Put 
$$Sp(r_i)^{\times i}_{\Delta} := \{(A, ..., A) \in Sp(r_i)^{\times i} \: \vert \: A \in Sp(r_i)\}$$ and   
$$O(r_i)^{\times i}_{\Delta} :=  \{(A, ..., A) \in O(r_i)^{\times i} \: \vert \: A \in O(r_i)\}.$$ 
By [C-M, Theorem (6.1.3)] we have  
$$G^{\phi} \cong \prod_{i: \mathrm{odd}}Sp(r_i)^{\times i}_{\Delta} \times 
\prod_{i: \mathrm{even}}O(r_i)^{\times i}_{\Delta}.$$  
Hence $$(G^{\phi})^0 \cong \prod_{i: \mathrm{odd}}Sp(r_i)^{\times i}_{\Delta} \times 
\prod_{i: \mathrm{even}}SO(r_i)^{\times i}_{\Delta}.$$ 
An important remark is that each factor of the right hand side is a simple Lie group 
except that $SO(2) \cong \mathbf{C}^*$ and $SO(4)$ is a semisimple Lie group of type 
$A_1 + A_1$.  
Since $\pi_1(O_{\mathbf p}) \cong G^{\phi}/(G^{\phi})^0$, (1) is clear from the above. 
If the condition of (2) holds, then $(G^{\phi})^0$ does not have $SO(2)$ as a factor;  hence $\chi ((G^{\phi})^0) = 0$.   
The etale covering $X^0$ of $O_{\mathbf p}$ can be written as $G/H$ for a 
suitable subgroup $H$ with $(G^{\phi})^0 \subset H \subset G^{\phi}$. By the same 
argument as in Proposition (1.1), (2) we see that $\mathrm{Pic}(G/H)$ is a finite group, 
which means that $X$ is {\bf Q}-factorial. $\square$

\begin{Exam}\label{Example (2.2)} \end{Exam} Let $O_{[2n]}$ be the regular nilpotent orbit. By Proposition \ref{Proposition (2.1)}, (1) 
we have $\pi_1(O_{[2n]}) \cong \mathbf{Z}/2\mathbf{Z}$. Let $X \to \bar{O}_{[2n]}$ be 
the double covering associated with the universal covering $X^0 \to O_{[2n]}$. 
Put the $n \times n$ matrix 
$$J_n = \left(\begin{array}{ccccc} 
0 & 0 & 0 & ... & 1 \\
0 & 0 & ... & 1 & 0 \\ 
... & ... & ... & ... & ... \\ 
0 & 1 & ... & 0 & 0 \\ 
1 & 0 & ... & 0 & 0 
\end{array}\right). $$   
Then $$Sp(2n) = \{ A \in GL(2n) \: \vert \: A^t \left(\begin{array}{cc}
0 & J_n \\ 
-J_n & 0 
\end{array}\right) A = \left(\begin{array}{cc}
0 & J_n \\ 
-J_n & 0 
\end{array}\right)\}. $$
Now let us consider the isotropic flag 
$$0 \subset \langle e_1 \rangle \subset \langle e_1, e_2 \rangle ... \subset \langle e_1, e_2, ..., 
e_{2n-1} \rangle \subset \mathbf{C}^{2n}$$ and let $Q_{1^{n-1}, 2, 1^{n-1}}$ be 
the parabolic subgroup of $Sp(2n)$ stabilizing the flag. Let $U$ be the unipotent 
radical of $Q_{1^{n-1}, 2, 1^{n-1}}$. 
One has a Levi decomposition $Q_{1^{n-1}, 2, 1^{n-1}} = U \cdot L$ with 
$$ L = \{ \left(\begin{array}{ccccccccc} 
t_1 & 0 & 0 & ... & ... & ... & ... & ... & 0 \\
0 & t_2 & 0 & ... & ... & ... & ... & ... & 0 \\ 
... & ... & ... & ... & ... & ... & ... & ... & ...\\
0 & ... & 0 & t_{n-1} & 0 & ... & ... & ... & 0\\
0 & ... & ... & 0 & A & 0 & ... & ... & 0 \\
0 & ... & ... & ... & 0 & t^{-1}_{n-1} & 0 & ... & 0\\
... & ... & ... & ... & ... & ... & ... & ... & ... \\ 
0 & ... & ... & ... & ... & ... & ... & t_2^{-1} & 0 \\ 
0 & ... & ... & ... & ... & ... & ... & 0 & t_1^{-1} 
\end{array}\right) \: \vert \: t_1, ..., t_{n-1} \in \mathbf{C}^*, \: A \in Sp(2)\}.$$
The Lie algebra $\mathfrak{l}$ decomposes into the direct sum 
$$\mathfrak{l} = gl(1)^{\oplus n-1} \oplus sp(2).$$ 
Take a nilpotent orbit $O_{[2]}$ of $sp(2)$. Then we have a crepant partial resolution  
$$Sp(2n) \times^{Q_{1^{n-1}, 2, 1^{n-1}}} (\mathfrak{n} + \bar{O}_{[2]}) \to 
\bar{O}_{[2n]}.$$ 
Let $\mathbf{C}^2 \to \bar{O}_{[2]}$ be the double covering associated with the 
universal covering of $O_{[2]}$. The adjoint action of $Sp(2)$ on $\bar{O}_{[2]}$ lifts to an 
action on $\mathbf{C}^2$. Since there is a natural projection $L \to Sp(2)$, $L$ acts in the adjoint way on $\bar{O}_{[2]}$, which lifts to an $L$-action on $\mathbf{C}^2$. 
This means that $Q_{1^{n-1}, 2, 1^{n-1}}$ acts on $\mathfrak{n} + \mathbf{C}^2$. 
Then we have a commutative diagram 
\begin{equation} 
\begin{CD} 
Sp(2n) \times^{Q_{1^{n-1}, 2, 1^{n-1}}} (\mathfrak{n} + \mathbf{C}^2) @>{\mu}>> X \\ 
@V{\pi'}VV @V{\pi}VV \\ 
Sp(2n) \times^{Q_{1^{n-1}, 2, 1^{n-1}}} (\mathfrak{n} + \bar{O}_{[2]}) @>>> \bar{O}_{[2n]}      
\end{CD} 
\end{equation}
The map $\mu$ is a crepant resolution of $X$. 
$\square$ \vspace{0.2cm}

Let us consider a partition $\mathbf{p} = [d^{r_d}, (d-1)^{r_{d-1}}, ..., 2^{r_2}, 1^{r_1}]$ of  
$2n$ such that 

(i) $r_i$ is even for each odd $i $, and 

(ii) $r_i \ne 0$ for each even $i$.   

For example, the partitions $[5^2, 4, 2]$ and $[6, 4, 2, 1^4]$ satisfy these conditions, but 
$[8, 4, 2, 1^2]$ does not satisfy (ii) because $6$ does not appear. 
If $d$ is even, then all $d$, $d-2$, $d-4$, ..., $2$ must appear in the partition. If $d$ is odd, then all $d-1$, $d-3$, ..., $2$ must appear in the partition. The following is a key proposition. 

\begin{Prop}\label{Proposition (2.3)} Let $\mathbf{p}$ be a partition satisfying (i) and (ii). 
Let $\pi: X \to \bar{O}_{\mathbf p}$ be the finite covering associated with the universal 
covering of $O_{\mathbf p}$. Then $\mathrm{Codim}_X\mathrm{Sing}(X) 
\geq 4$. \end{Prop}

{\em Proof}. Let $i > 1$ be a member of $\mathbf{p}$ such that $r_{i-1} = 0$. This means that $\mathbf{p}$ 
has the form $$[d^{r_d}, ..., i^{r_i}, (i-2)^{r_{i-2}}, (i-3)^{r_{i-3}}, ..., 1^{r_1}]$$ 
Such an $i$ is called a {\em gap member}. By the condition (ii) any gap member must be even.    
For a gap member $i$, one can find a nilpotent orbit $O_{\bar{\mathbf{p}}} \subset 
\bar{O}_{\mathbf p}$ with  
$$\bar{\mathbf{p}} = [d^{r_d}, ..., i^{r_i - 1}, (i-1)^2, (i-2)^{r_{i-2} - 1}, (i-3)^{r_{i-3}}, ..., 1^{r_1}].$$   
By Kraft and Procesi \cite{K-P}, 3.4 we see that $\mathrm{Codim}_{\bar{O}_{\mathbf p}}O_{\bar{\mathbf{p}}} = 2$ and the transversal slice $S$ for $O_{\bar{\mathbf p}} \subset \bar{O}_{\mathbf p}$ is an $A_1$-surface singularity. Note that $S$ is a quotient 
singularity $(\mathbf{C}^2/(\mathbf{Z}/2\mathbf{Z}), 0)$. Then we have a double cover 
$(\mathbf{C}^2, 0) \to S$.

Let $\{i_1, i_2, ..., i_k\}$ be the set of all gap members of $\mathbf{p}$. As defined above, for these gap members, 
we have nilpotent orbits 
$$O_{\bar{\mathbf p}_1}, ..., O_{\bar{\mathbf p}_k}.$$  
Notice that $\bar{O}_{\bar{\mathbf p}_{1}}$, ..., $\bar{O}_{\bar{\mathbf p}_{k}}$ are 
nothing but the irreducible components of $\mathrm{Sing}(\bar{O}_{\mathbf p})$ which have  
codimension 2 in $\bar{O}_{\mathbf p}$. To prove that $\mathrm{Codim}_X\mathrm{Sing}(X) 
\geq 4$, we only have to show that $X$ is smooth along $\pi^{-1}(O_{\bar{\mathbf p}_{j}})$ for 
each $1 \le j \le k$. 
Let $S_j$ be the transversal slices for $O_{\bar{\mathbf p}_j} \subset 
\bar{O}_{\mathbf p}$. Then it is equivalent to showing that $\pi^{-1}(S_j)$ are disjoint union of finite copies of  $(\mathbf{C}^2, 0)$. 

The gap member is closely related to the notion of {\em induced orbits}. 
There are two types of inductions: 

(Type I): Let $i$ be a gap member of $\mathbf{p}$ (which may possibly be odd or even).  
Put $r := r_d + ... + r_i$ and 
let $Q \subset Sp(2n)$ be a parabolic subgroup of flag type $(r, 2n-2r, r)$ with Levi 
decomposition $\mathfrak{q} = \mathfrak{n} \oplus \mathfrak{l}$. Notice that 
$\mathfrak{l} = \mathfrak{gl}(r) \oplus sp(2n-2r)$. 
There is a nilpotent 
orbit $O_{\mathbf{p}'}$ of $sp(2n-2r)$ with Jordan type 
$$\mathbf{p}' = [(d-2)^{r_d}, ..., (i -2)^{r_i + r_{i -2}}, (i-3)^{r_{i-3}},..., 1^{r_1}]$$ such that 
$O_{\mathbf p} = \mathrm{Ind}^{\mathfrak g}_{\mathfrak l}(O_{\mathbf{p}'}).$  Notice that 
$\mathbf{p}'$ satisfies (i). 

\begin{Claim}\label{Claim (2.3.1)}(cf. \cite{He}, Theorem 7.1, (d)). The generalized Springer map 
$$\mu: Sp(2n) \times^{Q}(\mathfrak{n} + \bar{O}_{\mathbf{p}'}) \to \bar{O}_{\mathbf p}$$
is a birational map. \end{Claim}

{\em Proof}.  This is true for any partition $\mathbf{p}$ with the condition (i). 
We write $\mathbf{p}$ as  $[d_1, d_2, ..., d_s]$ with 
$d_1 \geq d_2 \geq ... \geq d_s >0$. 
The partition $\mathbf{p}$ determines a Young diagram $Y(\mathbf{p})$. By definition  
$Y(\mathbf p)$ is a subset of $\mathbf{Z}_{>0}^2$ such that $(l,j) \in Y(\mathbf{p})$ if and only if 
$(l,j)$ satisfies $1 \le l \le d_j$.       
 
We can take a basis $\{e(l,j)\}_{(l,j) \in Y(\mathbf{p})}$ of $\mathbf{C}^{2n}$ so that 

(a) $\{e(l,j)\}$ is a Jordan basis for $x$, i.e. 
$x\cdot e(l,j) = e(l-1,j)$ for $l > 1$ and $x\cdot e(1,j) = 0$. 

(b) $\langle e(l,j), e(p,q) \rangle \ne 0$ if and only if $p = d_j - l + 1$ and 
$q = \beta(j)$. Here $\beta$ is a permutation of $\{1,2, ..., s\}$ such that 
$\beta^2 = id$, $d_{\beta(j)} = d_j$, and $\beta (j) \ne j$ if $d_j$ is odd.  The basis $\{e(l,j)\}$ are the same one as in 
\cite{He}, 5.1 (cf. \cite{S-S}, p.259, see also \cite{C-M}, 5.1). The notation in [S-S] is slightly different, but we here employ the 
notation in \cite{He}, 5.1. 

Put $F:= \sum_{1 \le j \le r} \mathbf{C} e(1,j)$. Then 
$F \subset F^{\perp}$ is an isotropic flag such that $x\cdot F = 0$ and $x \cdot \mathbf{C}^{2n} \subset F^{\perp}$ and $x$ is an endomorphism of 
$F^{\perp}/F$ with Jordan type $\mathbf{p}'$.  This is actually a unique isotropic 
flag of type $(r, 2n-2r, r)$ satisfying these properties. Hence $\mu^{-1}(x)$ consists 
of one element. $\square$  \vspace{0.2cm}

(Type II): Let $i$ be an even member of $\mathbf{p}$ with $r_i = 2$. Put 
$r := r_d + ... + r_{i-1} + 1$ and let $Q \subset Sp(2n)$ be a parabolic subgroup of flag type 
$(r, 2n-2r, r)$ with Levi 
decomposition $\mathfrak{q} = \mathfrak{n} \oplus \mathfrak{l}$. Notice that 
$\mathfrak{l} = \mathfrak{gl}(r) \oplus sp(2n-2r)$. 
There is a nilpotent 
orbit $O_{\mathbf{p}'}$ of $sp(2n-2r)$ with Jordan type 
$$\mathbf{p}' = [(d-2)^{r_d}, ..., i^{r_{i+2}}, (i-1)^{r_{i+1} + 2 + r_{i-1}}, (i-2)^{r_{i-2}}, ..., 1^{r_1}]$$ 
such that $O_{\mathbf p} = \mathrm{Ind}^{\mathfrak g}_{\mathfrak l}(O_{\mathbf{p}'}).$
Notice that $r_{i+1} + 2 + r_{i-1}$ is even because $i$ is even; hence $\mathbf{p}'$ satisfies 
the condition (i). 
 
\begin{Claim}\label{Claim (2.3.2)}(cf. \cite{He}, Theorem 7.1, (d)).  The generalized Springer map 
$$\mu: Sp(2n) \times^{Q}(\mathfrak{n} + \bar{O}_{\mathbf{p}'}) \to \bar{O}_{\mathbf p}$$
is generically finite of degree $2$. \end{Claim} 

{\em Proof}. Fix an element $x \in O_{\mathbf p}$ and 
take the same basis $\{e(l,j)\}$ of $\mathbf{C}^m$ as the previous claim. 
We may assume that the permutation $\beta$ satisfies $\beta (r) = r+1$ and $\beta (r+1) = r$ 
after a suitable change of the basis.  

We put $F := \sum_{1 \le j \le r} \mathbf{C} e(1,j)$. 
Then 
$F \subset F^{\perp}$ is an isotropic flag such that $x\cdot F = 0$ and $x \cdot \mathbf{C}^{2n} \subset F^{\perp}$ and $x$ is an endomorphism of 
$F^{\perp}/F$ with Jordan type $\mathbf{p}'$.
 
On the other hand, put $F' := \sum_{1 \le j \le r+1, \: j \ne r}\mathbf{C} e(1, j)$. 
Then 
$F' \subset (F')^{\perp}$ is an isotropic flag such that $x\cdot F' = 0$ and $x \cdot \mathbf{C}^{2n} \subset (F')^{\perp}$ and $x$ is an endomorphism of 
$(F')^{\perp}/F'$ with Jordan type $\mathbf{p}'$. 

The isotropic flags of type $(r, 2n-2r, r)$ with these properties are 
exactly two flags above.   
Indeed, at first, it can be checked that 
such a flag $F$ contains the subspace $\sum_{1 \le j \le r-1}\mathbf{C}e(1, j)$. 
Then $F$ is written as $$F = \sum_{1 \le j \le r-1}\mathbf{C}e(1, j) + 
\mathbf{C}(\alpha e(1, r+1) + \beta e(1, r))$$ for some $(\alpha, \beta) 
\ne 0$. Put $\langle e(1, r+1), e(i, r) \rangle = a_1$, 
$\langle e(2, r+1), e(i -1, r) \rangle = a_2$, ..., 
$\langle e(i, r+1), e(1, r) \rangle = a_i$.  Here $a_1, ..., a_i$ are all nonzero. 
Since $x \in sp(2n)$, we have 
$a_2 = -a_1$, $a_3 = -a_2$, ..., $a_i = -a_{i-1}$. Since $i$ is even, 
$a_i = -a_1$.  Then we see that 
$$F^{\perp} := \sum_{1 \le j \le r-1, \: 1 \le l \le d_j - 1} \mathbf{C}e(l, j)
+ \mathbf{C}(\alpha e(l, r+1) - \beta e(l, r)).$$ 
If $\alpha$ and $\beta$ are both nonzero, then 
$x^{i-1}(\alpha e(i, r+1) - \beta e(i, r)) \ne 0$. This contradicts that 
$x$ is an endomorphism of $F^{\perp}/F$ with Jordan type 
$\mathbf{p}'$. Therefore, $\alpha = 0$ or $\beta = 0$. 
 
Assume that $Q$ is the parabolic subgroup of $Sp(2n)$ stabilizing the flag $F \subset F^{\perp}$.
We have an $Sp(2n)$-equivariant (locally closed) immersion  $$\iota: Sp(2n) \times^Q (\mathfrak{n} + O_{\mathfrak{p}'}) 
\subset Sp(2n)/Q \times sp(2n), \:\:\: [g, y] \to (gQ, Ad_g(y)).$$ 
Consider $(F \subset F^{\perp}, x)$ and $(F' \subset (F')^{\perp}, x)$ as elements 
of $Sp(2n)/Q \times sp(2n)$. Note that $\iota ([1, x]) = (Q, x) = (F \subset F^{\perp}, x)$. 
 
We want to prove that $(F' \subset (F')^{\perp}, x)$ is also contained in 
$\mathrm{Im}(\iota)$. Let $Q'$ be the parabolic subgroup stabilizing the flag $F' \subset (F')^{\perp}$. Then one can write $Q' = gQg^{-1}$ for some $g \in Sp(2n)$. Let 
$\mathfrak{q}' = \mathfrak{n}' \oplus \mathfrak{l}'$ be a Levi decomposition. 
By definition $x \in \mathfrak{n}' + O'_{\mathbf{p}'}$, where $O'_{\mathbf{p}'}$ 
is a nilpotent orbit in $\mathfrak{l}'$ with Jordan type $\mathbf{p}'$. 
Then $g^{-1}xg \in \mathfrak{n} + g^{-1}O'_{\mathbf{p}'}g$, where $g^{-1}O'_{\mathbf{p}'}g 
\in g^{-1}\mathfrak{l}' g$. The Lie algebra $g^{-1}\mathfrak{l}'g$ is a Levi subalgebra 
of $\mathfrak{q}$. Hence $\mathfrak{l}$ and $g^{-1}O'_{\mathbf{p}'}g$ are conjugate 
by an element of $Q$. By changing $g$ by $gq'$ for a suitable $q' \in Q$, we may 
assume from the first that $g^{-1}O'_{\mathbf{p}'}g \subset \mathfrak{l}$. 
Since $g^{-1}O'_{\mathbf{p}'}g$ and $O_{\mathbf{p}'}$ have the same Jordan type 
and the nilpotent orbits of $sp(2n-2r)$ are completely determined by their Jordan types, 
we see that $g^{-1}O'_{\mathbf{p}'}g = O_{\mathbf{p}'}$. 
This means that $$x,\: g^{-1}xg \: \in \mathfrak{n} + O_{\mathbf{p}'}.$$
The $Q$-orbit of $x$ and the $Q$-orbit of $g^{-1}xg$ are both dense in 
$\mathfrak{n} + O_{\mathbf{p}'}$. Hence they intersects. In other words, there is an 
element $q \in Q$ such that $g^{-1}xg = qxq^{-1}$. Then $gq \in Z_{Sp(2n)}(x)$ and 
$Q' = (gq)Q(gq)^{-1}$. Now let us consider an element $[gq, x] \in Sp(2n) \times^Q(\mathfrak{n} + \bar{O}_{\mathbf{p}'})$. Then $\iota ([gq, x]) = 
(gQ, x) = (F' \subset (F')^{\perp}, x)$.   
$\square$
\vspace{0.2cm}

Let $\{i_1, i_2, ..., i_k\}$ be the set of all gap members of $\mathbf{p}$ and let  $i_j$ be 
one of them.   
Put $r := r_d + ... + r_{i_j}$ and 
let $Q_j \subset Sp(2n)$ be a parabolic subgroup of flag type $(r, 2n-2r, r)$ with Levi 
decomposition $\mathfrak{q}_j = \mathfrak{n}_j \oplus \mathfrak{l}_j$. There is a nilpotent 
orbit $O_{\mathbf{p}'_j}$ of $\mathfrak{l}_j$ with Jordan type 
$$\mathbf{p}'_j = [(d-2)^{r_d}, ..., (i_j - 1)^{r_{i_j + 1}}, (i_j-2)^{r_{i_j} + r_{i_j-2}}, (i_j-3)^{r_{i_j-3}},..., 1^{r_1}]$$ such that 
$O_{\mathbf p} = \mathrm{Ind}^{\mathfrak g}_{\mathfrak l}(O_{\mathbf{p}'_j})$.  
We have a generically finite morphism called the {\em generalized Springer map} 
$$\mu_j: Sp(2n) \times^{Q_j}(\mathfrak{n}_j + \bar{O}_{\mathbf{p}'_j}) \to \bar{O}_{\mathbf p}.$$   
By the previous claim, $\mu_j$ is a birational morphism. Let $b'_j$ be the number of distinct 
even members of $\mathbf{p}'_j$. Then $b'_j = b -1$. This, in particular, means that 
$$\pi_1(O_{\mathbf{p}'_j}) \cong (\mathbf{Z}/2\mathbf{Z})^{\oplus b-1}.$$  
Let $X_{\mathbf{p}'_j} \to \bar{O}_{\mathbf{p}'_j}$ be the finite covering associated with 
the universal covering of $O_{\mathbf{p}'_j}$. Then $\mathfrak{n}_j + X_{\mathbf{p}'_j}$ 
is a $Q_j$-space; hence we get a $(\mathbf{Z}/2\mathbf{Z})^{\oplus b-1}$-covering map   
$$\pi'_j: Sp(2n) \times^{Q_j}(\mathfrak{n}_j + X_{\mathbf{p}'_j}) \to 
Sp(2n) \times^{Q_j}(\mathfrak{n}_j + \bar{O}_{\mathbf{p}'_j}) $$ 
Let $X_j$ be the Stein factorization of $\mu_j \circ \pi'_j$. Then we have a commutative diagram    

\begin{equation} 
\begin{CD} 
Sp(2n) \times^{Q_j}(\mathfrak{n}_j + X_{\mathbf{p}'_j}) @>{\mu'_j}>> X_j \\ 
@V{\pi'_j}VV @V{\pi_j}VV \\ 
Sp(2n) \times^{Q_j}(\mathfrak{n}_j + \bar{O}_{\mathbf{p}'_j}) @>{\mu_j}>> \bar{O}_{\mathbf p}      
\end{CD} 
\end{equation}

By definition $\pi_j$ is a $(\mathbf{Z}/2\mathbf{Z})^{\oplus b-1}$-covering and 
$X_j$ factorizes $\pi$ as $$\pi: X \stackrel{\rho_j}\to X_j \stackrel{\pi_j}\to \bar{O}_{\mathbf p},$$ where $\rho_j$ is a $\mathbf{Z}/2\mathbf{Z}$-covering. 

\begin{Claim}\label{Claim (2.3.3)} (1) $\mu_j$ is a crepant resolution around $O_{\bar{\mathbf p}_{j}} \subset 
\bar{O}_{\mathbf p}$. In other words, $\mu_j^{-1}(S_j)$ is the minimal resolution of the $A_1$-surface singularity $S_j$.  
On the other hand, $\mu_j$ is an isomorphism over open neighborhoods of $O_{\bar{\mathbf p}_{1}}$, ..., $O_{\bar{\mathbf p}_{j-1}}$, $O_{\bar{\mathbf p}_{{j+1}}}$, ..., 
$O_{\bar{\mathbf p}_{k}}$. In other words, the maps $\mu_j^{-1}(S_1) \to S_1$, ..., 
$\mu_j^{-1}(S_{j-1}) \to S_{j-1}$, $\mu_j^{-1}(S_{j+1}) \to S_{j+1}$, ..., $\mu_j^{-1}(S_k) \to S_k$ 
are all isomorphism. 

(2) $\pi_j^{-1}(S_j)$ is a disjoint union of $2^{b-1}$ copies of $S_j$. $\pi^{-1}(S_j)$ is 
a disjoint union of $2^{b-1}$ copies of $(\mathbf{C}^2, 0)$. In other words, $\pi_j$ is an 
etale cover over an open neighborhood of $O_{\bar{\mathbf p}_j} \subset \bar{O}_{\mathbf p}$, and $\rho_j$ is a ramified double covering over an open neighborhood of $\pi_j^{-1}(O_{\bar{\mathbf p}_j}) \subset 
X_j$.   \end{Claim}

Notice that Claim \ref{Claim (2.3.3)}, (2) implies Proposition \ref{Proposition (2.3)}.
 
{\em Proof}. (1) The gap members of $\mathbf{p}'_j$ are $$i_1 -2, \: ..., \: i_{j-1}-2,\:  i_{j+1}, ...,  
i_k.$$ For the later convenience we put $$i'_1 := i_1 -2,\: ..., i'_{j-1} := i_{j-1} -2, \: 
i'_{j+1} := i_{j+1}, ..., i'_k := i_k.$$ Corresponding to these gap members, we get 
nilpotent orbits $O_{(\overline{\mathbf{p}'_j})_1}$, ..., $O_{(\overline{\mathbf{p}'_j})_{j-1}}$, 
$O_{(\overline{\mathbf{p}'_j})_{j+1}}$, ..., $O_{(\overline{\mathbf{p}'_j})_k}$ in $\bar{O}_{\mathbf{p}'_j}$. These are irreducible components of $\mathrm{Sing}(\bar{O}_{\mathbf{p}'_j})$ which have codimension 2 in $\bar{O}_{\mathbf{p}'_j}$.  
For each $1 \le l \le k$ $(l \ne j)$, we have a natural embedding 
$$Sp(2n) \times^{Q_j}(\mathfrak{n}_j + \bar{O}_{(\overline{\mathbf{p}'_j})_l}) \subset 
Sp(2n) \times^{Q_j}(\mathfrak{n}_j + \bar{O}_{\mathbf{p}'_j}).$$
One can check that $$\mu_j(\:\: Sp(2n) \times^{Q_j}(\mathfrak{n}_j + \bar{O}_{(\overline{\mathbf{p}'_j})_l})\:\:) = \bar{O}_{\bar{\mathbf p}_l}.$$   
From this fact we see that there is an irreducible component of $\mathrm{Sing}(\:\: 
Sp(2n) \times^{Q_j}(\mathfrak{n}_j + \bar{O}_{\mathbf{p}'_j})\:\:)$ which dominates 
$\bar{O}_{\bar{\mathbf p}_l}$ for each $l (\ne j)$, but there is no irreducible components 
of $\mathrm{Sing}(\:\: Sp(2n) \times^{Q_j}(\mathfrak{n}_j + \bar{O}_{\mathbf{p}'_j})\:\:)$ which 
dominates $\bar{O}_{\bar{\mathbf p}_j}$. 

Since $\mu_j$ is a crepant partial 
resolution of $\bar{O}_{\mathbf p}$ and $\bar{O}_{\mathbf p}$ has $A_1$-surface singularity 
along $O_{\bar{\mathbf p}_l}$, this means that $\mu_j$ is an isomorphism over an open neighborhood of $O_{\bar{\mathbf p}_l} \subset \bar{O}_{\mathbf p}$ for $l \ne j$, but  
$\mu_j$ is a crepant resolution around $O_{\bar{\mathbf p}_j} \subset \bar{O}_{\mathbf p}$. 

(2) Since $\mu_j$ is a crepant partial resolution and $\pi'_j$ is etale in codimension 1, 
the partial resolution $\mu'_j$ is a crepant partial resolution. By (1) there is a  
$\mu_j$-exceptional divisor $E$ of $Sp(2n) \times^{Q_j}(\mathfrak{n}_j + \bar{O}_{\mathbf{p}'_j})$ 
which dominates $\bar{O}_{\bar{\mathbf p}_j}$. Then $(\pi'_j)^{-1}(E)$ is a $\mu'_j$-exceptional 
divisor of $Sp(2n) \times^{Q_j}(\mathfrak{n}_j + X_{\mathbf{p}'_j})$ which dominates 
$\pi_j^{-1}(\bar{O}_{\bar{\mathbf p}_j})$. If $\pi_j$ is ramified over ${O}_{\bar{\mathbf p}_j}$, 
then $X_j$ is smooth along $\pi_j^{-1}({O}_{\bar{\mathbf p}_j})$. This contradicts that 
$\mu'_j$ is a crepant partial resolution. Hence $\pi$ is unramified over $\bar{O}_{\bar{\mathbf p}_j}$, which is nothing but the first statement of (2). 
Next suppose that the second statement of (2) does not hold. Then $\rho_j$ is etale over 
an open neighborhood of $\pi_j^{-1}({O}_{\bar{\mathbf p}_j}) \subset X_j$. Let $(\bar{O}_{\mathbf p})^0$ be an open set obtained from $\bar{O}_{\mathbf p}$ by excluding all irreducible 
components of $\mathrm{Sing}(\bar{O}_{\mathbf p})$ different from $\bar{O}_{\bar{\mathbf p}_j}$. We put $X_j^0 := \pi_j^{-1}((\bar{O}_{\mathbf p})^0)$. 
Then $\rho_j^{-1}(X_j^0) \to X_j^0$ is an etale cover. On the other hand, let $X_{\mathbf{p}'_j}^0$ be the universal covering of $O_{\mathbf{p}'_j}$. Then 
$Sp(2n) \times^{Q_j}(\mathfrak{n}_j + X_{\mathbf{p}'_j}^0)$ is simply connected. 
In fact, we have an exact sequence 
$$\pi_1(\mathfrak{n}_j + X_{\mathbf{p}'_j}^0) \to \pi_1(Sp(2n) \times^{Q_j}(\mathfrak{n}_j + X_{\mathbf{p}'_j}^0)) \to \pi_1(Sp(2n)/Q_j) \to 1$$
Since $\pi_1(\mathfrak{n}_j + X_{\mathbf{p}'_j}^0) = \{1\}$ and $\pi_1(Sp(2n)/Q_j) = \{1\}$, 
we have the result. 

For $l \ne j$, $\mu_j$ is an isomorphism over an open neighborhood of $O_{\bar{\mathbf p}_l} \subset \bar{O}_{\mathbf p}$ by (1). By Zariski's Main Theorem $\mu'_j$ is also an isomorphism over an open neighborhood of $\pi_j^{-1}(O_{\bar{\mathbf p}_l}) \subset X_j$. 
Moreover, $\mu'_j$ is a crepant resolution around $\pi_j^{-1}(O_{\mathbf{p}'_j}) \subset X_j$.  
Write 
$$\bar{O}_{\mathbf p} = O_{\mathbf{p}} \sqcup O_{\bar{\mathbf p}_1} \sqcup ... \sqcup O_{\bar{\mathbf p}_k} \sqcup F,$$ 
where $F$ is the union of all nilpotent orbits in $\bar{O}_{\mathbf p}$ with codimension $\geq 4$. 
Then 
$$Sp(2n) \times^{Q_j}(\mathfrak{n}_j + X_{\mathbf{p}'_j}) = (\pi_j \circ \mu'_j)^{-1}(O_{\mathbf{p}}) \sqcup (\pi_j \circ \mu'_j)^{-1}(O_{\bar{\mathbf p}_1}) \sqcup ... \sqcup (\pi_j \circ \mu'_j)^{-1}(O_{\bar{\mathbf p}_k}) \sqcup (\pi_j \circ \mu'_j)^{-1}(F).$$
For $l \ne j$, we see that $(\pi_j \circ \mu'_j)^{-1}(O_{\bar{\mathbf p}_l})$ has codimension $2$ in  
$Sp(2n) \times^{Q_j}(\mathfrak{n}_j + X_{\mathbf{p}'_j})$. Moreover, since $\mathrm{Codim}_{X_j}\pi_j^{-1}(F) \geq 4$, 
we also see that  $(\pi_j \circ \mu'_j)^{-1}(F)$ has codimension $\geq 2$ by Corollary \ref{Corollary (0.2)}. 
This means that $(\mu'_j)^{-1}(\pi_j^{-1}(O_{\mathbf p} \cup O_{\bar{\mathbf p}_j}))$ is obtained from 
$Sp(2n) \times^{Q_j}(\mathfrak{n}_j + X_{\mathbf{p}'_j})$ by removing a closed subset of codimension $\geq 2$.  
Since $(\mu'_j)^{-1}(\pi_j^{-1}(O_{\mathbf p} \cup O_{\bar{\mathbf p}_j}))$ is smooth, it is contained in $Sp(2n) \times^{Q_j}(\mathfrak{n}_j + X_{\mathbf{p}'_j}^{\mathrm{reg}})$.  
This implies that $(\mu'_j)^{-1}(X_j^0)$ is obtained from a smooth variety $Sp(2n) \times^{Q_j}(\mathfrak{n}_j + X_{\mathbf{p}'_j}^{\mathrm{reg}})$ by removing a closed subset of codimension 
$\geq 2$. There is a surjection map 
$$\pi_1(Sp(2n) \times^{Q_j}(\mathfrak{n}_j + X_{\mathbf{p}'_j}^0)) \to \pi_1(Sp(2n) \times^{Q_j}(\mathfrak{n}_j + X_{\mathbf{p}'_j}^{\mathrm{reg}})).$$  As already shown, the left hand side is trivial; hence, the right hand side is also trivial.   
Now we have $$\pi_1((\mu'_j)^{-1}(X_j^0)) \cong \pi_1(Sp(2n) \times^{Q_j}(\mathfrak{n}_j + X_{\mathbf{p}'_j}^{\mathrm{reg}}) = \{1\}.$$ 
Since $\mu'_j$ is birational, $\pi_1(X_j^0) = \{1\}$. 
This contradicts that $\rho_j^{-1}(X_j^0)$ is a (connected) etale cover of $X_j^0$ of degree 2. 
$\square$  \vspace{0.2cm}

We here consider an additional condition for $\mathbf{p}$: 

(iii) $r_i \ne 2$ for each even $i$. 

\begin{Cor}\label{Corollary (2.4)} Let $\mathbf{p}$ be a partition satisfying (i), (ii) and (iii). 
Let $\pi: X \to \bar{O}_{\mathbf p}$ be the finite covering associated with the universal 
covering of $O_{\mathbf p}$. Then $X$ has only {\bf Q}-factorial terminal singularities. \end{Cor} 

{\em Proof}. By the condition (iii), $X$ is {\bf Q}-factorial by Proposition \ref{Proposition (2.1)}, (2). 
Then the result follows from Proposition \ref{Proposition (2.3)}. $\square$ \vspace{0.2cm}

{\bf Construction of a Q-factorial terminalization}:  
Let $\bar{O}_{\mathbf p}$ be an arbitrary nilpotent orbit closure of $sp(2n)$, and 
let $X \to \bar{O}_{\mathbf p}$ be the finite covering associated with the universal 
covering of $O_{\mathbf p}$. 
We shall construct explicitly a {\bf Q}-factorial terminalization of $X$. Since $\mathbf{p}$ is a Jordan type of a nilpotent orbit, $\mathbf{p}$ satisfies the condition (i). Let $b$ be the number of distinct even members of $\mathbf{p}$.

By using the inductions of type (I) repeatedly for $\mathbf{p}$, we can finally find a parabolic subgroup 
$Q$ of $Sp(2n)$ and a nilpotent orbit $O_{\mathbf{p}'}$ of a Levi part 
$\mathfrak{l}$ of $\mathfrak{q}$ such that 

(a) $O_{\mathbf p} = \mathrm{Ind}^{\mathfrak{sp}(2n)}_{\mathfrak l}(O_{\mathbf{p}'}).$

(b) $\mathbf{p}'$ satisfies the condition (ii), and $b' = b$ (where $b'$ is the number of 
distinct even members of $\mathbf{p}'$).

Notice that $\mathfrak{l}$ is a direct sum of a simple Lie algebra $sp(2n')$, some  
simple Lie algebras of type $A$ and the center. $O_{\mathbf{p}'}$ is a nilpotent orbit of 
$sp(2n')$. \vspace{0.2cm}

For the partition $\mathbf{p}'$ thus obtained, we let $e$ be the number of even members 
$i$ of $\mathbf{p}'$ such that $r_i = 2$.
By using the inductions of type (II) repeatedly for $\mathbf{p}'$, we 
can finally find a parabolic subgroup 
$Q'$ of $Sp(2n')$ and a nilpotent orbit $O_{\mathbf{p}''}$ of a Levi part 
$\mathfrak{l}'$ of $\mathfrak{q}'$ such that 

(a') $O_{{\mathbf p}'} = \mathrm{Ind}^{\mathfrak{sp}(2n')}_{{\mathfrak l}'}(O_{\mathbf{p}''}),$

(b') $\mathbf{p}''$ satisfies the conditions (ii), (iii),  and $b'' = b' - e$ (where $b''$ is the number of 
distinct even members of $\mathbf{p}''$).  \vspace{0.2cm}

All together, we get a generalized Springer map 
$$\mu'': Sp(2n) \times^{Q''}(\mathfrak{n}'' + \bar{O}_{{\mathbf p}''}) \to \bar{O}_{\mathbf p},$$ 
where $\mu''$ is generically finite of degree $2^{b - b''}$. 
Let $\pi'': X'' \to \bar{O}_{{\mathbf p}''}$ be the finite covering associated with the universal 
covering of $O_{{\mathbf p}''}$. 
By Corollary \ref{Corollary (2.4)} $X''$ has only {\bf Q}-factorial terminal singularities.  
Note that $\mathrm{deg}(\pi'') = 2^{b''}$. There is a finite cover 
$$\pi'': Sp(2n) \times^{Q''}(\mathfrak{n}'' + X'') \to Sp(2n) \times^{Q''}(\mathfrak{n}'' + \bar{O}_{{\mathbf p}''})$$ of degree $2^{b''}$. Then the Stein factorization 
of $\mu'' \circ \pi''$ coincides with $X$. We have a commutative diagram  
 
\begin{equation} 
\begin{CD} 
Sp(2n) \times^{Q''}(\mathfrak{n}'' + X'') @>>> X \\ 
@V{\pi''}VV @V{\pi}VV \\ 
Sp(2n) \times^{Q''}(\mathfrak{n}'' + \bar{O}_{{\mathbf p}''}) @>{\mu''}>> \bar{O}_{\mathbf p}      
\end{CD} 
\end{equation}

The map $Sp(2n) \times^{Q''}(\mathfrak{n}'' + X'') \to X$ here obtained is a {\bf Q}-factorial terminalization of $X$. \vspace{0.2cm}

\begin{Exams}\label{Examples (2.5)} \end{Exams} (1) Let us consider the nilpotent orbit $O_{[6^2,4^2]} \subset 
sp(20)$ and let $\pi: X \to \bar{O}_{[6^2,4^2]}$ be the finite covering associated with 
the universal covering of $O_{[6^2,4^2]}$. We have $\mathrm{deg}(\pi) = 4$  by 
Proposition \ref{Proposition (2.1)}, (1). We shall construct a {\bf Q}-factorial terminalization of $X$ and 
we shall show that it is actually a crepant resolution. 

Let $Q_{4,12,4}$ be a parabolic subgroup of $Sp(20)$ with flag type 
$(4,12,4)$. Let $\mathfrak{q}_{4,12,4} = \mathfrak{n}_{4,12,4} \oplus \mathfrak{l}_{4,12,4}$ 
be a Levi decomposition. 
Then $\mathfrak{l}_{4,12,4} = sp(12) \oplus \mathfrak{gl}(4)$ 
and $O_{[6^2,4^2]}$ is induced from the nilpotent orbit $O_{[4^2, 2^2]} \subset sp(12)$. 
This is a type I induction, and the generalized Springer map $$  
Sp(20) \times^{Q_{4,12,4}}(\mathfrak{n}_{4,12,4} + \bar{O}_{[4^2, 2^2]}) \to 
\bar{O}_{[6^2,4^2]}$$ is a birational map. Let $Q_{1, 10, 1}$ be a parabolic subgroup of $Sp(12)$  
with the Levi decomposition $\mathfrak{q}_{1,10,1} = \mathfrak{n}_{1,10,1} \oplus \mathfrak{l}_{1,10,1}$. Then $\mathfrak{l}_{1,10,1} = sp(10) \oplus \mathfrak{gl}(1)$ and 
$O_{[4^2,2^2]}$ is induced from $O_{[3^2,2^2]} \subset sp(10)$. This is a type II induction. 
Hence we get a generically finite map of degree 2 $$Sp(20) \times^{Q_{4,1,10,1,4}}(\mathfrak{n}_{4,1,10,1,4} + 
\bar{O}_{[3^2,2^2]}) \to Sp(20) \times^{Q_{4,12,4}}(\mathfrak{n}_{4,12,4} + \bar{O}_{[4^2, 2^2]}).$$
Let $Q_{3,4,3}$ be a parabolic subgroup of $Sp(10)$ with flag type $(3,4,3)$ with a Levi 
decomposition $\mathfrak{q}_{3,4,3} = \mathfrak{n}_{3,4,3} \oplus \mathfrak{l}_{3,4,3}$. 
Then $\mathfrak{l}_{3,4,3} = sp(4) \oplus \mathfrak{gl}(3)$ and $O_{[3^2,2^2]}$ 
is induced from $O_{[1^4]} \subset sp(4)$. This is a type II induction, and we get a generically 
finite map of degree 2 
$$Sp(20) \times^{Q_{4,1,3,4,3,1,4}} \mathfrak{n}_{4,1,3,4,3,1,4} \to 
Sp(20) \times^{Q_{4,1,10,1,4}}(\mathfrak{n}_{4,1,10,1,4} + 
\bar{O}_{[3^2,2^2]}).$$ We can illustrate the induction step above by 
$$([1^4], sp(4)) \stackrel{\mathrm{Type II}}\to (O_{[3^2,2^2]}, sp(10)) \stackrel{\mathrm{Type II}}\to 
(O_{[4^2, 2^2]}, sp(12)) \stackrel{\mathrm{Type I}}\to (O_{[6^2,4^2]}, sp(20)).$$
Composing these 3 maps together, we have a generically 
finite map of degree 4 
$$Sp(20) \times^{Q_{4,1,3,4,3,1,4}} \mathfrak{n}_{4,1,3,4,3,1,4} \to \bar{O}_{[6^2,4^2]}.$$ 
This map factors through $X$ and $Sp(20) \times^{Q_{4,1,3,4,3,1,4}} \mathfrak{n}_{4,1,3,4,3,1,4}$ 
gives a crepant resolution of $X$. 

(2) Let $\pi: X \to \bar{O}_{[8, 5^2, 4, 3^2]}$ be the finite covering associated with the universal covering of $O_{[8, 5^2, 4, 3^2]}$. By Proposition \ref{Proposition (2.1)}, (1) we have $\mathrm{deg}(\pi) = 4$. 
We can take the following inductions
$$(O_{[4,3^2,2,1^2]}, sp(14)) \stackrel{\mathrm{Type I}}\to (O_{[6,5^2,4,3^2]}, sp(26)) \stackrel{\mathrm{Type I}}\to 
(O_{[8, 5^2, 4, 3^2]}, sp(28)).$$ Note that in each step the number of distinct even members 
does not change, and the partition $[4,3^2,2,1^2]$ satisfies the conditions (i), (ii) and (iii).  
Let $X_{[4,3^2,2,1^2]} \to \bar{O}_{[4,3^2,2,1^2]}$ be the finite covering associated 
with the universal covering of $O_{[4,3^2,2,1^2]}$. By Corollary \ref{Corollary (2.4)}, $X_{[4,3^2,2,1^2]}$ has 
only {\bf Q}-factorial terminal singularities. Then $Sp(28) \times^{Q_{1,6,14,6,1}}(\mathfrak{n}_{1,6,14,6,1} + X_{[4,3^2,2,1^2]})$ gives a {\bf Q}-factorial terminalization of 
$X$. $\square$

\section{$\mathfrak{g} = so(m)$}    

We write a partition $\mathbf{p}$ of $m$ as $[d^{r_d}, (d-1)^{r_{d-1}}, ..., 2^{r_2}, 1^{r_1}]$ with $r_d \ne 0$. Other $r_i$ may be possibly zero; in such a case $i$ does not appear in the partition. If $r_i > 0$, then we call $i$ a member of the partition.
Let $O$ be a nilpotent orbit of $\mathfrak{g}$. Then its Jordan type $\mathbf{p}$ is a partition 
of $m$  
such that all even members have even multiplicities (cf. \cite{C-M}, \S 5). When $\mathbf{p}$ consists of only even members, we call $\mathbf{p}$ {\em very even}. If $\mathbf{p}$ is not very even, then the orbit $O$ is uniquely determined by the Jordan type $\mathbf{p}$, and 
we denote by $O_{\mathbf{p}}$ the nilpotent orbit. On the other hand, if $\mathbf{p}$ is very 
even, there are two nilpotent orbits with Jordan type $\mathbf{p}$. The two orbits are conjugate to each other by an element of $O(m) \setminus SO(m)$. When we want to distinguish them, we denote them by $O^+_{\mathbf p}$, $O^{-}_{\mathbf p}$, but we usually 
denote by $O_{\mathbf{p}}$ one of them. A partition $\mathbf{p}$ is called {\em rather odd} 
if all odd members have multiplicity $1$. Note that a very even partition is rather odd. Let $a$ 
be the number of distinct odd members of $\mathbf{p}$. When $\mathfrak{g} = so(2n+1)$, we always have $a > 0$, but when $\mathfrak{g} = so(2n)$, we may possibly have $a = 0$. 

\begin{Prop}\label{Proposition (3.1)}  
(1) If $\mathbf{p}$ is not rather odd, then
$$\pi_1(O_{\mathbf p}) \cong (\mathbf{Z}/2\mathbf{Z})^{\oplus \mathrm{max}(a-1,0)}.$$ 
If $\mathbf{p}$ is rather odd, then there is a short exact sequence 
$$ 1 \to \mathbf{Z}/2\mathbf{Z} \to \pi_1(O_{\mathbf p}) \to (\mathbf{Z}/2\mathbf{Z})^{\oplus \mathrm{max}(a-1,0)} \to 1$$ 
so that $\mathbf{Z}/2\mathbf{Z}$ is contained in the center of 
$\pi_1(O_{\mathbf p})$. 

(2) Assume that $r_i \ne 2$ for all odd members $i$ of $\mathbf{p}$. 
Then $X$ is {\bf Q}-factorial for any etale covering $X^0 \to O_{\mathbf p}$. 
\end{Prop}

{\em Proof}. (1) Put $G = Spin(m)$. There is a double covering 
$\rho_m: G \to SO(m)$ and $\mathfrak{g} = so(m)$.  
Take an element $x$ from 
$O_{\mathbf p}$ and take an $sl(2)$-triple $\phi$ in $\mathfrak{g}$ containing $x$.  
By [C-M, Theorem (6.1.3)] $SO(m)^{\phi}$ is isomorphic to  
$$S(\prod_{i: \mathrm{even}}Sp(r_i)^{\times i}_{\Delta} \times 
\prod_{i: \mathrm{odd}}O(r_i)^{\times i}_{\Delta}) := $$
$$\{(\prod_{i: \mathrm{even}} A_{r_i}^{\times i}, \prod_{i: \mathrm{odd}} B_{r_i}^{\times i}) 
\in \prod_{i: \mathrm{even}}Sp(r_i)^{\times i}_{\Delta} \times 
\prod_{i: \mathrm{odd}}O(r_i)^{\times i}_{\Delta} \: \vert \: 
\prod_{i: \mathrm{odd}} \mathrm{det}(B_{r_i})^i = 1\}.$$ Then $G^{\phi} = \rho_m^{-1}(SO(m)^{\phi})$; hence $G^{\phi}$ is a double cover of   
$S(\prod_{i: \mathrm{even}}Sp(r_i)^{\times i}_{\Delta} \times 
\prod_{i: \mathrm{odd}}O(r_i)^{\times i}_{\Delta})$.  
When $\mathbf{p}$ has only even members, $S(\prod_{i: \mathrm{even}}Sp(r_i)^{\times i}_{\Delta} \times \prod_{i: \mathrm{odd}}O(r_i)^{\times i}_{\Delta}) = 
\prod_{i: \mathrm{even}}Sp(r_i)^{\times i}_{\Delta}$, which is connected. 
When $\mathbf{p}$ has some odd members, $S(\prod_{i: \mathrm{even}}Sp(r_i)^{\times i}_{\Delta} \times \prod_{i: \mathrm{odd}}O(r_i)^{\times i}_{\Delta})$ has $2^{a-1}$ connected 
components. Therefore, in any case, $S(\prod_{i: \mathrm{even}}Sp(r_i)^{\times i}_{\Delta} \times \prod_{i: \mathrm{odd}}O(r_i)^{\times i}_{\Delta})$ has $2^{\mathrm{max}(a-1,0)}$ 
connected components. 

If $\mathbf{p}$ is not rather odd, then the identity component 
of $S(\prod_{i: \mathrm{even}}Sp(r_i)^{\times i}_{\Delta} \times \prod_{i: \mathrm{odd}}O(r_i)^{\times i}_{\Delta})$ is  
$\prod_{i: \mathrm{even}}Sp(r_i)^{\times i}_{\Delta} \times \prod_{i: \mathrm{odd}}SO(r_i)^{\times i}_{\Delta}$. Note that $\pi_1(SO(r_i)) = \mathbf{Z}/2\mathbf{Z}$ for $r_i >2$ and $\pi_1(SO(2)) = \mathbf{Z}$. 
For an odd $i$ with $r_i \geq 2$, there is a unique surjective homomorphism 
$\phi_i: \pi_1(SO(r_i)^{\times i}_{\Delta}) \to \mathbf{Z}/2\mathbf{Z}$. 
We then have a surjection $$ \sum \phi_i : 
\prod_{i: odd, \: r_i \geq 2} \pi_1(SO(r_i)^{\times i}_{\Delta}) \to \mathbf{Z}/2\mathbf{Z}.$$
The left hand side is identified with $\pi_1(\prod_{i: \mathrm{even}}Sp(r_i)^{\times i}_{\Delta} \times \prod_{i: \mathrm{odd}}SO(r_i)^{\times i}_{\Delta})$. Therefore 
it determines a connected etale double covering of $\prod_{i: \mathrm{even}}Sp(r_i)^{\times i}_{\Delta} \times \prod_{i: \mathrm{odd}}SO(r_i)^{\times i}_{\Delta}$. One can check that $\rho_m^{-1}(\prod_{i: \mathrm{even}}Sp(r_i)^{\times i}_{\Delta} \times \prod_{i: \mathrm{odd}}SO(r_i)^{\times i}_{\Delta})$ is such an etale 
covering. 
Hence $\rho_m^{-1}(S(\prod_{i: \mathrm{even}}Sp(r_i)^{\times i}_{\Delta} \times \prod_{i: \mathrm{odd}}O(r_i)^{\times i}_{\Delta}))$ has $2^{\mathrm{max}(a-1,0)}$ connected components. If $\mathbf{p}$ is rather odd, then each connected component 
of $S(\prod_{i: \mathrm{even}}Sp(r_i)^{\times i}_{\Delta} \times \prod_{i: \mathrm{odd}}O(r_i)^{\times i}_{\Delta})$ is isomorphic to $\prod_{i: \mathrm{even}}Sp(r_i)^{\times i}_{\Delta}$, which is simply connected. 
Hence $\rho_m^{-1}(S(\prod_{i: \mathrm{even}}Sp(r_i)^{\times i}_{\Delta} \times \prod_{i: \mathrm{odd}}O(r_i)^{\times i}_{\Delta}))$ has 
$2\cdot 2^{\mathrm{max}(a-1,0)}$ connected components.  \vspace{0.2cm}

(2) If $\mathbf{p}$ is rather odd, $(G^{\phi})^0$ is isomorphic to 
$\prod_{i: \mathrm{even}}Sp(r_i)^{\times i}_{\Delta}$.
Then  $\chi ((G^{\phi})^0) = 0$.
If $\mathbf{p}$ is not rather odd, then $(G^{\phi})^0$ is a double covering of $\prod_{i: \mathrm{even}}Sp(r_i)^{\times i}_{\Delta} \times \prod_{i: \mathrm{odd}}SO(r_i)^{\times i}_{\Delta}$. Note that $SO(r)$ is a simple Lie group except that $SO(2) \cong \mathbf{C}^*$ and $SO(4)$ is a semisimple Lie group of type $A_1 + A_1$. This means that, if $r_i \ne 2$ for all odd $i$, then  
$$\chi (\: \prod_{i: \mathrm{even}}Sp(r_i)^{\times i}_{\Delta} \times \prod_{i: \mathrm{odd}}SO(r_i)^{\times i}_{\Delta} \:) = 0.$$
Then the short exact sequence $$1 \to \mathbf{Z}/2\mathbf{Z} \to 
(G^{\phi})^0 \to \prod_{i: \mathrm{even}}Sp(r_i)^{\times i}_{\Delta} \times \prod_{i: \mathrm{odd}}SO(r_i)^{\times i}_{\Delta} \to 1$$ yields an exact sequence 
of character groups $$\chi (\: \prod_{i: \mathrm{even}}Sp(r_i)^{\times i}_{\Delta} \times \prod_{i: \mathrm{odd}}SO(r_i)^{\times i}_{\Delta})\: ) \to \chi ((G^{\phi})^0) \to \chi (\mathbf{Z}/2\mathbf{Z}).$$ The 1-st term is zero by the observation above, and 
the 3-rd term is a finite group. Hence $\chi ((G^{\phi})^0)$ is also finite. As a result, 
$\chi ((G^{\phi})^0)$ is finite in any case under the assumption (2).  The remainder  
is the same as in the proof of Proposition (1.1), (2).   $\square$ \vspace{0.2cm}

Let $\pi: X \to \bar{O}_{\mathbf p}$ be the finite covering associated with the universal covering of $O_{\mathbf p}$ and $\tau: Y \to \bar{O}_{\mathbf p}$ 
a finite covering determined by the surjection $\pi_1(O_{\mathbf p}) \to 
(\mathbf{Z}/2\mathbf{Z})^{\oplus \mathrm{max}(a(\mathbf p)-1, 0)}$ in Proposition \ref{Proposition (3.1)}, (1).

\begin{Prop}\label{Proposition (3.2)} The adjoint action of $SO(m)$
on $\bar{O}_{\mathbf p}$ lifts to an $SO(m)$-action on $Y$.  If $\mathbf{p}$ is not rather odd, then $X = Y$. If $\mathbf{p}$ is rather odd, then $X$ is a double cover of 
$Y$ and the $SO(m)$ action on $Y$ does not lift to an $SO(m)$-action on $X$.  \end{Prop}

{\em Proof}. Put $X^0 := \pi^{-1}(O_{\mathbf p})$ and $Y^0 := \tau^{-1}(O_{\mathbf p})$. 
Choose $x \in O_{\mathbf p}$. By the proof of Proposition \ref{Proposition (3.1)}, the double cover $\rho_m: Spin(m) \to 
SO(m)$ induces a double cover $(Spin(m)^x)^0 \to (SO(m)^x)^0$ when  
$\mathbf{p}$ is not rather odd, and induces an isomorphism $(Spin(m)^x)^0 \cong  (SO(m)^x)^0$ when $\mathbf{p}$ is rather odd. In any case, $Y^0 = SO(m)/(SO(m)^x)^0$. 
Hence $SO(m)$ naturally acts on $Y^0$. Since $\Gamma (Y, \mathcal{O}_Y) 
= \Gamma (Y^0, \mathcal{O}_{Y^0})$, $SO(m)$ acts on $Y$.  
The natural $SO(m)$ action on $SO(m)/SO(m)^x$ is nothing but the 
adjoint action of $SO(m)$ on $O_{\mathbf p}$. This means that the $SO(m)$-action 
on $Y^0$ is a lift of the adjoint $SO(m)$-action on $O_{\mathbf p}$. Therefore 
the $SO(m)$-action on $Y$ is a lift of the adjoint action of $SO(m)$ on 
$\bar{O}_{\mathbf p}$. 

Assume that $\mathbf{p}$ is not rather odd. Then $$X^0 = Spin(m)/(Spin(m)^x)^0 
= SO(m)/(SO(m)^x)^0 = Y^0.$$ Hence $X = Y$. 

Assume that $\mathbf{p}$ is rather odd. Then $\pi^0$ factorizes as $X^0 \stackrel{\rho^0}\to Y^0 \stackrel{\tau^0}\to O_{\mathbf p}$. 
Let us consider the composite $SO(m) \times X^0 \to Y^0$ of the map 
$SO(m) \times X^0 \stackrel{id \times \rho^0}\to SO(m) \times Y^0$ and 
the map $SO(m) \times Y^0 \to Y^0$ determined by the $SO(m)$-action on $Y^0$.
The map $SO(m) \times X^0 \to Y^0$ lifts to a map to $X^0$ if and only 
if $\pi_1(SO(m) \times X^0) \to \pi_1(Y^0)$ is the zero map.   
Take a point $\tilde{x} \in X^0$ such that $\pi^0(\tilde{x}) = x$. Then the maps 
$$SO(m) \times \{\tilde{x}\} \to SO(m) \times X^0 \to Y^0$$ induces homomorphisms 
of fundamental groups 
$$\pi_1(SO(m) \times \{\tilde{x}\}) \to \pi_1(SO(m) \times X^0) \to \pi_1(Y^0) = \mathbf{Z}/2\mathbf{Z}.$$  
Since $\pi_1(X^0) = 1$, the first map is an isomorphism. Since  
$SO(m) \times \{\tilde{x}\}$ is a fibre bundle over $Y^0$ with a typical fiber 
$(SO(m)^x)^0$, the map $\pi_1(SO(m) \times \{\tilde{x}\}) \to \pi_1(Y^0)$ is a 
surjection. As a consequence, $\pi_1(SO(m) \times X^0) \to \pi_1(Y^0) = \mathbf{Z}/2\mathbf{Z}$ is a surjection. This means that the $SO(m)$-action 
on $Y^0$ does not lift to an $SO(m)$-action on $X^0$.  $\square$

\begin{Lem}\label{Lemma (3.3)} Assume that $\mathbf{p}$ is not very even. Then the 
adjoint action of $O(m)$ on $\bar{O}_{\mathbf p}$ lifts to an $O(m)$-action 
on $Y$. \end{Lem}

{\em Remark}. A lifting of an $O(m)$-action is not unique. 

{\em Proof}. We fix an element $x$ of $O_{\mathbf p}$ and an $sl(2)$-triple 
$\phi$ containing $x$. Then $O(m)^{\phi}$ is isomorphic to   
$$\prod_{i: \mathrm{even}}Sp(r_i)^{\times i}_{\Delta} \times 
\prod_{i: \mathrm{odd}}O(r_i)^{\times i}_{\Delta} := $$
$$\{(\prod_{i: \mathrm{even}} A_{r_i}^{\times i}, \prod_{i: \mathrm{odd}} B_{r_i}^{\times i}) 
\in \prod_{i: \mathrm{even}}Sp(r_i)^{\times i} \times 
\prod_{i: \mathrm{odd}}O(r_i)^{\times i}.$$  
Note that $$(O(m)^{\phi})^0 = 
\prod_{i: \mathrm{even}}Sp(r_i)^{\times i}_{\Delta} \times 
\prod_{i: \mathrm{odd}}SO(r_i)^{\times i}_{\Delta}.$$ 
Since $\mathbf{p}$ is not very even, there is an odd member $i_0$. 
We put $$H^{\phi} = \prod_{i: \mathrm{even}}Sp(r_i)^{\times i}_{\Delta} \times 
\prod_{i: \mathrm{odd} \ne i_0}SO(r_i)^{\times i}_{\Delta} \times O(r_{i_0})^{\times i_0}_{\Delta},$$ and $H^x  = U^x \cdot H^{\phi}$.  Then we have an isomorphism 
$$Y^0 := SO(m)/(SO(m)^x)^0 \cong O(m)/H^x.$$ The right hand side has a natural $O(m)$-action. This $O(m)$-action determines an $O(m)$-action on 
$Y$.  $\square$

Let us consider the following conditions for a partition $\mathbf{p}$ of $m$. 

(i) $r_i$ is even for each even $i$. 

(ii) $r_i \ne 0$ for every odd $i$.

\begin{Prop}\label{Proposition (3.4)} Let $\mathbf{p}$ be a partition of $m$ which is not rather odd. Assume that $\mathbf{p}$ satisfies the conditions  (i) and (ii). Let 
$X \to \bar{O}_{\mathbf p}$ be the finite covering associated with the universal covering of $O_{\mathbf p}$. Then $\mathrm{Codim}_X\mathrm{Sing}(X) \geq 4$. \end{Prop}

{\em Proof}. Notice that the conditions (i) and (ii) are replacements of the 
conditions (i) and (ii) in the previous section where the roles of 
odd and even members are reversed. In this sense, this proposition is an $SO(m)$-analogue 
of Proposition \ref{Proposition (2.3)}. By virtue of Proposition \ref{Proposition (3.2)}, this proposition is proved completely in the same way as 
Proposition \ref{Proposition (2.3)}.  $\square$ \vspace{0.2cm}

When $\mathbf{p}$ is rather odd, we encounter a different situation as the following 
example illustrates.   

\begin{Exam}\label{Example (3.5)} \end{Exam} In this example, we show that a usual induction step (cf. \S 2) does 
not work well for a rather odd partition $\mathbf{p}$. 
Put the $m \times m$ matrix 
$$J_m = \left(\begin{array}{ccccc} 
0 & 0 & 0 & ... & 1 \\
0 & 0 & ... & 1 & 0 \\ 
... & ... & ... & ... & ... \\ 
0 & 1 & ... & 0 & 0 \\ 
1 & 0 & ... & 0 & 0 
\end{array}\right). $$   
Then $$SO(m) = \{ A \in SL(m) \: \vert \: A^tJ A = J\}. $$ Fix positive integers 
$s_1, ..., s_k, q$ so that $m = 2\sum s_i + q$. Assume that $q \geq 3$. 
Let $Q'$ be a parabolic subgroup of $SO(m)$ fixing the isotropic flag of flag type $(s_1, ..., s_k,q,s_k, ..., s_1)$ 
$$0 \subset \langle e_1, ..., e_{s_1} \rangle \subset \langle e_1, ... e_{s_1 + s_2} \rangle \subset ... \subset \langle e_1, ..., e_{\sum_{1}^{k}s_i} \rangle \subset 
\langle e_1, ..., e_{\sum_{1}^{k}s_i + q} \rangle$$
$$\subset \langle e_1, ..., e_{\sum_{1}^{k}s_i + q +s_k} \rangle \subset ... \subset 
\langle e_1, ..., e_{\sum_{1}^{k}s_i + q +\sum_{1}^{k}s_i} \rangle = \mathbf{C}^m $$
One has a Levi decomposition $Q' = U' \cdot L'$ with 
$$ L' = \{ \left(\begin{array}{ccccccccc} 
A_1 & 0 & 0 & ... & ... & ... & ... & ... & 0 \\
0 & A_2 & 0 & ... & ... & ... & ... & ... & 0 \\ 
... & ... & ... & ... & ... & ... & ... & ... & ...\\
0 & ... & 0 & A_k & 0 & ... & ... & ... & 0\\
0 & ... & ... & 0 & B & 0 & ... & ... & 0 \\
0 & ... & ... & ... & 0 & A'_k & 0 & ... & 0\\
... & ... & ... & ... & ... & ... & ... & ... & ... \\ 
0 & ... & ... & ... & ... & ... & ... & A'_2 & 0 \\ 
0 & ... & ... & ... & ... & ... & ... & 0 & A_1' 
\end{array}\right) \: \vert \: A_i \in GL(s_i), \: A'_i = J_{s_i}(A_i^t)^{-1}J_{s_i}, \: 
B \in SO(q)\}.$$ In particular, $L' \cong \prod GL(s_i) \times SO(q)$.
Let $Q := \rho_m^{-1}(Q')$ and $L := \rho_m^{-1}(L')$ for $\rho_m: Spin(m) \to SO(m)$. 
The Lie algebra $\mathfrak{l}$ of $L$ is isomorphic to $\oplus \mathfrak{gl}(s_i) \oplus 
so(q)$. Let $\mathbf{p}'$ be a rather odd partition of $q$ such that 
$a(\mathbf{p}) = a(\mathbf{p}')$ and let $O_{\mathbf{p}'}$ be 
a nilpotent orbit of $so(q)$ with Jordan type $\mathbf{p}'$. 
Assume that $\mathbf{p}'$ is obtained from $\mathbf{p}$ by a succession of type I inductions and $O_{\mathbf{p}} = \mathrm{Ind}^{so(m)}_{\mathfrak l}(O_{\mathbf{p}'})$. 
For an odd $s_i$, we consider the map $SL(s_i) \times \mathbf{C}^* \to GL(s_i)$ 
determined by $(X, \lambda) \to \lambda^2X$. This is a cyclic covering of order 
$2s_i$. The cyclic group $\mu_{2s_i}$ acts on $SL(s_i) \times \mathbf{C}^*$ so that  
$(X, \lambda) \to (\zeta_{2s_i}^{-2}X, \zeta_{2s_i}\lambda)$ for a  
primitive $2s_i$-th root $\zeta_{2s_i}$ of unity. Let us consider 
a unique subgroup $\mu_{s_i}$ of $\mu_{2s_i}$ of order $s_i$. Then 
$SL(s_i) \times \mathbf{C}^*/\mu_{s_i} \to GL(s_i)$ is a double cover of 
$GL(s_i)$. For an even $s_i$, we consider the map  $SL(s_i) \times \mathbf{C}^* \to GL(s_i)$ determined by $(X, \lambda) \to \lambda X$. This is a cyclic covering of order 
$s_i$. The cyclic group $\mu_{s_i}$ acts on 
$SL(s_i) \times \mathbf{C}^*$ so that 
$(X, \lambda) \to (\zeta_{s_i}^{-1}X, \zeta_{s_i}\lambda)$ for a primitive 
$s_i$-th root $\zeta_{s_i}$ of unity. Let us consider 
a unique subgroup $\mu_{s_i/2}$ of $\mu_{s_i}$ of order $s_i/2$. Then 
$SL(s_i) \times \mathbf{C}^*/\mu_{s_i/2} \to GL(s_i)$ is a double cover of 
$GL(s_i)$.  Here we put $t_i = s_i$ when $s_i$ is odd, and put $t_i := s_i/2$ 
when $s_i$ is even.   
We then have a covering map 
$$(SL(s_1) \times \mathbf{C}^*)/\mu_{t_1} \times ... \times   
(SL(s_k) \times \mathbf{C}^*)/\mu_{t_k} \times Spin(q) \to GL(s_1) \times ... 
\times GL(s_k) \times SO(q).$$ The Galois group of this covering is 
$(\mathbf{Z}/2\mathbf{Z})^{\oplus k +1}$. 
Put $H := \mathrm{Ker}[(\mathbf{Z}/2\mathbf{Z})^{\oplus k +1} \stackrel{\sum}\to \mathbf{Z}/2\mathbf{Z}]$, where $\sum$ is defined by $(x_1, ..., x_{k+1}) \to \sum x_i$. 
\vspace{0.2cm}

\begin{Claim}\label{Claim (3.5.1)}  
$$L =  \{(SL(s_1) \times \mathbf{C}^*)/\mu_{t_1} \times ... \times   
(SL(s_k) \times \mathbf{C}^*)/\mu_{t_k} \times Spin(q)\}\: /\: H.$$ \end{Claim}

{\em Proof}. We first prove that $\rho_m^{-1}(GL(s_i)) = (SL(s_i) \times \mathbf{C}^*)/\mu_{t_i}$ and $\rho_m^{-1}(SO(q)) = Spin(q)$. 
For each $1 \le j \le s_1 + ... + s_k$, we consider the non-degenerate quadratic 
subspace $V_j := \langle e_j, \: e_{m+1-j} \rangle $ of $\mathbf{C}^m$. Then $SO(V_j)$ 
is a subgroup of $SO(m)$ and $\rho_m^{-1}(SO(V_j)) = Spin(V_j)$ (cf. \cite{F-H}, (20.31)). 
Note that $GL(s_i)$ contains some $SO(V_j)$. 
Since $\rho_m^{-1}(SO(V_j))$ is connected, we see that $\rho_m^{-1}(GL(s_i))$ is a connected double cover of $GL(s_i)$. Since $\pi_1(GL(s_i)) = \mathbf{Z}$, we have a unique surjective homomorphism $\pi_1(GL(s_i)) \to \mathbf{Z}/2\mathbf{Z}$. 
This means that $\rho_m^{-1}(GL(s_i)) = (SL(s_i) \times \mathbf{C}^*)/\mu_{t_i}$. 
Since the $q$ dimensional subspace $\langle e_{\sum_{1}^{k}s_i + 1}, ..., e_{\sum_{1}^{k}s_i +q} \rangle$ is a non-degenerate quadratic space, we have 
$\rho_m^{-1}(SO(q)) = Spin(q)$ by the same reason as above.   
  
We then have Cartesian diagrams 
\begin{equation} 
\begin{CD} 
(SL(s_i) \times \mathbf{C}^*)/\mu_{t_i} @>{\tilde{\iota}_i}>>  L \\ 
@VVV @VVV \\ 
GL(s_i) @>{\iota_i}>> L'      
\end{CD} 
\end{equation}
and 
\begin{equation} 
\begin{CD} 
Spin(q) @>{\tilde{\iota}_q}>>  L \\ 
@VVV @VVV \\ 
SO(q) @>{\iota_q}>> L'
\end{CD} 
\end{equation}
Here $\iota_i$ and $\iota_q$ are natural injections and all vertical maps 
are induced by $\rho_m$. 

By using the group structure of $L$, we have a commutative 
diagram 
\begin{equation} 
\begin{CD} 
\prod (SL(s_i) \times \mathbf{C}^*)/\mu_{t_i}  \times Spin(q) 
@>{\tilde{\iota}_1 \cdot ... \cdot \tilde{\iota}_k \cdot \tilde{\iota}_q}>>  L \\ 
@VVV @VVV \\ 
\prod GL(s_i) \times SO(q) @>{\cong}>> L'      
\end{CD} 
\end{equation}
The covering group of the left vertical map is $(\mathbf{Z}/2\mathbf{Z})^{\oplus k + 1}$. 
The etale covering $L \to L'$ is determined by a certain surjection 
$\phi: (\mathbf{Z}/2\mathbf{Z})^{\oplus k + 1} \to \mathbf{Z}/2\mathbf{Z}$. 
Let $\tau_i : \mathbf{Z}/2\mathbf{Z} \to (\mathbf{Z}/2\mathbf{Z})^{\oplus k + 1}$ be 
the inclusion map into the $i$-th factor. By the previous Cartesian diagrams, the 
composite $\phi \circ \tau_i$ must be isomorphisms for all $1 \le i \le k+1$. 
This means that $\phi$ is defined by $(x_1, ..., x_{k+1}) \to \sum x_i$. $\square$   \vspace{0.2cm} 

Let $\pi': X' \to \bar{O}_{{\mathbf p}'}$ be the finite covering associated with 
the universal covering of $O_{{\mathbf p}'}$. Let $\tau': Y' \to \bar{O}_{{\mathbf p}'}$ 
be a finite covering determined by the surjection $\pi_1(O_{{\mathbf p}'}) \to 
(\mathbf{Z}/2\mathbf{Z})^{\oplus \mathrm{max}(a(\mathbf{p})-1, 0)}$ in 
Proposition \ref{Proposition (3.1)}, (1).  Put $(Y')^0 := \tau'^{-1}(O_{{\mathbf p}'})$ and $(X')^0 := \pi'^{-1}(O_{{\mathbf p}'})$. The adjoint action of $SO(q)$ on ${O}_{{\mathbf p}'}$ lifts to 
an $SO(q)$-action on $(Y')^0$. Since $L' = \prod GL(s_i) \times SO(q)$, $L'$ acts on 
$(Y')^0$ (and on $Y'$) by the projection map $L' \to SO(q)$.  By the surjection $L \to L'$, $L$ also acts on $(Y')^0$ (and on $Y'$).  

We prove that this $L$-action never lifts to an $L$-action on $(X')^0$. Let 
$p: L \to SO(q)$ be the composite of the map $L \to L'$ and the projection 
map $L' \to SO(q)$.  
Take a point $y_0 \in (Y')^0$ and define a map $L \to (Y')^0$ by 
$g \to g\cdot y_0$. By definition this map factorizes as $L \stackrel{p}\to 
SO(q) \to (Y')^0$. It induces homomorphisms 
$$\pi_1(L) \stackrel{p_*}\to \pi_1(SO(q)) \to \pi_1((Y')^0).$$ 
If the $L$-action on $(Y')^0$ lifts to $(X')^0$, then the composite of these 
homomorphisms is the zero map (cf. the proof of Proposition \ref{Proposition (3.2)}). 
By Claim \ref{Claim (3.5.1)}, the map $p$ has connected fibers. This means that $p_*$ is 
a surjection. Hence the map $\pi_1(SO(q)) \to \pi_1((Y')^0)$ must be zero. 
But this map is not the zero map because the $SO(q)$-action on $(Y')^0$ 
does not lift to $(X')^0$ by Proposition \ref{Proposition (3.2)}.    
     
Since $L$ acts on $Y'$, $Q$ acts on $\mathfrak{n} + Y'$.
We have a commutative diagram 
\begin{equation} 
\begin{CD} 
Spin(m) \times^{Q}(\mathfrak{n} + Y') @>>> Z \\ 
@V{\hat{\tau}}VV @VVV \\ 
Spin(m) \times^{Q}(\mathfrak{n} + \bar{O}_{{\mathbf p}'}) @>{s}>> \bar{O}_{\mathbf p}      
\end{CD} 
\end{equation}

Here $Z$ is the Stein factorization of $s \circ \hat{\tau}$. 
The horizontal maps are both crepant partial resolutions. 
Let $\pi: X \to \bar{O}_{\mathbf p}$ be the finite covering of 
$\bar{O}_{\mathbf p}$ associated with the universal covering 
of $O_{\mathbf p}$. Then $\pi$ factorizes as $X \stackrel{\rho}\to Z \to \bar{O}_{\mathbf p}$. By the construction $\mathrm{deg}(\rho) = 2$. Since $SO(m)$ acts on $Spin(m) \times^{Q}(\mathfrak{n} + Y')$, the finite map $Z \to \bar{O}_{\mathbf p}$ is 
an $SO(m)$-cover. Therefore this cover is nothing but the cover $Y \to \bar{O}_{\mathbf p}$ determined by 
the surjection $\pi_1(O_{\mathbf p}) \to (\mathbf{Z}/2\mathbf{Z})^{\oplus \mathrm{max}(a-1,0)}$.  
Then we have the following claim. 

\begin{Claim}\label{Claim (3.5.2)} The $Q$-action on $\mathfrak{n} + (Y')^0$ does not 
lift to a $Q$-action on $\mathfrak{n} + (X')^0$. Moreover, we have  
$$\pi_1(Spin(m) \times^Q (\mathfrak{n} + (Y')^0)) = \{1\}. $$ \end{Claim} 
 
{\em Proof}. Applying the homotopy exact sequence to the $Q$-bundle: $Spin(m) \times (\mathfrak{n} + (Y')^0) \to Spin(m) \times^Q (\mathfrak{n} + (Y')^0)$, we get an exact sequence 
$$\pi_1(Q) \stackrel{ob}\to \pi_1(Spin(m) \times (\mathfrak{n} + (Y')^0)) \to  
\pi_1(Spin(m) \times^Q (\mathfrak{n} + (Y')^0)) \to 1.$$
The map $ob$ is induced from a map $$\beta_{y_0}: Q \to Spin(m) \times (\mathfrak{n} + (Y')^0) \:\: 
q \to (q^{-1}, q(0 + y_0))$$ with some $y_0 \in (Y')^0$. We shall prove that $ob$ is the obstruction map 
to lifting the $Q$ action on $\mathfrak{n} + (Y')^0$ to a $Q$-action on $\mathfrak{n} 
+ (X')^0$. Let us consider the diagram  

\begin{equation} 
\begin{CD} 
Q \times (\mathrm{n} + (X')^0) @.  \mathrm{n} + (X')^0\\ 
@VVV @VVV \\ 
Q \times (\mathfrak{n} + (Y')^0) @>>> \mathfrak{n} + (Y')^0      
\end{CD} 
\end{equation}

Here the horizontal map on the bottom row is given by 
$$Q \times (\mathfrak{n} + (Y')^0) \to \mathfrak{n} + (Y')^0 \:\: \:\: (q, n + y) 
\to q(n + y).$$  The map 
$$\alpha: Q \times (\mathrm{n} + (X')^0) \to \mathfrak{n} + (Y')^0$$ lifts to 
a map $$\tilde{\alpha}:  Q \times (\mathrm{n} + (X')^0) \to \mathfrak{n} + (X')^0$$ 
if and only if 
$$\alpha_*: \pi_1(Q \times (\mathrm{n} + (X')^0) \to \pi_1(\mathfrak{n} + (Y')^0)$$ is 
the zero map. If such a lift $\tilde{\alpha}$ exists, then we can find an element 
$\tau \in \mathrm{Aut}_{(Y')^0} (X')^0$ such that  
$$Q \times (\mathrm{n} + (X')^0) \stackrel{\tilde{\alpha}}\to \mathfrak{n} + (X')^0 
\stackrel{id + \tau}\to \mathfrak{n} + (X')^0$$ is a group 
action. We take a lift $\tilde{y}_0 \in (X')^0$ of $y_0 \in (Y')^0)$. Define a map 
$i_{\tilde{y}_0}: Q \to Q \times (\mathfrak{n} + (X')^0)$ by $i_{\tilde{y}_0}(q) = 
(q, 0 + \tilde{y}_0)$.  
Then we have a commutative diagram 

\begin{equation} 
\begin{CD} 
Q @>{\beta_{y_0}}>>  Spin(m) \times (\mathfrak{n} + (Y')^0)\\ 
@V{i_{{\tilde y}_0}}VV  @V{pr_2}VV \\ 
Q \times (\mathfrak{n} + (X')^0) @>{\alpha}>> \mathfrak{n} + (Y')^0     
\end{CD} 
\end{equation} 

Correspondingly we have a commutative diagram of fundamental groups 

\begin{equation} 
\begin{CD} 
\pi_1(Q) @>{ob \: = \: (\beta_{y_0})_*}>>  \pi_1(Spin(m) \times (\mathfrak{n} + (X')^0)\\ 
@V{\cong}VV  @V{\cong}VV \\ 
\pi_1(Q \times (\mathfrak{n} + (X')^0) @>{\alpha_*}>> \pi_1(\mathfrak{n} + (Y')^0)       
\end{CD} 
\end{equation} 

The vertical maps are isomorphisms because $\pi_1((X')^0) = \{1\}$ 
and $\pi_1(Spin(m)) = \{1\}$. Therefore  $ob$ is the obstruction map to a lift.

We shall prove that $ob$ is not the zero map. Let us consider the composite  
$$\bar{\beta}_{y_0}: Q \stackrel{\beta_{y_0}}\to Spin(m) \times (\mathfrak{n} + (Y')^0) \stackrel{pr_2}\to 
\mathfrak{n} + (Y')^0 \stackrel{pr_2}\to (Y')^0.$$ 
If we restrict $\bar{\beta}_{y_0}$ to $U \subset Q$, then it is a constant map 
sending all elements $g \in U$ to $y_0 \in (Y')^0$. 
By the homotopy exact sequence $$\pi_1(U) \to \pi_1(Q) \to \pi_1(L) \to 1,$$
we see that $\bar{\beta}_{y_0}$ induces homomorphisms 
$$\pi_1(L) \to \pi_1(Spin(m) \times (\mathfrak{n} + (Y')^0)) \to  
\pi_1(\mathfrak{n} + (Y')^0) \to \pi_1((Y')^0).$$
Here the last two maps are both isomorphisms because $\pi_1(Spin(m)) = 1$ 
and $\pi_1(\mathfrak{n}) = 1$. Since the $L$-action on $(Y')^0$ does not 
lift to $(X')^0$, the map $\pi_1(L) \to \pi_1((Y')^0)$ is not the zero map. 
This means that $ob = (\beta_{y_0})_*$ is not the zero map. 

Return to the original homotopy exact sequence. Then we have 
$\pi_1(Spin(m) \times (\mathfrak{n} + (Y')^0) = \mathbf{Z}/2\mathbf{Z}$. 
In our case the $Q$-action on $\mathfrak{n} + (Y')^0$ does not lift to a $Q$-action 
on $\mathrm{n} + (X')^0$. This means that $\pi_1(Spin(m) \times^Q (\mathfrak{n} + (Y')^0) = \{1\}$.   $\square$ \vspace{0.2cm}

We are now in a position to consider a rather odd partition $\mathbf{p}$. 
Let $\mathbf{p} = [d^{r_d}, (d-1)^{d_{r-1}}, ..., 2^{r_2}, 1^{r_1}]$ be a partition 
of $m$ such that 

(i) $r_i$ is even for each even $i$.     

Let $i$ and $j$ be different members of $\mathbf{p}$. 
Then $i$ and $j$ are called adjacent if there are no members of $\mathbf{p}$ 
between $i$ and $j$ except themselves.   
We consider the following conditions for $\mathbf{p}$:  

(iii) For any couple $(i,j)$ of adjacent members, $\vert i - j \vert \leq 4$. 
Moreover, $\vert i -j \vert = 4$ occurs only when $i$ and $j$ are both odd. 
The smallest member of $\mathbf{p}$ is smaller than $4$.    

For example, both $[11,8^2,7,3]$ and $[7,6^2,3^2]$ satisfy the conditions (i) and 
(iii). Then we have

\begin{Prop}\label{Proposition (3.6)} Assume that a partition $\mathbf{p}$ satisfies the conditions (i) and (iii). Moreover, assume that $\mathbf{p}$ is rather odd.  Let $\pi : X \to \bar{O}_{\mathbf{p}}$ be the finite covering associated with the 
universal covering of $O_{\mathbf p} \subset so(m)$. Then $\mathrm{Codim}_X\mathrm{Sing}(X) \geq 4$. In particular, $X$ has {\bf Q}-factorial terminal singularities.  \end{Prop}

$X$ is {\bf Q}-factorial by Proposition \ref{Proposition (3.1)}, (2).  If $\mathrm{Codim}_X\mathrm{Sing}(X) \geq 4$, then $X$ has terminal singularities by Corollary \ref{Corollary (0.3)}. 
We will prove Proposition \ref{Proposition (3.6)} after some preliminaries. 
Assume that $\mathbf{p}$ is a rather odd partition with (i) and (iii).  Recall that a gap member $i$ of $\mathbf{p}$ is a member of $\mathbf{p}$ 
such that $i > 1$ and $r_{i-1} = 0$. For each gap member $i$, we will take an orbit $O_{\bar{\mathbf p}} \subset \bar{O}_{\mathbf p}$ and look at the singularity of 
$\bar{O}_{\mathbf p}$ along $O_{\bar{\mathbf p}}$ by using \cite{K-P}, 3.4. 
Let us say that a gap member $i$ is {\em exceptional} 
if $i$ is an odd number with $i \geq 5$ and $r_{i-1} = r_{j-2} = r_{j-3} = 0$. In this case $r_i = 1$ and $r_{i-4} = 1$ by the rather oddness and the condition (iii). Otherwise we say that $i$ is {\em ordinary}. 
For an exceptional gap member $i$, one can find a nilpotent orbit $O_{\bar{\mathbf p}} 
\subset \bar{O}_{\mathbf p}$ with 
$$\bar{\mathbf p} = [d^{r_d}, ..., (i+1)^{r_{i+1}}, (i-2)^2, (i-5)^{r_{i-5}}, ..., 1^{r_1}].$$ 
Then $\mathrm{Codim}_{O_{\mathbf p}}O_{\bar{\mathbf p}} = 2$ and 
the transverse slice for $O_{\bar{\mathbf p}} 
\subset \bar{O}_{\mathbf p}$ is an $A_3$-surface singularity. 
 
If $i$ is an even gap member, then there are two cases. The first case is 
when $r_{i-1} = r_{i-2} = 0$ and $r_{i-3} \ne 0$. The second case is when 
$r_{i-1} = 0$, but $r_{i-2} \ne 0$ (the case $i = 2$ being in this case).  
In the first case, $r_i$ is a non-zero 
even number by the condition (i); hence, $r_i \geq 2$. Then one can find 
a nilpotent orbit $O_{\bar{\mathbf p}} 
\subset \bar{O}_{\mathbf p}$ with
$$\bar{\mathbf p} = [d^{r_d}, ..., i^{r_i -2}, (i-1)^3, (i-3)^{r_{i-3}-1}, ..., 1^{r_1}].$$
In the second case, $r_i$ and $r_{i-2}$ are both nonzero even numbers; hence 
$r_i$, $r_{i-2} \geq 2$.       
Then one can find 
a nilpotent orbit $O_{\bar{\mathbf p}} 
\subset \bar{O}_{\mathbf p}$ with
$$\bar{\mathbf p} = [d^{d_r}, ..., i^{r_i -2}, (i-1)^4, (i-2)^{r_{i-2} - 2}, ..., 1^{r_1}].$$ 
In both cases, $\mathrm{Codim}_{O_{\mathbf p}}O_{\bar{\mathbf p}} = 2$. 
In the first case, the transverse slice $S$ for $O_{\bar{\mathbf p}} 
\subset \bar{O}_{\mathbf p}$ is an $A_1$-surface singularity. 
Next let us consider the second case. If $\mathbf{p}$ is not very even, then  $S$ is isomorphic to a union of two $A_1$-surface singularities intersecting each other in the singular point. In particular, $\bar{O}_{\mathbf p}$ is not normal along  $O_{\bar{\mathbf p}}$. 
If $\mathbf{p}$ is very 
even, $S$ is an $A_1$-surface singularity. This can be understood more conceptually. 
In fact, there are two orbits $O_{\mathbf p}^{\pm}$ with Jordan type $\mathbf{p}$.  
The transverse slice for $${O}_{\bar{\mathbf{p}}} \subset \bar{O}_{\mathbf p}^{+} \cup \bar{O}_{\mathbf p}^{-}$$ is then isomorphic to the union of two $A_1$-surface singularities intersecting each other in the singular point. Since we denote by 
$O_{\mathbf p}$ one of the two orbits, the transverse slice for $O_{\bar{\mathbf p}} 
\subset \bar{O}_{\mathbf p}$ is an $A_1$-surface singularity.  

If $i$ is an ordinary,  odd gap member, then there are two cases. 
The first case is when $r_{i-1} = r_{i-2} = 0$ and $r_{i-3} \ne 0$. 
The second case is when $r_{i-1} = 0$, but $r_{j-2} \ne 0$.
In the first case, $r_{i-3}$ is a nonzero even number by the condition (i); hence $r_{i-3} \geq 2$. Then one can find a nilpotent orbit $O_{\bar{\mathbf p}} 
\subset \bar{O}_{\mathbf p}$ with 
$$\bar{\mathbf p} = [d^{r_d}, ..., i^{r_i -1}, (i-2)^3, (i-3)^{r_{i-3}-2}, ..., 1^{r_1}].$$ 
In the second case, one can find a nilpotent orbit $O_{\bar{\mathbf p}} 
\subset \bar{O}_{\mathbf p}$ with
$$\bar{\mathbf p} = [d^{r_d}, ..., i^{r_i-1}, (i-1)^2, (i-3)^{r_{i-3}}, ..., 1^{r_1}].$$
In this last case, $\bar{\mathbf p}$ may possibly be very even. In such a case, 
we will have two orbits $O_{\bar p}^{\pm}$.  
In both cases, $\mathrm{Codim}_{O_{\mathbf p}}O_{\bar{\mathbf p}} = 2$ and 
the transverse slices for $O_{\bar{\mathbf p}} 
\subset \bar{O}_{\mathbf p}$ are $A_1$-surface singularities. 

As a consequence, for an ordinary gap member $i$ of $\mathbf{p}$, 
the transverse slice for $O_{\bar{\mathbf p}} \subset \bar{O}_{\mathbf p}$ 
is an $A_1$-surface singularity or the union of two $A_1$-surface 
singularities. In the first case we call $i$ an {\em ordinary gap member of type $A_1$}. 
In the latter case we call $i$ an {\em ordinary gap member of type $A_1 \cup A_1$}.

Let $\{i_1, i_2, ..., i_k\}$ be the set of all gap members of $\mathbf{p}$. As defined above, for these gap members, 
we have nilpotent orbits 
$$O_{\bar{\mathbf p}_1}, ..., O_{\bar{\mathbf p}_k}.$$
More exactly, when $\bar{\mathbf{p}}_j$ is very even, we have two different orbits 
with Jordan type $\bar{\mathbf{p}}_j$. In such a case we understand that  
$O_{\bar{\mathbf{p}}_j}$ is both of them. 
  
To prove that $\mathrm{Codim}_X\mathrm{Sing}(X) 
\geq 4$, we only have to show that $X$ is smooth along $\pi^{-1}(O_{\bar{\mathbf p}_{j}})$ for 
each $1 \le j \le k$. 
Let $S_j$ be the transversal slices for $O_{\bar{\mathbf p}_j} \subset 
\bar{O}_{\mathbf p}$. More exactly, when $\bar{\mathbf p}_j$ is very even, 
we must consider two slices $S_j^{\pm}$ for $O_{\bar{\mathbf p}_j}^{\pm}$.  
The claim is then equivalent to showing that $\pi^{-1}(S_j)$ are disjoint union of finite copies of  $(\mathbf{C}^2, 0)$.  Note that $S_j$ is an $A_3$-surface 
singularity if $i_j$ is exceptional, and $S_j$ is an $A_1$-surface singularity or of type $A_1 \cup A_1$ if $i_j$ is ordinary. 

The gap member is closely related to the notion of {\em induced orbits}. 

(Type I)
Let $i$ be a gap member of $\mathbf{p}$. Let $O_{\mathbf p} \subset so(m)$ be a nilpotent 
orbit with Jordan type.  
Put $r := r_d + ... + r_{i}$ and 
let $\bar{Q} \subset SO(m)$ be a parabolic subgroup of flag type $(r, m-2r, r)$, and put $Q := \rho_m^{-1}(\bar{Q}) \subset Spin(m)$. Take a Levi 
decomposition $\mathfrak{q} = \mathfrak{n} \oplus \mathfrak{l}$.  Notice that 
$\mathfrak{l} = \mathfrak{gl}(r) \oplus so(m-2r)$. 
There is a nilpotent 
orbit $O_{\mathbf{p}'}$ of $so(m-2r)$ with Jordan type 
$$\mathbf{p}' = [(d-2)^{r_d}, ..., (i - 1)^{r_{i + 1}}, (i-2)^{r_{i} + r_{i-2}}, (i-3)^{r_{i-3}},..., 1^{r_1}]$$ such that 
$O_{\mathbf p} = \mathrm{Ind}^{\mathfrak g}_{\mathfrak l}(O_{\mathbf{p}'})$.  

\begin{Claim}\label{Claim (3.6.1)}(cf. \cite{He}, Theorem 7.1, (d)). The generalized Springer map 
$$\mu: Spin(m) \times^{Q}(\mathfrak{n} + \bar{O}_{\mathbf{p}'}) \to \bar{O}_{\mathbf p}$$ is a birational map. 
\end{Claim}

{\em Proof}. First notice that $Spin(m) \times^{Q}(\mathfrak{n} + \bar{O}_{\mathbf{p}'})
= SO(m) \times^{\bar{Q}}(\mathfrak{n} + \bar{O}_{\mathbf{p}'})$.  
It is enough to prove the claim by replacing $Spin(m)$ by $SO(m)$.
Then it is completely analogous to the corresponding statement for $Sp(2n)$ in 
the previous section.
 
For $x \in O_{\mathbf p}$, we take a basis $\{e(l,j)\}_{(l,j) \in Y(\mathbf{p})}$ of $\mathbf{C}^{m}$ so that 

(a) $\{e(l,j)\}$ is a Jordan basis for $x$, i.e. 
$x\cdot e(l,j) = e(l-1,j)$ for $l > 1$ and $x\cdot e(1,j) = 0$. 

(b) $\langle e(l,j), e(p,q) \rangle \ne 0$ if and only if $p = d_j -l + 1$ and 
$q = \beta(j)$. Here $\beta$ is a permutation of $\{1,2, ..., s\}$ such that 
$\beta^2 = id$, $d_{\beta(j)} = d_j$, and $\beta (j) \ne j$ if $d_j$ is even 
(cf. \cite{S-S}, p.259, see also \cite{C-M}, 5.1). 

Put $F:= \sum_{1 \le j \le r} \mathbf{C} e(1,j)$. Then 
$F \subset F^{\perp}$ is an isotropic flag such that $x\cdot F = 0$ and $x \cdot \mathbf{C}^{m} \subset F^{\perp}$ and $x$ is an endomorphism of 
$F^{\perp}/F$ with Jordan type $\mathbf{p}'$.  This is actually a unique isotropic 
flag of type $(r, m-2r, r)$ satisfying these properties. Hence $\mu^{-1}(x)$ consists 
of one element.  $\square$

For later convenience, we also introduce an induction of type II.  

(Type II)
Let $i$ be an odd member of $\mathbf{p}$ with $r_i = 2$. Put 
$r := r_d + ... + r_{i-1} + 1$ and let $\bar{Q} \subset SO(m)$ be a parabolic subgroup of flag type $(r, m-2r, r)$ and put $Q := \rho_m^{-1}(\bar{Q}) \subset Spin(m)$. Take a Levi 
decomposition $\mathfrak{q} = \mathfrak{n} \oplus \mathfrak{l}$. Notice that 
$\mathfrak{l} = \mathfrak{gl}(r) \oplus so(m-2r)$. 
There is a nilpotent 
orbit $O_{\mathbf{p}'}$ of $so(m-2r)$ with Jordan type 
$$\mathbf{p}' = [(d-2)^{r_d}, ..., i^{r_{i+2}}, (i-1)^{r_{i+1} + 2 + r_{i-1}}, (i-2)^{r_{i-2}}, ..., 1^{r_1}]$$ 
such that $O_{\mathbf p} = \mathrm{Ind}^{\mathfrak g}_{\mathfrak l}(O_{\mathbf{p}'}).$

\begin{Claim}\label{Claim (3.6.2)}(cf. \cite{He}, Theorem 7.1, (d)) The generalized Springer map 
$$\mu: Spin(m) \times^{Q}(\mathfrak{n} + \bar{O}_{\mathbf{p}'}) \to \bar{O}_{\mathbf p}$$ is a birational map if one of the following holds: 

(a) $\mathbf{p}'$ is very even, or   

(b) $m = 2r$. 

Otherwise $\mu$ is a generically finite map of degree $2$. \end{Claim} 

{\em Proof}. Since $Spin(m) \times^{Q}(\mathfrak{n} + \bar{O}_{\mathbf{p}'})
= SO(m) \times^{\bar{Q}}(\mathfrak{n} + \bar{O}_{\mathbf{p}'})$,   
it is enough to prove the claim by replacing $Spin(m)$ by $SO(m)$.
For $x \in O_{\mathbf p}$, we take the same basis 
$\{e(l,j)\}$ of $\mathbf{C}^m$ as the previous claim. After a suitable change of the basis, 
we may assume that $\beta (r) = r +1$ and $\beta (r+1) = r$. 
We put $F := \sum_{1 \le j \le r} \mathbf{C} e(1,j)$. 
Then $F \subset F^{\perp}$ is an isotropic flag such that $x\cdot F = 0$ and $x \cdot \mathbf{C}^{m} \subset F^{\perp}$ and $x$ is an endomorphism of 
$F^{\perp}/F$ with Jordan type $\mathbf{p}'$.
On the other hand, put $F' := \sum_{1 \le j \le r+1, \: j \ne r}\mathbf{C} e(1, j)$. 
Then $F' \subset (F')^{\perp}$ is an isotropic flag such that $x\cdot F' = 0$ and $x \cdot \mathbf{C}^{m} \subset (F')^{\perp}$ and $x$ is an endomorphism of 
$(F')^{\perp}/F'$ with Jordan type $\mathbf{p}'$. By a similar argument as in the proof of Claim \ref{Claim (2.3.2)}, the isotropic flags of type $(r, m-2r, r)$ with these properties are exactly two flags above. 

Assume that $\bar{Q}$ is the parabolic subgroup of $SO(m)$ stabilizing the flag 
$F \subset F^{\perp}$.
We have an $SO(m)$-equivariant (locally closed) immersion  $$\iota: SO(m) \times^{\bar{Q}} (\mathfrak{n} + O_{\mathfrak{p}'}) 
\subset SO(m)/\bar{Q} \times so(m), \:\:\: [g, y] \to (g\bar{Q}, Ad_g(y)).$$ 
Consider $(F \subset F^{\perp}, x)$ and $(F' \subset (F')^{\perp}, x)$ as elements 
of $SO(m)/\bar{Q} \times so(2n)$. Note that $\iota ([1, x]) = (\bar{Q}, x) = (F \subset F^{\perp}, x)$. 
Let $\bar{Q}'$ be the parabolic subgroup of $SO(m)$ stabilizing the flag 
$F' \subset (F')^{\perp}$. If $m \ne 2r$, then $\bar{Q}$ and $\bar{Q}'$ are conjugate 
parabolic subgroups. Moreover, if $\mathbf{p}'$ is not very even, a nilpotent orbit 
of $so(m-2r)$ with Jordan type $\mathbf{p}'$ is unique. 
Hence, if neither (a) nor (b) holds, we have $(\bar{Q}', x) \in \mathrm{Im}(\iota)$ 
by the same argument as in $Sp(2n)$.  
This means that $\mu^{-1}(x)$ consists of two elements. 

If (b) holds and $r \ne 1$, then $Q'$ is not conjugate to $Q$. In particular, 
$(\bar{Q}', x) \notin \mathrm{Im}(\iota)$. This means that 
$\mu^{-1}(x)$ consists of one element. When $m = 2$ and $r = 1$, 
two flags $F \subset F^{\perp}$ 
and $F' \subset (F')^{\perp}$ are different, but $\bar{Q} = \bar{Q}'$. 
Therefore $\mu^{-1}(x)$ consists of one element.

Finally, if (a) holds, we have one more nilpotent orbit $O_{\mathbf{p}'}^{-}$ of 
$so(m-2r)$ different from $O_{\mathbf{p}'}$ with Jordan type $\mathbf{p}'$. 
We have one more $SO(m)$-equivariant immersion 
$$\iota^{-}: SO(m) \times^{\bar{Q}} (\mathfrak{n} + O^{-}_{\mathfrak{p}'}) 
\subset SO(m)/\bar{Q} \times so(m).$$ 
We have $\mathrm{Im}(\iota) \cap \mathrm{Im}(\iota^-) = \emptyset$. 
Then $(F' \subset (F')^{\perp}, x) \in \mathrm{Im}(\iota^-)$. This implies that 
$\mu^{-1}(x)$ consists of one element. 
$\square$ \vspace{0.2cm}

For each gap member $i_j$, consider an induction of type I:  $$\mathbf{p}'_j = [(d-2)^{r_d}, ..., (i_j - 1)^{r_{i_j + 1}}, (i_j-2)^{r_{i_j} + r_{i_j-2}}, (i_j-3)^{r_{i_j-3}},..., 1^{r_1}].$$ We then have a generalized Springer map 
$$\mu_j: Spin(m) \times^{Q_j}(\mathfrak{n}_j + \bar{O}_{\mathbf{p}'_j}) \to \bar{O}_{\mathbf p},$$ which is a birational map.

\begin{Lem}\label{Lemma (3.7)}  
(1) Assume that $i_j$ is an ordinary gap member of $\mathbf{p}$. 
If $S_j$ is an $A_1$-surface singularity, then $\mu_j^{-1}(S_j)$ is the minimal resolution of $S_j$. If $S_j$ is of type $A_1 \cup A_1$, then  $\mu_j^{-1}(S_j)$ is the minimal resolution 
of the normalization of $S_j$. 
On the other hand, $\mu_j^{-1}(S_l) \to S_l$ is an isomorphism for each 
$l \ne j$. 

(2) Assume that $i_j$ is an exceptional gap member of $\mathbf{p}$. 
Then $\mu_j^{-1}(S_j)$ is a crepant partial  resolution of $S_j$  
with one exceptional curve $C_j \cong \mathbf{P}^1$. $\mu_j^{-1}(S_j)$ 
has $A_1$-surface singularities at two points $p_j^{\pm} \in C_j$. 
On the other hand, $\mu_j^{-1}(S_l) \to S_l$ is an isomorphism for each 
$l \ne j$. \end{Lem} 

{\em Proof}. (1) Assume that $i_j$ is ordinary. Then the gap members of $\mathbf{p}'_j$ are $$i_1 -2, \: ..., \: i_{j-1}-2,\:  i_{j+1}, ...,  
i_k.$$ For the later convenience we put $$i'_1 := i_1 -2,\: ..., i'_{j-1} := i_{j-1} -2, \: 
i'_{j+1} := i_{j+1}, ..., i'_k := i_k.$$ Corresponding to these gap members, we get 
nilpotent orbits $O_{(\overline{\mathbf{p}'_j})_1}$, ..., $O_{(\overline{\mathbf{p}'_j})_{j-1}}$, 
$O_{(\overline{\mathbf{p}'_j})_{j+1}}$, ..., $O_{(\overline{\mathbf{p}'_j})_k}$ in $\bar{O}_{\mathbf{p}'_j}$. 
Here we follow the notation explained below Proposition \ref{Proposition (3.6)}. 
These are irreducible components of $\mathrm{Sing}(\bar{O}_{\mathbf{p}'_j})$ which have codimension 2 in $\bar{O}_{\mathbf{p}'_j}$. Take  
$l$ so that $1 \le l \le k$ and $l \ne j$. $\overline{(\mathbf{p}'_j)}_l$ is very even if and 
only if $\bar{\mathbf p}_l$ is very even.  
If $i_l$ is ordinary gap member of $\mathbf{p}$ of type $A_1$ (resp. of type $A_1 \cup A_1$), then $i'_l$ is also an ordinary gap member of $\mathbf{p}'_j$ of type $A_1$ (resp. of type $A_1 \cup A_1$). If $i_l$ is an 
exceptional gap member of $\mathbf{p}$, then $i'_l$ is an exceptional gap member of $\mathbf{p}'_j$.  
For each $l$, we have a natural embedding 
$$Spin(m) \times^{Q_j}(\mathfrak{n}_j + \bar{O}_{(\overline{\mathbf{p}'_j})_l}) \subset 
Spin(m) \times^{Q_j}(\mathfrak{n}_j + \bar{O}_{\mathbf{p}'_j}).$$
One can check that $$\mu_j(\:\: Spin(m) \times^{Q_j}(\mathfrak{n}_j + \bar{O}_{(\overline{\mathbf{p}'_j})_l})\:\:) = \bar{O}_{\bar{\mathbf p}_l}.$$
This implies the last statement of (1). There are no singularities of 
$Spin(m) \times^{Q_j}(\mathfrak{n}_j + \bar{O}_{\mathbf{p}'_j})$ lying over $O_{\bar{\mathbf p}_j}$. This implies the first and the second statements of (1).  

(2) Assume that $i_j$ is exceptional. Then $i_j - 2$ is still a gap member of $\mathbf{p}'_j$. Hence the gap members of $\mathbf{p}'_j$ are 
$$i_1 -2, \: ..., \: i_{j-1}-2, i_j -2, i_{j+1}, ...,  i_k.$$ 
We put $$i'_1 := i_1 -2,\: ..., i'_{j-1} := i_{j-1} -2, i'_j := i_j -2, 
i'_{j+1} := i_{j+1}, ..., i'_k := i_k.$$ Corresponding to these gap members, 
we get nilpotent orbits $O_{(\overline{\mathbf{p}'_j})_1}$, ..., $O_{(\overline{\mathbf{p}'_j})_{j}}$, ..., $O_{(\overline{\mathbf{p}'_j})_k}$. Note that we have an additional orbit 
$O_{(\overline{\mathbf{p}'_j})_{j}}$ in the exceptional case.  Take $l$ so that 
$1 \le l \le k$. For each $l$, we have a natural embedding 
$$Spin(m) \times^{Q_j}(\mathfrak{n}_j + \bar{O}_{(\overline{\mathbf{p}'_j})_l}) \subset 
Spin(m) \times^{Q_j}(\mathfrak{n}_j + \bar{O}_{\mathbf{p}'_j}).$$
One can check that $$\mu_j(\:\: Spin(m) \times^{Q_j}(\mathfrak{n}_j + \bar{O}_{(\overline{\mathbf{p}'_j})_l})\:\:) = \bar{O}_{\bar{\mathbf p}_l}.$$
For $l \ne j$, we apply the same argument as in (1), 
and we see that the last statement of (2) holds true. 
Let us consider the case $l = j$. 
We have $$\bar{\mathbf p}_j = [d^{r_d}, ..., (i_j + 1)^{r_{i_j + 1}}, (i_j -2)^2, (i_j - 5)^{r_{i_j-5}}, ..., 1^{r_1}],$$
$$\mathbf{p}'_j = [(d-2)^{r_d}, ..., (i_j-1)^{r_{i_j +1}}, i_j -2, i_j-4, (i_j-5)^{r_{i_j-5}}, ..., 1^{r_1}],$$ 
$$(\overline{{\mathbf p}'_j})_j = [(d-2)^{r_d}, ..., (i_j-1)^{r_{i_j +1}}, (i_j -3)^2, (i_j-5)^{r_{i_j-5}}, ..., 1^{r_1}].$$ We then see that $(\overline{{\mathbf p}'_j})_j \to \bar{\mathbf p}_j$ is an induction of type II. 
When $(\overline{{\mathbf p}'_j})_j $ is not very even, the generalized Springer map 
$$Spin(m) \times^{Q_j}(\mathfrak{n}_j + \bar{O}_{(\overline{\mathbf{p}'_j})_j}) \to 
\bar{O}_{\bar{\mathbf p}_j}$$ is a generically finite map of degree 2. Moreover, it is an 
etale double cover over $O_{\bar{\mathbf p}_j}$. Since $\bar{O}_{\mathbf{p}'_j}$ has 
an $A_1$-surface singularity along $O_{(\overline{\mathbf{p}'_j})_j}$, $Spin(m) \times^{Q_j}(\mathfrak{n}_j + \bar{O}_{\mathbf{p}'_j})$ has an 
$A_1$-singularity along 
$Spin(m) \times^{Q_j}(\mathfrak{n}_j + O_{(\overline{\mathbf{p}'_j})_j})$. This implies 
the first and second statements of (2).  Assume that $(\overline{{\mathbf p}'_j})_j $ is very even. 
Then there are two orbits $O^{\pm}_{(\overline{{\mathbf p}'_j})_j }$. Then each generalized Springer 
map $$Spin(m) \times^{Q_j}(\mathfrak{n}_j + \bar{O}^{\pm}_{(\overline{\mathbf{p}'_j})_j}) \to \bar{O}_{\bar{\mathbf p}_j}$$ is a birational map. Moreover, it is an isomorphism over 
$O_{\bar{\mathbf p}_j}$. $Spin(m) \times^{Q_j}(\mathfrak{n}_j + \bar{O}_{\mathbf{p}'_j})$ has an $A_1$-surface singularity along 
two disjoint subvarieties $Spin(m) \times^{Q_j}(\mathfrak{n}_j + O^{+}_{(\overline{\mathbf{p}'_j})_j})$ and $Spin(m) \times^{Q_j}(\mathfrak{n}_j + O^{-}_{(\overline{\mathbf{p}'_j})_j})$. 
Therefore the first and the second statements of (2) still holds in this case.  $\square$
\vspace{0.2cm}

In the above, we associate a nilpotent orbit $O_{\bar{\mathbf p}} \subset 
\bar{O}_{\mathbf p}$ with a gap member $i$ of $\mathbf{p}$. Assume that 
$\mathbf{p}$ is not very even and $i$ is an even gap member with $r_{i-2} \ne 0$. We already remarked that the transverse slice 
$S$ for ${O}_{\bar{p}} \subset \bar{O}_{\mathbf p}$ is of type $A_1 \cup A_1$.
In this case, $\bar{O}_{\mathbf p}$ is not normal. So let us consider the normalization map  $\nu: \tilde{O}_{\mathbf p} \to \bar{O}_{\mathbf p}$.    
Before going to the proof of Proposition \ref{Proposition (3.6)}, we will give a description of 
$\nu$. 

\begin{Lem}\label{Lemma (3.8)} $\nu^{-1}(O_{\bar{\mathbf p}})$ is a connected, etale double cover 
of $O_{\bar{\mathbf p}}$. \end{Lem}

{\em Proof}. By the description of $S$, the statement is clear except that  
$\nu^{-1}(O_{\bar{\mathbf p}})$ is connected.
 
We apply inductions to $\mathbf{p}$ repeatedly by using the gap members different 
from $i$ and finally get a partition $\mathbf{p}'$ such that $\mathbf{p}'$ has a unique 
gap member $i'$, which comes from the originally fixed gap member $i$. For example, start with $\mathbf{p} = [10^2, 7, 5, 4^2, 2^2]$ and $i = 4$. 
Then $10$, $7$, $2$ are gap members different from $4$. Take an induction for $2$ 
to get new partition $\mathbf{p}_1 = [8^2, 5, 3, 2^2]$. The partition $\mathbf{p}_1$ has 
gap members $8$, $5$ except $2$, which comes from the originally fixed gap member 
$4$. We next take an induction for $5$ to get $\mathbf{p}_2 = [6^2, 3^2, 2^2]$. 
Finally we get $\mathbf{p}' = [4^2, 3^2, 2^2]$ by taking an induction for $6$. 

Write the partition $\mathbf{p}'$ as 
$$\mathbf{p}' = [{d'}^{r_{d'}}, ..., (i'+1)^{r'_{i' + 1}}, i'^{r'_{i'}}, (i'-2)^{r'_{i'-2}}, (i'-3)^{r'_{i'-3}}, ..., 1^{r'_1}].$$ 
By the construction, $r'_{d'}$, ..., $r'_{i' + 1}$, $r'_{i'}$, $r'_{i'-2}$, $r'_{i'-3}$, ... $r'_1$ are all 
nonzero.  There is a generalized Springer map 
$$\mu': Spin(m) \times^{Q'}(\mathfrak{n}' + \bar{O}_{{\mathbf p}'}) 
\to \bar{O}_{\mathbf p}.$$ 
Put 
$$\bar{{\mathbf p}'} = [{d'}^{r'_{d'}}, ..., {i'}^{r'_{i'} -2}, (i'-1)^4, (i'-2)^{r'_{i'-2} - 2}, ..., 1^{r'_1}].$$
Then $\bar{O}_{{\mathbf p}'}$ has $A_1 \cup A_1$-singularity along 
$O_{\bar{{\mathbf p}'}}$. Let $\nu': \tilde{O}_{{\mathbf p}'} \to \bar{O}_{{\mathbf p}'}$ 
be the normalization map.  
The subvariety $Spin(m) \times^{Q'}(\mathfrak{n}' + \bar{O}_{\bar{{\mathbf p}'}}) \subset  
Spin(m) \times^{Q'}(\mathfrak{n}' + \bar{O}_{{\mathbf p}'})$ is birationally  
mapped to $\bar{O}_{\bar{\mathbf p}}$ by $\mu'$. Put $\bar{\mu'} := \mu'\vert_{Spin(m) \times^{Q'}(\mathfrak{n}' + 
\bar{O}_{\bar{{\mathbf p}'}})}$. Then $(\bar{\mu}')^{-1}
(O_{\bar{\mathbf p}}) \subset Spin(m) \times^{Q'}(\mathfrak{n}' + O_{\bar{{\mathbf p}'}})$ 
and $(\bar{\mu}')^{-1}
(O_{\bar{\mathbf p}}) \to O_{\bar{\mathbf p}}$ is an isomorphism. 
Moreover, $\mu'$ induces a birational map 
of the normalizations of both sides: 
$${\mu'}^n: Spin(m) \times^{Q'}(\mathfrak{n}' + \tilde{O}_{{\mathbf p}'}) 
\to \tilde{O}_{\mathbf p},$$ which induces a dominating map 
$$Spin(m) \times^{Q'} (\mathfrak{n}' + (\nu')^{-1}(\bar{O}_{\bar{{\mathbf p}'}})) 
\to \nu^{-1}(\bar{O}_{\bar{\mathbf p}}).$$ 
Therefore it suffices to show that 
$(\nu')^{-1}(O_{\bar{{\mathbf p}'}})$ is connected in order to show that 
$\nu^{-1}(O_{\bar{\mathbf p}})$ is connected. In fact, if $(\nu')^{-1}(O_{\bar{{\mathbf p}'}})$ is connected, then $(\nu')^{-1}(\bar{O}_{\bar{{\mathbf p}'}})$ is irreducible. 
Then $\nu^{-1}(\bar{O}_{\bar{\mathbf p}})$ is also irreducible by the dominating 
map. This means that $\nu^{-1}(O_{\bar{\mathbf p}})$ is connected.  
By the argument above, we may assume that $\mathbf{p}$ contains all 
numbers from $1$ to $d$ except $i-1$. Apply the induction for the unique gap 
member $i$ of $\mathbf{p}$. We get a partition 
$$\mathbf{p}_{full} := [(d-2)^{r_d}, ..., (i-1)^{r_{i+1}}, (i-2)^{r_i + r_{i-2}}, (i-3)^{r_{i-3}}, ..., 
1^{r_1}].$$ Notice that $\mathbf{p}_{full}$ has full members, i.e. all numbers from 
$1$ to $d-2$. Put $r := r_d + ... + r_i$ and let $Q_{r, m-2r, r}$ be a parabolic subgroup 
of $SO(m)$ with flag type $(r, m-2r, r)$. Put $Q := \rho_m^{-1}(Q_{r,m-2r,r})$ for 
$\rho_m: Spin(m) \to SO(m)$.Then the Levi part $\mathfrak{l}$ of $\mathfrak{q}$ 
contains $so(m-2r)$ as a direct factor. In particular, 
the nilpotent orbit $O_{\mathbf{p}_{full}}$ is contained in $\mathfrak{l}$.
By \cite{Na 1}, Corollary (1.4.3), the normalization $\tilde{O}_{\mathbf{p}_{full}}$
of $\bar{O}_{\mathbf{p}_{full}}$ has {\bf Q}-factorial 
terminal singularities except when $\mathbf{p}_{full} = [2^{(m-2r-2)/2}, 1^2]$. 
 
Therefore the generalized Springer map 
$$\mu: Spin(m) \times^Q(\mathfrak{n} + \tilde{O}_{\mathbf{p}_{full}}) \to 
\bar{O}_{\mathbf p}$$ gives a {\bf Q}-factorial terminalization of 
$\bar{O}_{\mathbf p}$ except when 

(a) $\mathbf{p}_{full} = [2^{(m-2r-2)/2}, 1^2]$.

In case (a), $\tilde{O}_{\mathbf{p}_{full}}$ is not yet {\bf Q}-factorial. To construct a 
{\bf Q}-factorial terminalization of $\bar{O}_{\mathbf p}$, we take 
a parabolic subgroup $Q_{r, (m-2r)/2, (m-2r)/2, r}$ of $SO(m)$ and 
put $\tilde{Q} := \rho_m^{-1}(Q_{r, (m-2r)/2, (m-2r)/2, r})$. Then 
$$\tilde{\mu}: T^*(Spin(m)/{\tilde Q}) = Spin (m) \times^{\tilde Q} \tilde{\mathfrak{n}} 
\to \bar{O}_{\mathbf p}$$ gives a crepant resolution of $\bar{O}_{\mathbf p}$. 

Assume that (a) does not occur. Suppose that $\nu^{-1}(O_{\bar{\mathbf p}})$ is not connected,  Then 
we have two different $\mu$-exceptional divisors $E^{\pm}$ which respectively dominate the closures of two connected components of $\nu^{-1}(O_{\bar{\mathbf p}})$. 
We consider the following two cases separately 

(Case I) $\mathbf{p}_{full} = [1^2]$ 

(Case II) otherwise. 

When Case (II) occurs, $Q$ is a maximal parabolic subgroup.  Hence, the 
relative Picard number $\rho(\mu)$ of $\mu$ must be one. This contradicts that 
$\mathrm{Exc}(\mu)$ contains $2$ irreducible divisors. 
We next consider Case (I). In this case $Q$ is not maximal and $\rho(\mu) = 2$. 
It is easily checked that Case (I) happens only when 
$\mathbf{p} = [3^2,2^2]$. In this case $\pi_1(O_{\mathbf p}) = \{1\}$. Moreover, since 
the odd member $3$ has multiplicity $2$, $\tilde{O}_{\mathbf p}$ is 
not {\bf Q}-factorial. In fact, take $x \in O_{\mathbf p}$ and put $G := Spin (m)$. By the proof of Proposition \ref{Proposition (3.1)}, (1), we see that 
$\chi ((G^\phi)^0)$ is infinite. Now let us look at the proof of Proposition (1.1), (2).  
Then we see that $\chi ((G^x)^0)$ is infinite and $\mathrm{Pic}(G/(G^x)^0)$ is infinite. 
Since $\pi_1(O_{\mathbf p}) = \{1\}$, we have $G^x = (G^x)^0$. This means that $\mathrm{Pic}(O_{\mathbf p})$ is infinite; hence, 
$\tilde{O}_{\mathbf p}$ is not {\bf Q}-factorial.  
 
On the other hand, since $\rho(\mu) =2$, $[E^+]$ and $[E^{-}]$ 
span $\mathrm{NS}(\mu) \otimes \mathbf{Q}$, which means that $\tilde{O}_{\mathbf p}$ 
is {\bf Q}-factorial by \cite{Na 1}, Lemma (1.1.1). This is a contradiction. 

Finally assume that (a) occurs. In this case $\rho(\tilde{\mu}) = 2$. 
It is easily checked that (a) happens only when 
$\mathbf{p}  = [4^2,3^2,2^2]$. Then $\pi_1(O_{\mathbf p}) = \{1\}$. 
Moreover, since the odd member $3$ has multiplicity $2$, $\tilde{O}_{\mathbf p}$ is 
not {\bf Q}-factorial by the same argument as above. On the other hand, since $\rho(\tilde{\mu}) =2$, $[E^+]$ and $[E^{-}]$ 
span $\mathrm{NS}(\tilde{\mu}) \otimes \mathbf{Q}$, which means that $\tilde{O}_{\mathbf p}$  
is {\bf Q}-factorial by [ibid, Lemma (1.1.1)]. This is a contradiction. $\square$
\vspace{0.2cm} 
 
Let $\tau: Y \to \bar{O}_{\mathbf p}$ be the finite covering 
determined by the surjection $\pi_1(O_{\mathbf p}) \to (\mathbf{Z}/2\mathbf{Z})^{\oplus \mathrm{max}(a(\mathbf{p}) -1,0)}$ in 
Proposition \ref{Proposition (3.1)}, (1).  Then $\tau$ is the finite covering associated with the $SO(m)$-universal covering of $O$. The map 
$\pi$ factors through $Y$: 
$$X \stackrel{\rho}\to Y \stackrel{\tau}\to \bar{O}_{\mathbf p}.$$ 
Assume that $\mathbf{p}$ is a {\em rather odd} partition of $m$ 
Let $\{i_1, ..., i_k\}$ be the set of gap members of $\mathbf{p}$. 
For $1 \le j \le k$, we take a transverse slice $S_j$ for $O_{\bar{\mathbf p}_j} 
\subset \bar{O}_{\mathbf p}$. 

\begin{Prop}\label{Proposition (3.9)}
 
 (1) Assume that $i_j$ is an odd, ordinary gap member of $\mathbf{p}$ such that $r_{i_j - 2} \ne 0$. Then 
$a(\mathbf p) \geq 2$ and $S_j$ is an $A_1$-surface singularity.  
Moreover,  $\tau^{-1}(S_j)$ 
is a disjoint union of $2^{a(\mathbf p) -2}$ copies of $(\mathbf{C}^2,0)$.  
 
 (2-i) Assume that $i_j$ is an ordinary gap member of $\mathbf{p}$ such that 
$r_{i_j -2} = 0$. Then $S_j$ is an $A_1$-surface singularity. Moreover,  
$\tau^{-1}(S_j)$ is a disjoint union of $2^{\mathrm{max}(a(\mathbf p) -1,0)}$ copies of $A_1$-surface singularity. 

 (2-ii) Assume that $i_j$ is an even, ordinary gap member of $\mathbf{p}$ with $r_{i_j -2} \ne 0$. If $\mathbf{p}$ is very even, 
 then $S_j$ is an $A_1$-surface singularity. If $\mathbf{p}$ is not very even, then $S_j$ is a union of two $A_1$-singularities 
 $S_j^+$ and $S_j^-$ intersecting in one point. In the former case, $\tau^{-1}(S_j)$ is an $A_1$-surface 
 singularity. In the latter case, $\tau^{-1}(S_j^+)$ and $\tau^{-1}(S_j^-)$ are respectively disjoint unions of  
 $2^{\mathrm{max}(a(\mathbf p) -1,0)}$ copies of $A_1$-surface singularity. 
 
 (3) Assume that $i_j$ is an exceptional gap member of $\mathbf{p}$. Then $a(\mathbf p) \geq 2$ and $S_j$ is an $A_3$-surface singularity. Moreover, $\tau^{-1}(S_j)$ is a disjoint union of $2^{a(\mathbf p) -2}$ copies of an $A_1$-surface 
 singularity. 
 \end{Prop}
 
 {\em Proof}.
 For simplicity we respectively write $\bar{\mathbf p}$ for $\bar{\mathbf p}_j$, $\mathbf{p}'$ for $\mathbf{p}'_j$ and 
 $i$ for $i_j$. The generalized Springer map $\mu_j$ is simply denoted by $\mu$, and put $Q := Q_j$ ($\mu_j$ and 
 $Q_j$ are defined just before Lemma \ref{Lemma (3.7)}).  Moreover, we put $S := S_j$.

 (1)  It remains to prove the last statement.  We can write 
 $$\mathbf{p} = [d^{r_d}, ..., (i+1)^{r_{i+1}}, i, i-2, (i-3)^{r_{i-3}}, ..., 1^{r_1}],$$
 $$\bar{\mathbf p} := [d^{r_d}, ..., (i+1)^{r_{i+1}}, (i-1)^2, (i-3)^{r_{i-3}}, ..., 1^{r_1}]$$ and 
 $$\mathbf{p}' := [(d-2)^{r_d}, ..., (i-1)^{r_{i+1}}, (i-2)^2, (i-3)^{r_{i-3}}, ..., 1^{r_1}].$$
 Then $\mathbf{p}'$ is a partition of $m' := m - 2(r_d + ... + r_{i+1} + 1)$, which is not  rather odd because the odd member $i-2$ has multiplicity $2$. Let us consider the nilpotent orbit $O_{\mathbf{p}'} \subset so(m')$. Since $a(\mathbf{p}') = a(\mathbf{p}) -1$ and 
$\mathbf{p}'$ is not rather odd,  $\pi_1(O_{\mathbf{p}'})$ has order $2^{a(\mathbf{p}) - 2}$. Let $Q_{r_d + ... + r_{i+1} + 1, m', r_d + ... + r_{i+1} + 1} \subset SO(m)$ be a 
parabolic subgroup of flag type  $(r_d + ... + r_{i+1} + 1, m', r_d + ... + r_{i+1} + 1)$ and 
put $Q := \rho_m^{-1}(Q_{r_d + ... + r_{i+1} + 1, m', r_d + ... + r_{i+1} + 1})$ for 
the double covering $\rho_m: Spin(m) \to SO(m)$. We have a generalized 
Springer map 
$$\mu: Spin(m) \times^Q (\mathfrak{n} + \bar{O}_{\mathbf{p}'}) \to \bar{O}_{\mathbf p}.$$ 

Assume that $\bar{\mathbf p}$ is not very even. We calculate $\pi_1(O_{\mathbf p} 
\cup O_{\bar{\mathbf p}})$.
We can write $\mathrm{Sing}(\bar{O}_{\mathbf p})$ as a union of ${O}_{\bar{\mathbf p}}$ and finite number of other nilpotent orbits: 
$$\mathrm{Sing}(\bar{O}_{\mathbf p}) = {O}_{\bar{\mathbf p}} 
\sqcup {O}_1 \sqcup ... \sqcup {O}_k \sqcup {O}_{k+1} \sqcup ... \sqcup {O}_{k+l}$$ Here $\mathrm{Codim}_{{O}_{\mathbf p}}O_{\alpha}
= 2$ for $1 \le \alpha \le k$ and $\mathrm{Codim}_{{O}_{\mathbf p}}O_{\alpha} \geq 4$ for 
$k+1 \le \alpha \le k+l$.      
By Lemma \ref{Lemma (3.7)}, (1), the map $\mu$ is an isomorphism over $O_{\mathbf p} \cup 
O_1 \cup ... \cup O_k$, and $\mu$ induces a crepant resolution of an open neighborhood of $O_{\bar{\mathbf p}} \subset \bar{O}_{\mathbf p}$.  Write $$Spin(m) \times^Q (\mathfrak{n} + \bar{O}_{\mathbf{p}'}) 
= \mu^{-1}(O_{\mathbf p}) \sqcup \mu^{-1}({O}_{\bar{\mathbf p}}) 
\sqcup \mu^{-1}({O}_1) \sqcup ... \sqcup \mu^{-1}({O}_k) \sqcup \mu^{-1}({O}_{k+1}) \sqcup ... \sqcup \mu^{-1}({O}_{k+l})$$ 
Then $\mu^{-1}(O_1)$, ..., $\mu^{-1}(O_k)$ have codimension 2 in $Spin(m) \times^Q (\mathfrak{n} + \bar{O}_{\mathbf{p}'})$. 
On the other hand, since $\mathrm{Codim}_{O_{\mathbf p}}O_{\alpha}  \geq 4$ for 
$k+1 \le \alpha \le k+l$, the inverse images $\mu^{-1}(O_{k+1})$, ..., $\mu^{-1}(O_{k+l})$ have codimension $\geq 2$ by 
Corollary \ref{Corollary (0.2)}. This means that $\mu^{-1}(O_{\mathbf p} 
\cup O_{\bar{\mathbf p}})$ is obtained from $Spin(m) \times^Q (\mathfrak{n} + \bar{O}_{\mathbf{p}'})$ by removing a closed 
subset of codimension $\geq 2$.   
There is an isomorphism $$\mu_* : 
\pi_1(\mu^{-1}(O_{\mathbf p} 
\cup O_{\bar{\mathbf p}})) \to \pi_1(O_{\mathbf p} 
\cup O_{\bar{\mathbf p}})$$ by \cite{Ko}, Theorem 7.8 because $O_{\mathbf p} 
\cup O_{\bar{\mathbf p}}$ has only quotient singularities. 
Since $\mu^{-1}(O_{\mathbf p} 
\cup O_{\bar{\mathbf p}})$ is smooth, it is contained in $Spin(m) \times^Q (\mathfrak{n} + O_{\mathbf{p}'})$. 
Therefore, $\mu^{-1}(O_{\mathbf p} 
\cup O_{\bar{\mathbf p}})$ is obtained from the smooth variety $Spin(m) \times^Q (\mathfrak{n} + O_{\mathbf{p}'})$ by removing a closed subset of codimension at least $2$.   
Hence we have  $$\pi_1(\mu^{-1}(O_{\mathbf p} 
\cup O_{\bar{\mathbf p}})) \cong \pi_1(Spin(m) \times^Q (\mathfrak{n} + O_{\mathbf{p}'})).$$ 
By the exact sequence 
$$\pi_1(\mathfrak{n} + O_{\mathbf{p}'}) \to \pi_1(Spin(m) \times^Q (\mathfrak{n} + O_{\mathbf{p}'})) \to \pi_1(Spin(m)/Q) \to 1$$ we 
see that $\pi_1(Spin(m) \times^Q (\mathfrak{n} + O_{\mathbf{p}'}))$ has order 
at most $2^{a(\mathbf{p}) -2}$ because $\pi_1(O_{\mathbf{p}'})$ has order 
$2^{a(\mathbf{p}) -2}$. Let us consider the finite covering $\tau$. By definition, 
$\mathrm{deg}(\tau) = 2^{a(\mathbf{p}) -1}$. Since $S$ is an $A_1$-surface singularity, 
there are two possibilities. The first case is when $\tau^{-1}(S)$ is a disjoint union of $2^{a(\mathbf{p}) -1}$ copies of an $A_1$-surface singularity. The second case is 
when $\tau^{-1}(S)$ is a disjoint union of $2^{a(\mathbf{p}) -2}$ copies of $(\mathbf{C}^2,0)$. We only have to show that the first case does not occur. 
Actually, if the first case occurs, $\tau$ induces an etale covering of $O_{\mathbf p} 
\cup O_{\bar{\mathbf p}}$. Since $\mathrm{deg}(\tau) = 2^{a(\mathbf{p}) -1}$,  
this implies that $\vert \pi_1(O_{\mathbf p} \cup O_{\bar{\mathbf p}})\vert = 2^{a(\mathbf{p}) -1}$. This is a contradiction. 

Next assume that $\bar{\mathbf p}$ is very even. In this case, there are two 
nilpotent orbits $O_{\bar{\mathbf p}}^{\pm}$ with Jordan type $\bar{\mathbf p}$. We can write 
$$\mathrm{Sing}(\bar{O}_{\mathbf p}) = \bar{O}_{\bar{\mathbf p}}^{+} \cup \bar{O}_{\bar{\mathbf p}}^{-} 
\cup \bar{O}_1 \cup ... \cup \bar{O}_k \cup \bar{O}_{k+1} \cup ... \cup \bar{O}_{k+l}$$ Here $\mathrm{Codim}_{\bar{O}_{\mathbf p}}O_{\alpha} 
= 2$ for $1 \le \alpha \le k$ and $\mathrm{Codim}_{\bar{O}_{\mathbf p}}O_{\alpha} 
\geq 4$ for $k+1 \le \alpha \le k +l$. 
By the same reasoning as the case when $\bar{\mathbf p}$ is not very even, 
we see that $\pi_1(O_{\mathbf p} \cup O_{\bar{\mathbf p}}^{+} \cup 
O_{\bar{\mathbf p}}^{-})$ has order at most $2^{a(\mathbf{p}) -2}$. 
Let $S^{\pm}$ be respectively transverse slices for $O_{\bar{\mathbf p}}^{\pm} 
\subset \bar{O}_{\mathbf p}$. The adjoint action of $O(m)$ interchanges 
$O^{+}_{\bar{\mathbf p}}$ and $O^{-}_{\bar{\mathbf p}}$. Since $\tau$ is 
$O(m)$-equivariant by Lemma \ref{Lemma (3.3)}, $\tau^{-1}(S^+)$ and $\tau^{-1}(S^-)$ have the same 
splitting type. Therefore, there are two possibilities. The first case is 
when $\tau^{-1}(S^{\pm})$ are both disjoint unions of $2^{a(\mathbf{p}) -1}$ copies of an $A_1$-surface singularity. The second case is 
when $\tau^{-1}(S^{\pm})$ are both disjoint unions of $2^{a(\mathbf{p}) -2}$ copies of $(\mathbf{C}^2,0)$. Then we see that the first case does not occur by the same reasoning as in  the case $\bar{\mathbf p}$ is not very even. This completes the proof of (1). 
\vspace{0.2cm}

(2)  Assume that $i$ is a gap member of $\mathbf{p}$ except the case (1).  
Then $\mathbf{p}'$ is also rather odd. 
Moreover, $a(\mathbf{p}') = a(\mathbf{p})$.   
When $\mathbf{p}$ is very even, $\mathbf{p}'$ is also very even. 
In this case there are two orbits with Jordan type $\mathbf{p}'$. But there 
is a unique nilpotent orbit $O_{\mathbf{p}'}$ with Jordan type $\mathbf{p}'$ such that $O_{\mathbf p}$ is induced from $O_{\mathbf{p}'}$. 
Let $X_{\mathbf{p}'} \to 
\bar{O}_{\mathbf{p}'}$ be the finite covering associated with the universal covering 
of $O_{\mathbf{p}'}$. Then the $Q$-action on $\mathfrak{n} + \bar{O}_{\mathbf{p}'}$ 
does not lift to a $Q$-action on $\mathfrak{n} + X_{\mathbf{p}'}$ 
by Claim \ref{Claim (3.5.2)}. Instead, we 
take a cover $\tau': Y_{\mathbf{p}'} \to \bar{O}_{{\mathbf p}'}$ corresponding to 
the surjection $\pi_1(O_{{\mathbf p}'}) \to (\mathbf{Z}/2\mathbf{Z})^{\oplus \mathrm{max}(a(\mathbf{p}')-1, 0)}$ in 
Proposition \ref{Proposition (3.1)}, (1). Then 
the $Q$-action on $\mathfrak{n} + \bar{O}_{\mathbf{p}'}$ lifts to a 
$Q$-action on $\mathfrak{n} + Y_{\mathbf{p}'}$.  
Therefore we have a 
commutative diagram 
\begin{equation} 
\begin{CD} 
Spin(m) \times^{Q}(\mathfrak{n}+ Y_{\mathbf{p}'}) @>{\mu'}>> Y'\\ 
@V{\pi'}VV @V{\tau'}VV \\ 
Spin(m) \times^{Q}(\mathfrak{n} + \bar{O}_{\mathbf{p}'}) @>{\mu}>> \bar{O}_{\mathbf p}      
\end{CD} 
\end{equation}
Here $Y'$ is the Stein factorization of $\mu \circ \pi'$.  
$\pi$ factorizes as $$X \stackrel{\rho'}\to Y' \stackrel{\tau'}\to \bar{O}_{\mathbf p},$$ where $\mathrm{deg}(\rho') = 2$ and $\mathrm{deg}(\tau') = 2^{\mathrm{max}(a(\mathbf{p}) - 1,0)}$. Notice that $Q = \rho_m^{-1}(\bar{Q})$ for a parabolic subgroup 
$\bar{Q} \subset SO(m)$. Hence  $Spin(m) \times^{Q}(\mathfrak{n}+ Y_{\mathbf{p}'})
= SO(m) \times^{\bar Q}(\mathfrak{n}+ Y_{\mathbf{p}'})$. In particular, $SO(m)$ acts on $Y'$ and $\tau'$ is a 
$SO(m)$-finite cover. Since $\mathrm{deg}(\tau') = \mathrm{deg}(\tau)$, we see that $\tau' = \tau$ and $Y' = Y$.

(2-i):  We only have to prove the last statement. 
The transverse slice $S$ is an $A_1$-surface singularity. By Lemma \ref{Lemma (3.7)}, (1), $\mu$ induces a crepant resolution of an open 
neighborhood of $O_{\bar{\mathbf p}} \subset \bar{O}_{\mathbf p}$. 
Put $E:= \mu^{-1}(O_{\bar{\mathbf p}})$. Then $E$ is a divisor and 
$Spin(m) \times^{Q}(\mathfrak{n} + \bar{O}_{\mathbf{p}'})$ is smooth around 
$E$. Since $\mu'$ is a crepant partial resolution, $\pi'$ must be unramified at 
$E$. In fact, suppose that $\pi'$ is ramified around $E$. Let us consider the commutative diagram above over an open neighborhood 
of $O_{\bar{\mathbf p}} \subset \bar{O}_{\mathbf p}$. Put $E' := (\pi')^{-1}(E)$. 
Then $K_{Spin(m) \times^{Q}(\mathfrak{n}+ Y_{\mathbf{p}'})} = 
(\pi')^*K_{Spin(m) \times^{Q}(\mathfrak{n} + \bar{O}_{\mathbf{p}'})} + rE'$ for some $r > 0$. Since 
$Spin(m) \times^{Q}(\mathfrak{n} + \bar{O}_{\mathbf{p}'}) = \mu^*K_{\bar{O}_{\mathbf p}}$, 
we have $$K_{Spin(m) \times^{Q}(\mathfrak{n}+ Y_{\mathbf{p}'})} = (\pi')^*\mu^*K_{\bar{O}_{\mathbf p}} + rE'.$$ 
On the other hand, since $\mu'$ is crepant, $$K_{Spin(m) \times^{Q}(\mathfrak{n}+ Y_{\mathbf{p}'})} = (\mu')^*K_Y 
= (\mu')^*\tau^*K_{\bar{O}_{\mathbf p}} = (\pi')^*\mu^*K_{\bar{O}_{\mathbf p}}.$$ This is a contradiction; hence, 
$\pi'$ is unramified along $E$.  

This means that $\tau^{-1}(S)$ is a disjoint union of $2^{\mathrm{max}(a(\mathbf{p}) - 1,0)}$ copies of an $A_1$-surface singularity.
 
(2-ii) :  
If $\mathbf{p}$ is very even, then $\mathrm{deg}(\tau) = 1$. Hence $\tau^{-1}(S)$ is an $A_1$-surface singularity. 
Assume that $\mathbf{p}$ is not very even. 
Then $S$ is a union of two $A_1$-surface singularities.  
Let $\nu: \tilde{O}_{\mathbf p} \to \bar{O}_{\mathbf p}$ be the normalization. 
$SO(m)$ naturally acts on $\tilde{O}_{\mathbf p}$. 
Then $\nu^{-1}(O_{\bar{\mathbf p}}) \to O_{\bar{\mathbf p}}$ is an etale double cover.  Moreover, $\nu^{-1}(O_{\bar{\mathbf p}})$ is connected by Lemma \ref{Lemma (3.8)}.  Hence $SO(m)$ acts on $\nu^{-1}(O_{\bar{\mathbf p}})$ transitively.  
$\tau$ induces a finite covering $\tilde{\tau}: Y \to \tilde{O}_{\mathbf p}$. 
Then $\tilde{S} := \nu^{-1}(S)$ is a disjoint union of two $A_1$-surface 
singularities: $\tilde{S} = \tilde{S}^+ \sqcup \tilde{S}^{-}$.  
We may assume that $\tilde{S}^+$ and $\tilde{S}^{-}$ are interchanged each 
other by a suitable element of $SO(m)$. 

Let us consider ${\tau}^{-1}(S)$. Notice that ${\tau}^{-1}(S) = 
\tilde{\tau}^{-1}(\tilde{S}^+) \sqcup \tilde{\tau}^{-1}(\tilde{S}^-)$.
By Lemma \ref{Lemma (3.7)}, (i), $\mu$ induces a crepant resolution of an open 
neighborhood of $O_{\bar{\mathbf p}} \subset \bar{O}_{\mathbf p}$. 
Put $E:= \mu^{-1}(O_{\bar{\mathbf p}})$. Then $E$ is a divisor and 
$Spin(m) \times^{Q}(\mathfrak{n} + \bar{O}_{\mathbf{p}'})$ is smooth around 
$E$. Since $\mu'$ is a crepant partial resolution, $\pi'$ must be unramified at 
$E$ by the same argument as in (2-i). This means that $\tilde{\tau}^{-1}(\tilde{S}^+)$ and $\tilde{\tau}^{-1}(\tilde{S}^-)$ 
are respectively disjoint unions of $2^{\mathrm{max}(a(\mathbf{p}) - 1,0)}$ copies of an 
$A_1$-surface singularity.  

(3) The transverse slice $S$ for $O_{\bar{\mathbf p}} \subset 
\bar{O}_{\mathbf p}$ is an $A_3$-surface singularity. 
By Lemma \ref{Lemma (3.7)}, $\mu^{-1}(S)$ is a crepant partial  resolution of $S$  
with one exceptional curve $C \cong \mathbf{P}^1$. $\mu^{-1}(S)$ 
has $A_1$-surface singularities at two points $p^{\pm} \in C$. 
We prove that ${\pi'}^{-1}(\mu^{-1}(S))$ is a disjoint union of $2^{a(\mathbf{p})-2}$ 
copies of the minimal resolution of an $A_1$-surface singularity. This would mean 
that $\tau^{-1}(S)$ is a disjoint union of $2^{a(\mathbf{p})-2}$ 
copies of $A_1$-surface singularity. For this purpose we must look at the finite cover 
$\tau': Y_{\mathbf{p}'} \to \bar{O}_{{\mathbf p}'}$ induced from the $SO$-universal covering 
of $O_{{\mathbf p}'}$. 
Recall that 
$$\mathbf{p}' = [(d-2)^{r_d}, ..., (i-1)^{r_{i +1}}, i -2, i-4, (i-5)^{r_{i-5}}, ..., 1^{r_1}].$$
The nilpotent orbit closure  $\bar{O}_{\mathbf{p}'}$ contains nilpotent orbits of 
Jordan type $(\bar{{\mathbf p}'})_j $. By definition 
$$(\bar{{\mathbf p}'})_j = [(d-2)^{r_d}, ..., (i-1)^{r_{i +1}}, (i -3)^2, (i -5)^{r_{i -5}}, ..., 1^{r_1}].$$ 
Notice that we write $i$  for $i_j$ and the subscript of $(\bar{{\mathbf p}'})_j$ indicates this $j$.  

The partition $(\bar{{\mathbf p}'})_j $ may possibly be very even. In such a case, there are two orbits with Jordan type  $(\bar{{\mathbf p}'})_j $. If $(\bar{{\mathbf p}'})_j $ is not very even, such an orbit is unique. In any case, let 
$O_{(\bar{{\mathbf p}'})_j }$ be one of such orbits, and let $T$ be a transverse slice for $O_{(\bar{{\mathbf p}'})_j } 
\subset \bar{O}_{\mathbf{p}'}$. $T$ is an $A_1$-surface singularity. 
By applying Lemma (3.9) to $\tau'$, we see that 
$(\tau')^{-1}(T)$ is a disjoint union of some copies of $(\mathbf{C}^2,0)$. 
Look at the map ${\pi'}^{-1}(\mu^{-1}(S)) \to \mu^{-1}(S)$ induced by 
$\pi'$. Recall that there are two points $p^{\pm} \in C$ such that $\mu^{-1}(S)$ 
has $A_1$-surface singularities at these points. By the observation above, this map 
is ramified at $p^{\pm}$. This means that ${\pi'}^{-1}(\mu^{-1}(S))$ is a disjoint union of $2^{a(\mathbf{p})-2}$ 
copies of the minimal resolution of an $A_1$-surface singularity. $\square$ 
\vspace{0.2cm}

{\em Proof of Proposition \ref{Proposition (3.6)}}. 
The map $\pi$ factors through $Y$: 
$$X \stackrel{\rho}\to Y \stackrel{\tau}\to \bar{O}_{\mathbf p}.$$ 
By Proposition \ref{Proposition (3.9)}, $\pi^{-1}(S_j)$ is a disjoint union of the copies of 
$(\mathbf{C}^2, 0)$ or a disjoint union of the copies of $A_1$-surface singularity.  We prove that 
the latter case does not happen. 
 
(a) Assume that Case (1) in Proposition \ref{Proposition (3.9)} occurs. 
By the proposition, $\tau^{-1}(S_j)$ is already a disjoint union of $(\mathbf{C}^2, 0)$. Then 
$\rho$ is etale over an open neighborhood of $\tau^{-1}(S_j)$; hence $\pi^{-1}(S_j)$ is a disjoint union of the copies of 
$(\mathbf{C}^2, 0)$. 

(b) Next consider one of the following cases of Proposition 3.9:  

Case (2-i), 

Case (2-ii)  and $\mathbf{p}$ is very even, 

Case (3). 

When $\mathbf{p}$ is not very even in Case (2-ii), we will need an additional care. Such a case will be treated in (c). 

Since $\pi$ is a Galois cover, there are two possibilities for 
$\pi^{-1}(S_j)$. The first case is when $\pi^{-1}(S_j)$  is a disjoint union of the copies of $(\mathbf{C}^2,0)$. The second case is when 
$\pi^{-1}(S_j)$ is a disjoint union of the copies of an $A_1$-surface singularity.
If the second case occurs, then 
$\rho$ induces an etale cover of $\tau^{-1}(O_{\mathbf p} \cup O_{\bar{\mathbf p}_j})$ 
of degree 2. In particular, $\pi_1(\tau^{-1}(O_{\mathbf p} \cup O_{\bar{\mathbf p}_j})) 
\ne \{1\}$. Let $\tau_j^0: Y_{\mathbf{p}'_j} \to O_{\mathbf{p}'_j}$ be the $SO$-universal covering and 
let $\tau_j: Y_{\mathbf{p}'_j} \to \bar{O}_{\mathbf{p}'_j}$ be its associated covering. 
Look at the commutative diagram 
\begin{equation} 
\begin{CD} 
Spin(m) \times^{Q_j}(\mathfrak{n}_j + Y_{\mathbf{p}'_j}) @>{\mu'_j}>> Y \\ 
@V{\pi'}VV @V{\tau}VV \\ 
Spin(m) \times^{Q_j}(\mathfrak{n}_j + \bar{O}_{\mathbf{p}'_j}) @>{\mu_j}>> \bar{O}_{\mathbf p}      
\end{CD} 
\end{equation}
The birational map 
$$(\mu'_j)^{-1}(\tau^{-1}(O_{\mathbf p} \cup O_{\bar{\mathbf p}_j})) \to \tau^{-1}(O_{\mathbf p} \cup O_{\bar{\mathbf p}_j})$$ induces an isomorphism of the fundamental groups of both sides by \cite{Ko}, Theorem 7.8 
because $\tau^{-1}(O_{\mathbf p} \cup O_{\bar{\mathbf p}_j})$ 
has only quotient singularities.    
We next calculate $\pi_1((\mu'_i)^{-1}(\tau^{-1}(O_{\mathbf p} \cup O_{\bar{\mathbf p}_j})))$.
The closure $\bar{O}_{\mathbf p}$ contains $O_{\mathbf{p}}$ as an open dense orbit and contains 
$O_{\bar{\mathbf p}_1}$, ..., $O_{\bar{\mathbf p}_k}$ as codimension 2 orbits. Other nilpotent orbits in $\bar{O}_{\mathbf p}$ 
have codimension $\geq 4$. Therefore, one can write   
$$\bar{O}_{\mathbf p} = {O}_{\mathbf p} \sqcup O_{\bar{\mathbf p}_1} \sqcup ... \sqcup O_{\bar{\mathbf p}_j} \sqcup ... \sqcup
O_{\bar{\mathbf p}_k} \sqcup F,$$ where $F$ is the union of all nilpotent orbits with codimension $\geq 4$. 
Hence $$Y = \tau^{-1}({O}_{\mathbf p}) \sqcup \tau^{-1}(O_{\bar{\mathbf p}_1}) \sqcup ... \sqcup \tau^{-1}(O_{\bar{\mathbf p}_j}) \sqcup ... \sqcup \tau^{-1}(O_{\bar{\mathbf p}_k}) \sqcup \tau^{-1}(F).$$
For $l \ne j$, the birational map $\mu_j$ is an isomorphism over an open neighborhood of $O_{\bar{\mathbf p}} \subset 
\bar{O}_{\bar{\mathbf p}}$ by Lemma 3.7, (1). By Zariski's Main Theorem, $\mu'_j$ is also an isomorphism over an 
open neighborhood of $\tau^{-1}(O_{\bar{\mathbf p}_l}) \subset Y$. Moreover, $\mu'_j$ gives a crepant resolution of an open 
neighborhood of $\tau^{-1}(O_{\bar{\mathbf p}_j}) \subset Y$. Here write 
$$Spin(m) \times^{Q_j}(\mathfrak{n}_j + Y_{\mathbf{p}_j}) =  
(\tau \circ \mu'_j)^{-1}({O}_{\mathbf p}) \sqcup (\tau \circ \mu'_j)^{-1}(O_{\bar{\mathbf p}_1}) \sqcup ... \sqcup (\tau \circ \mu'_j)^{-1}(O_{\bar{\mathbf p}_k}) 
\sqcup (\tau \circ \mu'_j)^{-1}(F).$$
Then $(\tau \circ \mu'_j)^{-1}(O_{\bar{\mathbf p}_l})$ has codimension 2 in   
$Spin(m) \times^{Q_j}(\mathfrak{n}_j + Y_{\mathbf{p}_j})$ for $l \ne j$. Since $\mathrm{Codim}_Y\tau^{-1}(F) \geq 4$, we see that 
$(\tau \circ \mu'_j)^{-1}(F)$ has codimension $\geq 2$ in $Spin(m) \times^{Q_j}(\mathfrak{n}_j + Y_{\mathbf{p}_j})$ by Corollary \ref{Corollary (0.2)}.  
Therefore, $(\mu'_j)^{-1}(\tau^{-1}(O_{\mathbf p} \cup O_{\bar{\mathbf p}_j}))$ is obtained from 
$Spin(m) \times^{Q_j}(\mathfrak{n}_j + Y_{\mathbf{p}_j})$ by removing the subset of codimension $\geq 2$. 
Let $Y^{\mathrm{reg}}_{\mathbf{p}'_j}$ be the smooth part of $Y_{\mathbf{p}'_j}$. 
Since $(\mu'_j)^{-1}(\tau^{-1}(O_{\mathbf p} \cup O_{\bar{\mathbf p}_j}))$ is smooth, 
there is an inclusion  
$$(\mu'_j)^{-1}(\tau^{-1}(O_{\mathbf p} \cup O_{\bar{\mathbf p}_j}))  \subset Spin(m) \times^{Q_j}(\mathfrak{n}_j + Y^{\mathrm{reg}}_{\mathbf{p}'_j})$$
Hence, $(\mu'_j)^{-1}(\tau^{-1}(O_{\mathbf p} \cup O_{\bar{\mathbf p}_j})$ is obtained from a smooth variety $Spin(m) \times^{Q_j}(\mathfrak{n}_j + Y^{\mathrm{reg}}_{\mathbf{p}'_j})$ by removing a closed subset of codimension 
$\geq 2$. Therefore we have an isomorphism 
$$\pi_1((\mu'_j)^{-1}(\tau^{-1}(O_{\mathbf p} \cup O_{\bar{\mathbf p}_j}))) \cong 
\pi_1(Spin(m) \times^{Q_j}(\mathfrak{n}_j + Y^{\mathrm{reg}}_{\mathbf{p}'_j})).$$
Put $Y^0_{\mathbf{p}'_j} := \tau_j^{-1}(O_{\mathbf{p}'_j})$. Then,   
by Claim \ref{Claim (3.5.2)}, $\pi_1(Spin(m) \times^{Q_j}(\mathfrak{n}_j + Y^0_{\mathbf{p}'_j})) = \{1\}$.
Since there is a surjection $$\pi_1(Spin(m) \times^{Q_j}(\mathfrak{n}_j + Y^0_{\mathbf{p}'_j})) \to 
\pi_1(Spin(m) \times^{Q_j}(\mathfrak{n}_j + Y^{\mathrm{reg}}_{\mathbf{p}'_j})),$$  
we see that $\pi_1(Spin(m) \times^{Q_j}(\mathfrak{n}_j + Y^{\mathrm{reg}}_{\mathbf{p}'_j})) = \{1\}$. This is a contradiction. 
As a consequence, 
$\pi^{-1}(S_j)$  is a disjoint union of the copies of $(\mathbf{C}^2,0)$.

(c) We finally assume that Case (2-ii) occurs and $\mathbf{p}$ is not very even. Then $S_j$ is of type $A_1 \cup A_1$. Let $\nu: \tilde{O}_{\mathbf p} \to \bar{O}_{\mathbf p}$ be the normalization map. 
Then $\tilde{S}_j := \nu^{-1}(S_j)$ is a disjoint union of two $A_1$-surface singularities: 
$\tilde{S}_j = \tilde{S}_j^+ \sqcup \tilde{S}_j^-$. $\tilde{S}_j^+$ and $\tilde{S}_j^-$ are interchanged each other by a 
suitable element of $SO(m)$.   
Let us look at the double cover $\rho$. Since $\rho$ is $Spin(m)$-equivariant, there are two possibilities. The first case is when $(\tilde{\tau} \circ \rho)^{-1}(\tilde{S}^{\pm})$ are both disjoint unions of the copies of $(\mathbf{C}^2,0)$. The second case is when $(\tilde{\tau} \circ \rho)^{-1}(\tilde{S}^{\pm})$ are both disjoint unions of the copies of  
$A_1$-surface singularity.  Here $\tilde{\tau}$ is the map from $Y$ to the normalization $\tilde{O}_{\mathbf p}$ of 
$\bar{O}_{\mathbf p}$ induced from $\tau: Y \to \bar{O}_{\mathbf p}$. 
By the same argument as in (b), the second case does not occur. As a consequence, $\pi^{-1}(S_j)$ is 
a disjoint union of the copies of $(\mathbf{C}^2,0)$.  $\square$ 
\vspace{0.2cm}

{\bf Construction of a Q-factorial terminalization}

Let $\mathbf{p}$ be a partition of $m$ with Condition (i). Let $O_{\mathbf p} \subset so(m)$ be a nilpotent orbit with Jordan type $\mathbf{p}$, and let $X \to \bar{O}_{\mathbf p}$ be the finite covering associated with the universal covering 
of $O_{\mathbf p}$. We are now in a position to construct a {\bf Q}-factorial terminalization of $X$. \vspace{0.2cm}

{\bf (A) The case when $\mathbf{p}$ is rather odd}.  

We will use the following induction step. 

(Double Induction of type I)  Let $i$ be a gap member of $\mathbf{p}$. 
When $i$ is odd, we assume that $i \geq 5$ and $r_{i-1} = ... = r_{i-4} = 0$. 
When $i$ is even, we assume that $i \geq 4$ and $r_{i-1} = ... = r_{i-3} = 0$. 
We remark that such a gap member exists if $\mathbf{p}$ does not satisfy the 
condition (iii) (which is defined just before Proposition \ref{Proposition (3.6)}). 
Put $$ \mathbf{p}' = [(d-4)^{r_d}, ..., (i-3)^{r_{i+1}}, (i-4)^{r_i + r_{i-4}}, (i-5)^{r_{i-5}}, ..., 1^{r_1}].$$ Then $\mathbf{p}'$ satisfies (i) and is still rather odd.  
If we put $s := r_d + ... + r_i$, then $\mathbf{p}'$ is a partition of 
$m-4s$. Let $Q_{(2s, m-4s, 2s)} \subset SO(m)$ be a parabolic subgroup 
of flag type $(2s, m-4s, 2s)$ and put $Q := \rho_m^{-1}(Q_{(2s, m-4s, 2s)})$ for 
the double cover $\rho_m: Spin(m) \to SO(m)$.  
We write the Levi decomposition of $\mathfrak{q}$ as $\mathfrak{q} = 
\mathfrak{l} \oplus \mathfrak{n}$. We have  $\mathfrak{l} = \mathfrak{gl}(2s)  
\oplus so(m-4s)$. One can take a nilpotent orbit $O_{[2^s]} \times O_{{\mathbf p}'}$ in 
$\mathfrak{gl}(sr) \oplus so(m-4s)$ so that 
$O_{\mathbf p} = \mathrm{Ind}^{so(m)}_{\mathfrak l}(O_{[2^s]} \times O_{{\mathbf p}'})$. 
The generalized Springer map 
$$\mu: Spin(m) \times^Q(\mathfrak{n} + \bar{O}_{[2^s]} \times \bar{O}_{{\mathbf p}'}) \to \bar{O}_{\mathbf p}$$ 
is a birational map. We indicate this process simply by 
$$\mathbf{p} \stackrel{\mathrm{Type\: I}^2}\leftarrow \mathbf{p}'.$$

Recall that $a(\mathbf{p})$ is the number of distinct odd members of 
$\mathbf{p}$. Now assume that $\mathbf{p}$ does not satisfy the condition (iii).
Then we can repeat double induction step of type I to  $\mathbf{p}$ so that $a$ does not change, and finally get a partition $\mathbf{p}'$ 
with conditions (i) and (iii). 
$$\mathbf{p} \stackrel{\mathrm{Type\: I}^2}\leftarrow \mathbf{p}_1 \stackrel{\mathrm{Type\: I}^2}\leftarrow ... \stackrel{\mathrm{Type\: I}^2}\leftarrow  \mathbf{p}_k = \mathbf{p}' $$  
In each step, we record the number $s_i$ and put $q := m - 4\sum s_i$. 

{\bf Case 1}: $q \geq 3$ 
 
Put the $m \times m$ matrix 
$$J_m = \left(\begin{array}{ccccc} 
0 & 0 & 0 & ... & 1 \\
0 & 0 & ... & 1 & 0 \\ 
... & ... & ... & ... & ... \\ 
0 & 1 & ... & 0 & 0 \\ 
1 & 0 & ... & 0 & 0 
\end{array}\right). $$   
Then $$SO(m) = \{ A \in SL(m) \: \vert \: A^tJ A = J\}. $$
Now let us consider the isotropic flag of type $(2s_1, 2s_2, ..., 2s_k, q,  2s_k, ..., 2s_2, 2s_1)$:  
$$0 \subset \langle e_1, ..., e_{2s_1} \rangle \subset \langle e_1, ... e_{2s_1 + 2s_2} \rangle \subset ... \subset  \langle e_1, ..., e_{\sum_{1}^{k}2s_i} \rangle \subset 
\langle e_1, ..., e_{\sum_{1}^{k}2s_i + q} \rangle$$
$$\subset \langle e_1, ..., e_{\sum_{1}^{k}2s_i + q +2s_k} \rangle \subset ... \subset 
\langle e_1, ..., e_{\sum_{1}^{k}2s_i + q +\sum_{1}^{k}2s_i} \rangle = \mathbf{C}^m $$ 
Let $Q'$ be 
the parabolic subgroup of $SO(m)$ stabilizing this flag.  
One has a Levi decomposition $Q' = U' \cdot L'$ with 
$$ L' = \{ \left(\begin{array}{ccccccccc} 
A_1 & 0 & 0 & ... & ... & ... & ... & ... & 0 \\
0 & A_2 & 0 & ... & ... & ... & ... & ... & 0 \\ 
... & ... & ... & ... & ... & ... & ... & ... & ...\\
0 & ... & 0 & A_k & 0 & ... & ... & ... & 0\\
0 & ... & ... & 0 & B & 0 & ... & ... & 0 \\
0 & ... & ... & ... & 0 & A'_k & 0 & ... & 0\\
... & ... & ... & ... & ... & ... & ... & ... & ... \\ 
0 & ... & ... & ... & ... & ... & ... & A'_2 & 0 \\ 
0 & ... & ... & ... & ... & ... & ... & 0 & A_1' 
\end{array}\right) \: \vert \: A_i \in GL(2s_i), \: A'_i = J_{2s_i}(A_i^t)^{-1}J_{2s_i}, \: 
B \in SO(q)\}.$$ In particular, $L' \cong \prod GL(2s_i) \times SO(q)$.  
We define a parabolic subgroup $Q$ of $Spin(m)$ by 
$Q := \rho_m^{-1}(Q')$ for $\rho_m: Spin(m) \to SO(m)$.  
Then the Levi part $L$ of $Q$ is a double cover of $L$, which is described as 
follows. There is a natural map $SL(2s_i) \times \mathbf{C}^* \to GL(2s_i)$ defined by 
$(X, \lambda) \to \lambda X$. The kernel of this map is the cyclic group  
$$\mu_{2s_i} := \{ (\lambda^{-1}I_{2s_i}, \lambda)\: \vert \: \lambda^{2s_i} = 1\}.$$
Let $\mu_{s_i}$ be the subgroup of $\mu_{2s_i}$ of order $s_i$.    
Then the covering map $SL(2s_i) \times \mathbf{C}^*/\mu_{s_i} \to GL(2s_i)$ 
has degree $2$. We then have a covering map 
$$(SL(2s_1) \times \mathbf{C}^*)/\mu_{s_1} \times ... \times   
(SL(2s_k) \times \mathbf{C}^*)/\mu_{s_k} \times Spin(q) \to GL(2s_1) \times ... 
\times GL(2s_k) \times SO(q).$$ The Galois group of this covering is 
$(\mathbf{Z}/2\mathbf{Z})^{\oplus k +1}$. 
Put $H := \mathrm{Ker}[(\mathbf{Z}/2\mathbf{Z})^{\oplus k +1} \stackrel{\sum}\to \mathbf{Z}/2\mathbf{Z}]$, where $\sum$ is defined by $(x_1, ..., x_{k+1}) \to \sum x_i$. 
By Claim \ref{Claim (3.5.1)} we have  
$$L =  \{(SL(2s_1) \times \mathbf{C}^*)/\mu_{s_1} \times ... \times   
(SL(2s_k) \times \mathbf{C}^*)/\mu_{s_k} \times Spin(q)\}\: /\: H.$$

The Levi part $L$ has the same Lie algebra of that of $L'$. Hence we have
$\mathfrak{l} = \oplus gl(2s_i) \oplus so(q)$. We consider the nilpotent orbit 
$\prod O_{[2^{s_i}]} \times O_{\mathbf{p}'}$ in $\mathfrak{l}$. By the construction 
$$O_{\mathbf p} = \mathrm{Ind}^{so(m)}_{\mathfrak l}(\prod O_{[2^{s_i}]} \times O_{\mathbf{p}'}).$$
Let $X_{[2^{s_i}]} \to \bar{O}_{[2^{s_i}]}$ be the double covering associated with 
the universal covering of $O_{[2^{s_i}]}$. By Proposition \ref{Proposition (1.8)}, (2) and Proposition \ref{Proposition (1.9)}, 
$X_{[2^{s_i}]}$ has only {\bf Q}-factorial terminal singularities. 
Let $X_{{\mathbf p}'} \to \bar{O}_{{\mathbf p}'}$ 
be the finite covering associated with the universal covering of $O_{{\mathbf p}'}$. 
By Proposition \ref{Proposition (3.6)} $X_{{\mathbf p}'}$ has only {\bf Q}-factorial terminal singularities. 
$(SL(2s_1) \times \mathbf{C}^*)/\mu_{s_1} \times ... \times (SL(2s_k) \times \mathbf{C}^*)/\mu_{s_k} \times Spin(q)$ acts on 
$\prod X_{[2^{s_i}]} \times X_{\mathbf{p}'}$ because 
each $SL(2s_i) \times \mathbf{C}^*/\mu_{s_i}$ acts on $X_{[2^{s_i}]}$, and 
$Spin(q)$ acts on $X_{\mathbf{p}'}$.  By Proposition \ref{Proposition (3.1)}, (1) we have an exact sequence $$1 \to \mathbf{Z}/2\mathbf{Z} \to \pi_1(O_{\mathbf{p}'}) \to (\mathbf{Z}/2\mathbf{Z})
^{\oplus \mathrm{max}(a(\mathbf{p}')-1, 0)} \to 1.$$ Then the generator of  $\mathbf{Z}/2\mathbf{Z}$ induces a covering involution $-1$ for $X_{\mathbf{p}'} \to 
\bar{O}_{\mathbf{p}'}$. Hence $(\mathbf{Z}/2\mathbf{Z})^{\oplus k+1}$ acts on 
$\prod X_{[2^{s_i}]} \times X_{\mathbf{p}'}$. Let $H$ be the subgroup of  
$(\mathbf{Z}/2\mathbf{Z})^{\oplus k+1}$ defined as above. The quotient space 
$\prod X_{[2^{s_i}]} \times X_{\mathbf{p}'}/H$ has only terminal singularities because the fixed locus of every nonzero element of $H$ has codimension $\geq 4$. Moreover, since  $\prod X_{[2^{s_i}]} \times X_{\mathbf{p}'}$ is {\bf Q}-factorial, $\prod X_{[2^{s_i}]} \times X_{\mathbf{p}'}/H$ is also {\bf Q}-factorial by Lemma \ref{Lemma (1.4)}. The Levi part 
$L$ acts on $(\prod X_{[2^{s_i}]} \times X_{\mathbf{p}'})/H$. This action is a lifting of the 
adjoint $L$-action on $\prod O_{[2^{s_i}]} \times O_{\mathbf{p}'}$. 
Let us consider the Levi decomposition $\mathfrak{q} = \mathfrak{l} \oplus 
\mathfrak{n}$. Then the observation above means that $\mathfrak{n} + (\prod X_{[2^{s_i}]} \times X_{\mathbf{p}'})/H$ is a $Q$-space. There is a finite cover 
$$\pi': Spin(m) \times^{Q}(\mathfrak{n} + (\prod X_{[2^{s_i}]} \times X_{\mathbf{p}'})/H)
\to Spin(m) \times^{Q}(\mathfrak{n} + \prod \bar{O}_{[2^{s_i}]} \times \bar{O}_{\mathbf{p}'}) $$ and $\mathrm{deg}(\pi) = \mathrm{deg}(\pi')$.  
We now have a commutative diagram 
\begin{equation} 
\begin{CD} 
Spin(m) \times^{Q}(\mathfrak{n} + (\prod X_{[2^{s_i}]} \times X_{\mathbf{p}'})/H) @>{\mu'}>> X \\ 
@V{\pi'}VV @V{\pi}VV \\ 
Spin(m) \times^{Q}(\mathfrak{n} + \prod \bar{O}_{[2^{s_i}]} \times \bar{O}_{\mathbf{p}'}) @>{\mu}>> \bar{O}_{\mathbf p}      
\end{CD} 
\end{equation} 
Here $X$ coincides with the Stein factorization of $\mu \circ \pi'$. 
The map $\mu'$ gives a {\bf Q}-factorial terminalization of $X$.
\vspace{0.2cm}

{\bf Case 2}: $q \leq 2$

When $q = 0$, $\mathbf{p}' = [0]$. In this case, $\mathbf{p}$ is very even. Hence $\pi: X \to \bar{O}_{\mathbf p}$ is a 
double covering.  The cases $q = 1, 2$ do not occur.    
We define $\sum : \mathbf{Z}/2\mathbf{Z}^{\oplus k} \to \mathbf{Z}/2\mathbf{Z}$ 
by $\sum (x_1, ..., x_k) := \sum x_i$ and put $H := \mathrm{Ker}(\sum)$. 
We have a commutative diagram 
\begin{equation} 
\begin{CD} 
Spin(m) \times^{Q}(\mathfrak{n} + (\prod X_{[2^{s_i}]}/H) @>{\mu'}>> X \\ 
@V{\pi'}VV @V{\pi}VV \\ 
Spin(m) \times^{Q}(\mathfrak{n} + \prod \bar{O}_{[2^{s_i}]}) @>{\mu}>> \bar{O}_{\mathbf p}      
\end{CD} 
\end{equation} 
Here $X$ coincides with the Stein factorization of $\mu \circ \pi'$. 
The map $\mu'$ gives a {\bf Q}-factorial terminalization of $X$.
\vspace{0.2cm}

{\bf (B) The case when $\mathbf{p}$ is not rather odd}

We will use two induction steps of type (I) and of type (II).

By using the inductions of type (I) repeatedly for $\mathbf{p}$, we can finally find a parabolic subgroup 
$Q$ of $Spin(m)$ and a nilpotent orbit $O_{\mathbf{p}'}$ of a Levi part 
$\mathfrak{l}$ of $\mathfrak{q}$ such that 

(a) $O_{\mathbf p} = \mathrm{Ind}^{\mathfrak{so}(m)}_{\mathfrak l}(O_{\mathbf{p}'}).$

(b) $a(\mathbf{p}') = a(\mathbf{p})$, and 

(c) $\mathbf{p}'$ satisfies the condition (ii).   

For the partition $\mathbf{p}'$ thus obtained, we let $e$ be the number of odd members 
$i$ of $\mathbf{p}'$ such that $r_i = 2$.
By using the inductions of type (II) repeatedly for $\mathbf{p}'$, we 
can finally find a parabolic subgroup 
$Q'$ of $Spin(m)$ and a nilpotent orbit $O_{\mathbf{p}''}$ of the Levi part 
$\mathfrak{l}'$ of $\mathfrak{q}'$ such that 

(a') $O_{{\mathbf p}'} = \mathrm{Ind}^{\mathfrak{sp}(2n')}_{{\mathfrak l}'}(O_{\mathbf{p}''}),$

(b') $a(\mathbf{p}'') = a(\mathbf{p}') - e$, 

(c') $\mathbf{p}''$ satisfies the condition (ii) and 

(d') the multiplicity of any odd member of $\mathbf{p}''$ is not $2$. 

All together, we get a generalized Springer map 
$$\mu'': Spin(m) \times^{Q''}(\mathfrak{n}'' + \bar{O}_{{\mathbf p}''}) \to \bar{O}_{\mathbf p}.$$ 

Hereafter we consider two cases separately. 

{\bf Case 1}:  $r_i  \geq 3$ for some odd member $i$ of 
$\mathbf{p}$ 

In this case, $\mathbf{p}''$ is not rather odd. By the assumption $a(\mathbf{p}) \geq 
e + 1$. 
Let $X'' \to \bar{O}_{{\mathbf p}''}$ be the finite covering associated with the universal covering of $O_{{\mathbf p}''}$. Then $X''$ has only terminal 
singularities by Proposition \ref{Proposition (3.4)}. On the other hand, $X$ is {\bf Q}-factorial by (d'). 
Since $\mathbf{p}''$ is not rather odd, 
the $Q''$-action on $\mathfrak{n}'' + \bar{O}_{\mathbf{p}''}$ lifts to a $Q''$-action 
on $\mathfrak{n}'' + X''$.
There is a finite cover 
$$\pi'': Spin(m) \times^{Q''}(\mathfrak{n}'' + X'') \to Spin(m) \times^{Q''}(\mathfrak{n}'' + \bar{O}_{{\mathbf p}''}).$$ 
We have $$\mathrm{deg}(\mu'') = 2^e, \:\:
\mathrm{deg}(\pi'') = 2^{\mathrm{max}(a(\mathbf{p}'') -1,0)} 
= 2^{\mathrm{max}(a(\mathbf{p}) -e -1,0)} = 2^{a(\mathbf{p})-e-1}$$ and 
$\mathrm{deg}(\pi) = 2^{a(\mathbf{p})-1}$. It can be checked that $\mathrm{deg}(\pi) 
= \mathrm{deg}(\pi'') \cdot \mathrm{deg}(\mu'')$. 
Therefore we have a commutative diagram  
 \begin{equation} 
\begin{CD} 
Spin(m) \times^{Q''}(\mathfrak{n}'' + X'') @>>> X \\ 
@V{\pi''}VV @V{\pi}VV \\ 
Spin(m) \times^{Q''}(\mathfrak{n}'' + \bar{O}_{{\mathbf p}''}) @>{\mu''}>> \bar{O}_{\mathbf p}      
\end{CD} 
\end{equation}
where $X$ is obtained as the Stein factorization of the map $\mu'' \circ \pi''$.  
The map $Spin(m) \times^{Q''}(\mathfrak{n}'' + X'') \to X$ here obtained is a {\bf Q}-factorial terminalization of $X$. \vspace{0.2cm}

{\bf Case 2}: $r_i \leq 2$ for any odd member $i$ of $\mathbf{p}$   

In this case $\mathbf{p}''$ is rather odd. 
If there is an odd $i$ such that $r_i = 1$, then $a(\mathbf{p}) \geq e +1 $. 
If $r_i = 2$ for all odd $i$, then $a(\mathbf{p}) = e \geq 1$. The last inequality 
holds because if $e = 0$, then $\mathbf{p}$ is rather odd and this is not our case.
Then we have a short exact sequence 
$$ 1 \to \mathbf{Z}/2\mathbf{Z} \to \pi_1(O_{\mathbf{p}''}) \to 
(\mathbf{Z}/2\mathbf{Z})^{\oplus \mathrm{max}(a(\mathbf{p}'')-1, 0)} \to 1.$$ 
Let $Y'' \to  \bar{O}_{{\mathbf p}''}$ be the finite covering determined by 
the surjection from $\pi_1(O_{\mathbf{p}''})$ to  
$(\mathbf{Z}/2\mathbf{Z})^{\oplus \mathrm{max}(a(\mathbf{p}'')-1, 0)}$.
Since $\mathbf{p}''$ has only odd, ordinary gap members, we can apply 
Proposition \ref{Proposition (3.9)}, (1) for each gap member of $\mathbf{p}''$ to conclude that 
$\mathrm{Codim}_{Y''}\mathrm{Sing}(Y'') \geq 4$. Since $\mathbf{p}''$ is rather odd, 
$Y''$ is {\bf Q}-factorial by Proposition \ref{Proposition (3.1)}, (2). Therefore $Y''$ has only {\bf Q}-factorial 
terminal singularities.  
Then the $Q''$-action on $\mathfrak{n}'' + \bar{O}_{\mathbf{p}''}$ lifts to a $Q''$-action 
on $\mathfrak{n}'' + Y''$. 
There is a finite cover 
$$\tau'': Spin(m) \times^{Q''}(\mathfrak{n}'' + Y'') \to Spin(m) \times^{Q''}(\mathfrak{n}'' + \bar{O}_{{\mathbf p}''}).$$  
We have 
$$\mathrm{deg}(\mu'') = \left\{ 
\begin{array}{rl}    
2^e & \quad (r_i = 1 \:\: \mathrm{for}\: \mathrm{some}\: \mathrm{odd}\: i )\\ 
2^{e-1}& \quad (r_i = 2 \:\: \mathrm{for}\: \mathrm{all}\: \mathrm{odd}\: i) 
\end{array}\right.$$
$$\mathrm{deg}(\tau'') = 2^{\mathrm{max}(a({\mathbf p}'')-1,0)} 
= 2^{\mathrm{max}(a(\mathbf{p})-e-1,0)}$$ and 
$\mathrm{deg}(\pi) = 2^{\mathrm{max}(a({\mathbf p})-1,0)}$. 
It can be checked that $\mathrm{deg}(\pi) = \mathrm{deg}(\tau'') \cdot 
\mathrm{deg}(\mu'')$ in any case.  
Therefore we have a commutative diagram  
 
\begin{equation} 
\begin{CD} 
Spin(m) \times^{Q''}(\mathfrak{n}'' + Y'') @>>> X \\ 
@V{\tau''}VV @V{\pi}VV \\ 
Spin(m) \times^{Q''}(\mathfrak{n}'' + \bar{O}_{{\mathbf p}''}) @>{\mu''}>> \bar{O}_{\mathbf p}      
\end{CD} 
\end{equation}
where $X$ is obtained as the Stein factorization of the map $\mu'' \circ \tau''$.  
The map $Spin(m) \times^{Q''}(\mathfrak{n}'' + Y'') \to X$ here obtained is a {\bf Q}-factorial terminalization of $X$. 

\begin{Exams}\label{Examples (3.10)} \end{Exams} (1) Let us construct a {\bf Q}-factorial terminalization of 
the finite covering $X \to \bar{O}_{[15, 8^2, 3]} \subset so(34)$ associated with the universal covering of 
$O_{[15, 8^2, 3]}$. By double inductions of type I 
$$ [15, 8^2, 3]  \stackrel{Type I^2}\leftarrow [11, 8^2, 3] \stackrel{Type I^2}\leftarrow [7, 4^2, 3]$$ 
we get a partition $[7, 4^2, 3]$ of $18$ with conditions (i) and (iii). 
Let $\bar{Q} \subset SO(34)$ be a parabolic subgroup stabilizing an 
isotropic flag of type $(2,6,18,6,2)$ and put $Q := \rho_{34}^{-1}(\bar{Q})$ for 
$\rho_{34}: Spin (34) \to SO(34)$. Let $X_{[2]}$ (resp. $X_{[2^3]}$) be a finite 
covering of $\bar{O}_{[2]} \subset sl(2)$ (resp. $\bar{O}_{[2^3]} \subset sl(6)$) 
associated with the universal covering of $O_{[2]}$ (resp. $O_{[2^3]}$). 
Moreover, let $X_{[7, 4^2, 3]}$ be the finite covering of $\bar{O}_{[7, 4^2, 3]} \subset 
so(18)$ associated with the universal covering of $O_{[7, 4^2, 3]}$. 
Then $$Spin (34) \times^Q (\mathfrak{n} + (X_{[2]} \times X_{[2^3]} \times X_{[7, 4^2,3]})/H)$$ is a {\bf Q}-factorial terminalization of $X$. Here $H := 
\mathrm{Ker}[\sum : (\mathbf{Z}/2\mathbf{Z})^{\oplus 3} \to \mathbf{Z}/2\mathbf{Z}]$. 

(2) Let $X$ be the finite covering of $\bar{O}_{[11^3, 3^2, 1]} \subset so(40)$ 
associated with the universal covering of $O_{[11^3, 3^2, 1]}$. 
By inductions of type I and of type II 
$$[11^3, 3^2, 1] \stackrel{Type I}\leftarrow [9^3, 3^2, 1] \stackrel{Type I}\leftarrow  
[7^3, 3^2, 1] \stackrel{Type I}\leftarrow [5^3, 3^2, 1] \stackrel{Type II}\leftarrow  
[3^3, 2^2, 1]$$ we finally get a partition $[3^3, 2^2, 1]$ of $14$. 
Let $\bar{Q} \subset SO(40)$ be a parabolic subgroup stabilizing an 
isotropic flag of type $(3,3,3,4,14,4,3,3,3)$ and put $Q := \rho_{40}^{-1}(\bar{Q})$ for 
$\rho_{40}: Spin (40) \to SO(40)$. Let us consider the nilpotent orbit 
$O_{[3^3, 2^2, 1]} \subset so(14)$. 
Let $X_{[3^3,2^2,1]} \to \bar{O}_{[3^3, 2^2, 1]}$ be the finite covering associated with 
the universal covering of $O_{[3^3, 2^2, 1]}$.  
Then $$Spin (40) \times^Q (\mathfrak{n} + X_{[3^3,2^2,1]})$$ is a {\bf Q}-factorial 
terminalization of $X$. 

(3) Let $X$ be the finite covering of $\bar{O}_{[13^2, 3, 1]} \subset so(30)$ associated 
with the universal covering of $O_{[13^2, 3, 1]}$. 
By inductions of type I and of type II 
$$[13^2, 3, 1] \stackrel{Type I}\leftarrow 
[11^2, 3, 1] \stackrel{Type I}\leftarrow [9^2, 3, 1] \stackrel{Type I}\leftarrow 
[7^2, 3, 1] \stackrel{Type I}\leftarrow [5^2, 3, 1] \stackrel{Type II}\leftarrow 
[4^2, 3, 1]$$ we finally get a partition $[4^2, 3, 1]$ of $12$. 
Let $\bar{Q} \subset SO(30)$ be a parabolic subgroup stabilizing an 
isotropic flag of type $(2,2,2,2,1, 12, 1, 2,2,2,2)$ and put $Q := \rho_{30}^{-1}(\bar{Q})$ for 
$\rho_{30}: Spin (30) \to SO(30)$.
Let us consider the nilpotent orbit $O_{[4^2, 3, 1]} \subset so(12)$.
Since the partition $[4^2, 3, 1]$ is rather odd,  there is an exact sequence 
$$1 \to \mathbf{Z}/2\mathbf{Z} \to \pi_1(O_{[4^2, 3, 1]}) \to \mathbf{Z}/2\mathbf{Z} 
\to 1.$$ Let $Y_{[4^2, 3, 1]} \to \bar{O}_{[4^2, 3, 1]}$ be a double covering determined 
by the surjection $\pi_1(O_{[4^2, 3, 1]}) \to \mathbf{Z}/2\mathbf{Z}$. 
Then $$Spin (30) \times^Q (\mathfrak{n} + Y_{[4^2, 3, 1]})$$ is a {\bf Q}-factorial 
terminalization of $X$.

\begin{center} 
Research Institute for Mathematical Science, Kyoto University, Japan

namikawa@kurims.kyoto-u.ac.jp 
\end{center}




\end{document}